\documentclass[a4paper,11pt,twoside]{amsart}
\usepackage{cite}
\usepackage{xcolor}
\usepackage[sc]{mathpazo}
\usepackage{eulervm}
\usepackage[utf8]{inputenc}
\usepackage{mathtools}
\usepackage{braket}
\usepackage{amsmath, amssymb, stmaryrd}
\usepackage{enumerate}
\usepackage{float, hyperref, circuitikz}
\hypersetup{hypertexnames=false}
\usepackage{tikz, tikz-cd}
\usetikzlibrary{arrows.meta}
\usetikzlibrary{shapes.geometric, positioning}
\usetikzlibrary{shapes.geometric, positioning}
\usepackage{orcidlink}         

\newcommand{\tip}{\widehat{\pi}}
\newcommand{\spmap}{\nu}
\newcommand{\Spmap}{\widehat{\nu}}
\newcommand{\acyc}{\mathrm{Acyc}}
\newcommand{\Gr}{G}
\newcommand{\rk}{\mathrm{rk}}
\newcommand{\C}{\mathbb{C}}
\newcommand{\Q}{\mathbb{Q}}
\newcommand{\PP}{\mathbb{P}}
\newcommand{\OO}{\mathcal{O}}
\newcommand{\sgn}{\mathbb{S}}
\newcommand{\hbull}{H^{\bullet}}
\newcommand{\Sprod}{\mathcal{S}}
\newcommand{\rag}{\overset{e}\rightarrow}
\newcommand{\tap}{\overset{\Gr}\rightsquigarrow}

\newcommand{\Fiso}{\mathcal{S}}
\newcommand{\ParExtSp}{\mathcal{E}}
\newcommand{\wV}{\overline{V}}
\newcommand{\ipat}{\mathcal{I}}
\newcommand{\lle}{\prec}
\newcommand{\GL}{{\mathrm G\mathrm L}}

\newcommand{\ch}{\operatorname{ch}}

\newcommand{\Ext}{\mathrm{Ext}}
\newcommand{\hilb}{\mathrm{Hilb}}
		
\newcommand{\vect}{M}
\newcommand{\ms}{\mathcal{M}}
\newcommand{\pms}{\mathcal{P}}
\newcommand{\pat}{\mathcal{T}}
\newcommand{\patf}{\mathcal{F}}
\newcommand{\aut}{\mathrm{Aut}}
\newcommand{\ukad}{\underline{\epsilon}_k}
\newcommand{\Aut}{\mathrm{Aut}}

\newcommand{\modtor}{\vect_\rho^{\mathrm{ab}}}
\newcommand{\modtir}{\mathcal{S}_\rho^{\mathrm{ab}}}
\newcommand{\br}{\bar{r}}
\newcommand{\bro}{\overline{\rho}}

\newcommand{\jmap}{\vartheta}
\newcommand{\augf}{\mathbb{F}_\rho}
\newcommand{\bugf}{\mathbb{G}_\rho}

\newcommand{\ula}{\underline{\lambda}}
\newcommand{\la}{\lambda}
\newcommand{\uka}{\underline{k}}
\newcommand{\Idm}{I}

\newcommand{\slice}{\mathcal{C}}

\newcommand{\ls}{\mathcal{L}}
\newcommand{\locmod}{U}
\newcommand{\VV}{\mathcal{V}}
\newcommand{\nparts}{A}

\newcommand{\ccirc}{\mathcal{C}^\circ}
\newcommand{\autbro}{\Sigma[\bro]}

\newcommand\isomto{\stackrel{\textstyle\sim}{\smash{\longrightarrow}\rule{0pt}{0.4ex}}}

\theoremstyle{plain}
\newtheorem{theorem}{Theorem}[section]
\newtheorem{proposition}[theorem]{Proposition}
\newtheorem{lemma}[theorem]{Lemma}

\newtheorem{corollary}[theorem]{Corollary}
\theoremstyle{definition}
\newtheorem{definition}[theorem]{Definition}
\newtheorem{notation}[theorem]{Notation}

\newtheorem{remark}[theorem]{Remark}

\title[The intersection cohomology of moduli
spaces]{Parabolic bundles and the intersection cohomology of moduli spaces of
vector bundles on curves}

\author[Camilla Felisetti]{Camilla Felisetti \orcidlink{0000-0001-9255-7866}
	}
\address{Dipartimento di Scienze Fisiche, Informatiche e Matematiche, 
	Università di Modena e Reggio Emilia}
\email{Camilla.Felisetti@unimore.it}

\author[Andras Szenes]{Andras Szenes \orcidlink{0000-0002-2877-680X}
	}
\address{Section de mathématiques, Université de Genève}
\email{Andras.Szenes@unige.ch}

\author[Olga Trapeznikova]{Olga Trapeznikova \orcidlink{0000-0002-1034-3969}
	}
\address{Faculty of Mathematics, IST Austria}
\email{Olga.Trapeznikova@ist.ac.at}

\subjclass[2020]{Primary 14H60, 14D20; Secondary 55N33, 32S60}
\keywords{Vector bundles, Parabolic bundles, Intersection cohomology}
\begin{document}

\begin{abstract}
The study of the intersection cohomology of moduli spaces of semistable bundles was initiated by Frances Kirwan in the 1980s. In this paper, we develop a self-contained, purely geometric framework to study the intersection cohomology of the moduli spaces of semistable degree-0 vector bundles of arbitrary rank on Riemann surfaces. Motivated by the 2015 work of Mozgovoy and Reineke, our approach applies the Beilinson--Bernstein--Deligne--Gabber Decomposition Theorem to the forgetful map from a parabolic moduli space. We give a detailed description of the topology of this map, the geometry of its fibers, and the precise structure of the relevant local systems, making systematic use of the multiplicative structures of the theory. As a global application of this machinery, we obtain a geometric proof of a recursive formula that reduces the calculation of the intersection Betti numbers of the degree-0 moduli spaces to the known formulas for the smooth, degree-1 moduli spaces.
\end{abstract}
	\maketitle
		\tableofcontents

\section{Introduction: the intersection cohomology of moduli spaces}
Let $C$ be a smooth complex projective curve of genus $g 
\geqslant 2$, and denote by $\ms_d(r)$ the moduli space of semistable, rank-$r$, 
degree-$d$ vector bundles on $C$. When $r$ and $d$ are coprime, $\ms_d(r)$ is a 
smooth projective variety whose cohomology is well understood and has been 
computed using a variety of geometric, gauge-theoretic, and arithmetic 
techniques \cite{Newstead67,HN74,DR75, AB83, EK00, BGL94}. When $\gcd(d, r) > 
1$, the moduli space $\ms_d(r)$ is singular, and the natural replacement for ordinary cohomology is its intersection cohomology $IH^\bullet(\ms_d(r))$ \cite{GM83}. 

The study of the intersection cohomology of the moduli spaces $\ms_d(r)$ was 
initiated by Kirwan in the  1980s \cite{K85,K86,K86vect}. Of particular 
interest is the degree-zero space $\ms_0(r)$, which is homeomorphic to the 
character variety of unitary representations of the fundamental group of $C$
\cite{NS65}. 
For many years, explicit results were largely restricted to  the rank-2 case, until 2015, when Mozgovoy and Reineke \cite{MR15} (see also \cite{Meinhardt})  presented an elegant  recursive formula for the intersection Betti numbers of $\ms_d(r)$ in all  ranks and degrees, relating the intersection cohomology of moduli stacks to  Donaldson--Thomas invariants in the framework of hereditary categories.
 The article \cite{MR15} was finalized at the end of 2024, while
the present paper was in the process of being written up. Recently, some of
these results were extended to the case of $G$-bundles for
an arbitrary reductive group $G$ in \cite{BDIKP}.

	In this paper, we study the Beilinson--Bernstein--Deligne--Gabber Decomposition Theorem \cite{BBD82} for the forgetful morphism $$\pi\colon \mathcal{P}_0(r)\to \ms_0(r),$$ 
where $\pms_0(r)$ is a smooth moduli space of vector bundles on $C$ equipped with a parabolic structure. The Decomposition Theorem expresses the cohomology of $\mathcal{P}_0(r)$ as a direct sum of contributions indexed by the strata of a natural stratification of $\ms_0(r)$.
We show that the open dense stable stratum contributes
$$IH^\bullet(\ms_0(r)) \otimes H^\bullet(\mathbb{P}^{r-1}),$$
while the remaining contributions are expressed in terms of the intersection cohomology of products of lower-rank moduli spaces $\ms_0(r')$  (with $r' < r$), twisted by certain local systems. We explicitly determine the local systems carried by these strata and the corresponding summands in the Decomposition Theorem.

The tautological Hecke correspondence identifies $\mathcal{P}_0(r)$ with a $\mathbb{P}^{r-1}$-bundle over the smooth moduli space $\ms_1(r)$, and therefore gives an isomorphism of graded vector spaces
$$ H^\bullet(\mathcal{P}_0(r))\cong  H^\bullet(\ms_1(r))\otimes  H^\bullet(\PP^{r-1}).$$
Comparing this identity with the decomposition described above and isolating 
the contribution of the open stratum yields a recursive formula for the 
intersection Betti numbers of $\ms_0(r)$ in terms of the ordinary Betti numbers 
of $\ms_1(r)$. The resulting formula differs from that of Mozgovoy and Reineke, 
although it has a similar recursive structure.

Our geometric approach is self-contained, explicit and works entirely with 
projective coarse moduli spaces, avoiding the use of  virtual or 
stack-theoretic invariants.  Beyond the numerical recursion, it unlocks several 
further structural properties of $IH^\bullet(\ms_0(r))$.

In particular, our constructions are naturally compatible with mixed Hodge 
structures, 
allowing the decomposition to be refined to compute the \textit{intersection 
Hodge 
numbers} of $\ms_0(r)$. More generally, the geometric framework developed here 
indicates that our methods are amenable to motivic generalizations, offering a 
way to access finer motivic structures of these spaces.

Also, our explicit description of the Decomposition Theorem for 
$\pi\colon\pms_0(r)\to\ms_0(r)$ identifies $IH^\bullet(\ms_0(r))$ with a 
canonical subspace of $\hbull(\pms_0(r))$, whose ring structure is well-known. 
This could lead to a better understanding of the multiplicative structures on 
$IH^\bullet(\ms_0(r))$. We note that in \cite{FT}, this identification was 
already used to compare the Poincar\'e pairing on 
$\hbull(\pms_0(r))$ to the Poincar\'e--Verdier pairing on 
$IH^\bullet(\ms_0(r))$, and thereby obtain explicit formulas (cf. \cite{JKKW}) 
for the latter.

Finally, our arguments make systematic use of the multiplicative structure on cohomology, leading to several techniques that may be useful in related problems, and  offering a fresh perspective on the Decomposition Theorem itself.

The paper is organized as follows.
\begin{itemize}
	\item In \S\ref{sec:dt-purity-main}, we review the basic properties of 
	intersection cohomology, recall the Decomposition Theorem in a form adapted 
	to our setting, collect the mixed-Hodge-theoretic results used throughout 
	the paper, and state our main results; in \S\ref{sec:plan}, we sketch
	our strategy of the proofs.
 \item In \S\ref{sec:rank2}, we illustrate the main ideas of our approach in the case $r=2$.
	\item In \S\ref{sec:moduliandparmap} we introduce the map
	$\pi:\pms_0(r)\to\ms_{0}(r)$, whose decomposition is the central object of our study.
		\item In \S\ref{sec:parabolicfiber} we calculate the cohomology of the
	relevant fibers of our map $\pi$.
    \item Next, in \S\ref{sec:nosupport}, we prove a key vanishing theorem, 
    which 
  restricts the strata contributing to the Decomposition
    Theorem for \(\pi\).
	\item \S\ref{sec:combsub} and \S\ref{sec:locsys}
	are devoted to the two methods we propose for the identification of  the
	local
	systems over the relevant strata.
	\item Finally, in \S\ref{sec:final}, we derive the recursion mentioned
	above.
\end{itemize}

\textbf{Acknowledgments.} We are deeply indebted to Donu Arapura and Claude
Sabbah for sharing their expertise with us. Useful discussions with Ben
Davison, Mark De Cataldo, Tamas Hausel, Anton Mellit, Luca Migliorini, Sergey
Mozgovoy,  and Vivek Shende are gratefully acknowledged.

C.F. was supported by IN$\delta$AM GNSAGA and by the project: “Discrete Methods
in Combinatorial Geometry and Geometric Topology” of the University of Modena
and Reggio Emilia;
 A. Sz. was supported by SNF grant S110100 and the NCCR SwissMAP;
 O.T. was supported by the European Union's Horizon 2020 research and innovation programme under the Marie Sk\l{}odowska-Curie Grant Agreement No. 101034413.
	\section{The Decomposition Theorem and statement of the main
	results}\label{sec:dt-purity-main}

\subsection{Intersection cohomology and IC-sheaves}
We begin with a quick introduction to intersection
cohomology; for general background, we refer the reader to \cite{dCM09,GM83}.

Some of what follows holds in greater generality, but, for simplicity, we
assume that $Y$ is a possibly singular, complex
quasiprojective algebraic variety.
The \textit{rational intersection cohomology groups} of $Y$ form a graded vector
space $$IH^{\bullet}(Y)=\bigoplus_{i=0}^{2\dim Y}IH^i(Y),$$ and if $Y$ has at 
worst 
finite quotient singularities, then its intersection cohomology agrees with its 
ordinary cohomology $H^\bullet(Y)$. Although intersection cohomology does not 
carry a  product structure in general, it is naturally a module over 
$H^\bullet(Y)$.

When $Y$ is projective with an ample class $\eta\in H^2(Y;\mathbb Q)$, the 
following key properties hold:
\begin{enumerate}[(i)]
    \item (\textit{Hard Lefschetz}) For each $0 \le i \le \dim Y$, multiplying by $\eta^i$ induces an isomorphism:
    $$\cup \eta^i\colon IH^{\dim Y-i}(Y)\isomto IH^{\dim Y+i}(Y).$$
    \item (\textit{Poincar\'e--Verdier duality})
For each $0 \le i \le 2\dim Y$, there is a non-degenerate pairing:
	\begin{equation*}
		IH^i(Y)\otimes IH^{2\dim Y-i}(Y) \to \Q.
	\end{equation*} 
\end{enumerate}

Just as usual cohomology can be described as the cohomology of the constant sheaf 
$\Q_Y$,
$H^{\bullet}(Y)={H}^{\bullet}(Y,\Q_Y),$ there
is a constructible complex of sheaves $IC_Y$ on $Y$, called the
\textit{IC-sheaf}, whose hypercohomology computes the intersection cohomology 
of $Y$:
$$IH^i(Y):= \mathbb{H}^{i}(Y,IC_Y).$$

\textit{Convention}: In this paper, we set the cohomology sheaves of the complex $IC_Y$ to be graded
from 0 to $2\dim Y$, so that if $Y$ is nonsingular, then $IC_Y=\Q_Y$.

As is the case for ordinary cohomology, one may consider 
intersection
cohomology with coefficients in a local system. Given a  Whitney  
stratification by smooth quasi-projective varieties, $Y = \sqcup_\alpha 
Y_\alpha$, one can associate to each local system $\ls$ on a stratum 
$Y_\alpha$  a constructible complex $IC(\overline{Y}_{\alpha}, \mathcal{L})$ on 
$Y$, called the \textit{intermediate extension} of $\mathcal{L}$. This complex 
is supported on the closure $\overline{Y}_{\alpha}$, restricts to $\mathcal{L}$ 
on $Y_\alpha$:
$$IC(\overline{Y}_{\alpha},\ls)|_{Y_{\alpha}}=\ls,$$ and  satisfies standard support and cosupport conditions (see \cite{GM83}).
The corresponding  \textit{twisted intersection cohomology} groups are defined by
$$IH^i(\overline{Y}_{\alpha},\ls):=
\mathbb{H}^{i}(Y,IC(\overline{Y}_\alpha,\ls)).$$
In particular, the IC sheaf is the intermediate extension of the constant sheaf 
on $Y_0\subset Y$, the dense open stratum of smooth points: 
$IC_{Y}=IC(Y,\Q_{Y_0})$.

For cones over projective varieties, the following statement is useful for explicit calculations.

\begin{proposition}[\cite{M23,D95}]\label{prop:ICcones}
	Let $Z$ be a complex projective variety of dimension $n-1$, endowed with an
	ample line bundle $L$.  Let $Y$ be the affine cone
	over $Z$ given by
	contracting the zero-section in the total space of $L^{-1}$. Then
	$$IH^d(Y)\cong \begin{cases} IH_{\mathrm{prim}}^d(Z) &\mbox{ if $d<n$},\\
		0 &\mbox{ otherwise,}\end{cases}$$
	where the primitive cohomology is given by
	$$IH^d_{\mathrm{prim}}(Z)=\operatorname{ker}\left\lbrace \cup
	c_1(L)^{n-d}\colon IH^d(Z)\rightarrow IH^{2n-d}(Z)\right\rbrace.$$
\end{proposition}

\textit{Convention regarding graded local systems.} If
$f:X\rightarrow Y$ is a morphism and $f$ is a fibration over $U\subseteq Y$,
then we place the local system $R^jf_*\Q_{f^{-1}(U)}$ in cohomological degree 
$j$.
More generally, for a graded local system $\ls=\bigoplus_j\ls^j$, with
$\ls^j$ placed in degree $j$, we set
\(IC(Y,\ls):=\bigoplus_j IC(Y,\ls^j).\)

\subsection{The Decomposition Theorem}\label{sec:decomp}\label{sec:strategy}

An effective tool for computing the intersection cohomology of a variety $Y$ is the Decomposition Theorem \cite{BBD82}, applied to a proper morphism
$f\colon X\to Y,$
with $X$ smooth. Below, we present a streamlined version of the theorem adapted 
to our needs. 
For a more comprehensive treatment, we refer the reader to the seminal paper of 
de Cataldo and Migliorini \cite{dCM05} (see also \cite{dCM09} and Williamson's 
Bourbaki exposé \cite{W17}).

\begin{theorem}[Decomposition theorem]\label{thm:DT}
	Let $f\colon X\to Y $ be a proper map from a smooth variety $X$ to a possibly
	singular variety $Y$, and assume that 	$ Y = \sqcup_\alpha Y_\alpha$ is a
	Whitney stratification adapted to $f$; in particular,
	$f^{-1}(Y_\alpha)\to Y_\alpha$ is a topologically locally trivial
	fibration over $Y_\alpha$ for each $\alpha$.
	Let $y\in Y_\beta$, $U\subset Y$ a small open neighborhood of $y$, and
	consider the \emph{cohomological restriction map}
	\begin{equation}\label{eq:cohrestr}
		\jmap_y^{*}\colon H^{\bullet}(f^{-1}(U))\to H^{\bullet}(f^{-1}(U\setminus
		Y_\beta)),\quad\text{where }\jmap_y\colon f^{-1}(U\setminus
		Y_\beta)\hookrightarrow
		f^{-1}(U)
	\end{equation}
	is the inclusion. Set $ L_{y}:= \ker \jmap_y^{*}$, and introduce the local system
	\begin{equation}
		\ls_{\beta}\longrightarrow Y_{\beta}
		\text{ with stalk } L_{y},\, y\in Y_{\beta}.
	\end{equation}
	Then
	\begin{enumerate}[I.]
		\item[0.] $L_{y}$ is a graded
		ideal in $H^{\bullet}(f^{-1}(U))\cong H^{\bullet}(f^{-1}(y))$, whose
		graded pieces are in the degree range $[\dim X-\dim
		Y_\beta-\boldsymbol{r}(f), \dim X-\dim Y_\beta+\boldsymbol{r}(f)]$, where $\boldsymbol{r}(f)=\dim\bigl(X\times_Y X\bigr)-\dim X$ is the
		\textit{defect of semismallness} of $f$.
		\item In the derived category of bounded constructible complexes, there
		is a direct sum decomposition, depending on the choice of an $f$-ample
		class:
		\begin{equation}\label{eq:dt1intro}
			Rf_*\Q_X\cong \bigoplus_\alpha
			IC(\overline{Y}_{\alpha},\ls_{\alpha}).\tag{2.2.A}
		\end{equation}
		\item  In particular, passing to the stalks of the cohomology sheaves 
		of the complexes
		on both
		sides of this decomposition, we obtain the decomposition of the 
		cohomology of the fibers
		\begin{equation}\label{eq:cohsheaves}
			H^{\bullet}(f^{-1}(y)) \cong L_{y}\oplus\bigoplus_{
			\substack{\overline{Y}_\alpha\ni y\\ \alpha\ne\beta} }
			\mathcal{H}^{\bullet}\left(IC(\overline{Y}_{\alpha},\ls_{\alpha})\right)_{y}.
			  \tag{2.2.B}
		\end{equation}
		\item
		Finally, taking global cohomology over $Y$, we obtain the decomposition
		\begin{equation}\label{eq:globalDT}
			H^{\bullet}(X) \cong \bigoplus_\alpha
			IH^{\bullet}(\overline{{Y}}_{\alpha},\ls_{\alpha}).\tag{2.2.C}
		\end{equation}
	\end{enumerate}
\end{theorem}

\vskip 3mm

The strata $Y_\alpha$ appearing in \eqref{eq:dt1intro} for which   $\ls_\alpha$ does not vanish are called
\textit{supports} of $f$.

This theorem suggests the following \textit{strategy} for computing 
$IH^\bullet(Y)$. 
Assume that the local system $\ls_0$ on the open dense stratum $Y_0$   is trivial, and may be 
identified with the cohomology of the 
generic fiber $F_0$ of $f$. Then the contribution of $Y_0$ to the 
right-hand side of \eqref{eq:globalDT} is
\begin{equation}\label{eq:genfiber}
	IH^\bullet(Y,\ls_0)
	\cong IH^{\bullet}(Y)\otimes
	H^{\bullet}(F_0),\end{equation}
and, if $H^{\bullet}(X)$ is known, the equality \eqref{eq:globalDT} may be
	used to start a recursive procedure to calculate $IH^\bullet(Y)$.

We will describe this strategy in more detail at the end  of 
this Section, in \S\ref{sec:plan}, after presenting our results in 
\S\ref{sec:summary}.
First, however, we make a brief detour to record a
technical refinement of the Decomposition Theorem for later reference.  

\subsection{Mixed Hodge structures and the Decomposition 
Theorem}\label{sec:hodge}
 We begin by recalling the necessary background on mixed Hodge
structures (see \cite{D71,PS08}).

\subsubsection{Purity and weight polynomials}

\label{sec:wp}

The cohomology groups of a complex algebraic variety $Y$ are endowed with an 
increasing \textit{weight filtration} $W_\bullet$ such that
\[
\mathrm{Gr}^W_k H^{m}(Y):=W_kH^{m}(Y)/W_{k-1}H^{m}(Y)=0,\quad \text{ for 
}k<0,k>2m.
\]
We say that the cohomology of $Y$ is \emph{pure} if
\[
\mathrm{Gr}^W_k H^m(Y)=0 \qquad\text{whenever }k\neq m.
\]
Purity may be characterized in terms of resolutions of singularities.
\begin{proposition}\label{prop:critpure}
	Let $Y$ be a projective algebraic variety. Then $Y$ has pure cohomology if 
	and 
	only if, for some (and hence for any) birational morphism 
	$\widetilde{Y}\to Y$ from a smooth projective variety $\widetilde{Y}$, the 
	induced pullback map
	$\hbull(Y)\to\hbull(\widetilde{Y})$ is injective.
\end{proposition}
\begin{corollary}\label{cor:pureinject}
	If $X\to Y$ is a proper birational morphism and the projective variety
	$Y$ has pure cohomology, then the pullback
	$\hbull(Y)\to\hbull(X)$ is injective.
\end{corollary}

Proposition \ref{prop:purecones} below describes an important class of 
varieties with pure cohomology. 
\begin{definition} \label{def:semiproj}
	A \textit{semiprojective variety} is a quasi-projective variety $X$ 
	equipped with a $\C^*$-action such that the fixed-point locus is proper 
	and, for every $x\in X$, the limit $\lim_{t\to 0}t\cdot x$ exists in $X$.
\end{definition}
\begin{proposition}\label{prop:purecones}
	Let $X$ be a semiprojective variety. Then the intersection cohomology
	$I\hbull(X)$ has pure Hodge structure.
\end{proposition}

Similarly, the compactly supported cohomology $H_c^\bullet(Y)$ carries a natural mixed Hodge structure and weight filtration, with the same conventions and notion of purity as above.
We define the \emph{weight polynomial} of $Y$ by
\[
E_t(Y)=\sum_{j,k}(-1)^{j}\dim \mathrm{Gr}^W_k(H_c^j(Y)) t^k.
\]
This invariant has the following important motivic property.
\begin{proposition}\label{prop:motivicwp}
	Let $Z\subset Y$ be a locally closed subvariety.
	Then
	\[
	E_t(Y)=E_t(Z)+E_t(Y\setminus Z).
	\]
\end{proposition}

\subsubsection{Purity and the Decomposition Theorem}
\label{sec:mhsdecomp}

With these preparations, we are ready to revisit the Decomposition Theorem with
the goal of providing a criterion for the purity of the cohomology of the fibers
\(f^{-1}(y)\) appearing in \eqref{eq:cohsheaves}. The  statement below was
communicated to us by Claude Sabbah.

We call a \textit{graded local system pure} if the weight of each homogeneous 
component
coincides with its degree. Under the total grading, tensor
products, direct sums, and direct summands of pure graded objects remain pure.

\begin{proposition}\label{prop:purefibers}
	Assume the setup of Theorem \ref{thm:DT}, and fix \(y\in Y\).
	Suppose that, for every stratum \(Y_\alpha\) such that
	\(
	y\in\overline Y_\alpha\setminus Y_\alpha,
	\)
	there exists a finite algebraic morphism
	\[
	\gamma_\alpha\colon
	Z_\alpha\longrightarrow\overline Y_\alpha
	\]
	with the following properties:
	
	\begin{enumerate}[(i)]
		\item the open subset
		\(
		\widetilde Z_\alpha:=\gamma_\alpha^{-1}(Y_\alpha)
		\)
		is dense in \(Z_\alpha\), and the restriction
		\[
		\widetilde\gamma_\alpha
		:=
		\gamma_\alpha|_{\widetilde Z_\alpha}
		\colon
		\widetilde Z_\alpha\longrightarrow Y_\alpha
		\]
		is finite étale;
		
		\item the pulled-back local system
		\(
		\widetilde\gamma_\alpha^*
		\mathcal L_\alpha^\bullet
		\)
		is constant on each connected component of
		\(\widetilde Z_\alpha\);
		
		\item for every \(z\in\gamma_\alpha^{-1}(y)\), the local intersection
		cohomology of \(Z_\alpha\) at \(z\) is pure.
	\end{enumerate}
	
	Then the mixed Hodge structure on
	\(
	H^\bullet(f^{-1}(y),\mathbb Q)
	\)
	is pure.
\end{proposition}

\begin{proof}
	By Saito's theory \cite{S89}, \(\mathcal L_\alpha^\bullet\) underlies a
	polarizable variation of pure Hodge structures on \(Y_\alpha\), and hence 
	the same is true for the pulled back local system
	\(\widetilde\gamma_\alpha^*\mathcal L_\alpha^\bullet\).
	Combined with assumption \textup{(ii)}, this means that 
	\(\widetilde\gamma_\alpha^*\mathcal L_\alpha^\bullet\) is a constant local 
	system underlying a polarizable variation of pure Hodge structures on
	the algebraic variety \(\widetilde Z_\alpha\). According to the
	Griffiths--Schmid fixed-part theorem \cite[\S6(a)]{GS75}, this implies
	that  on each connected
	component of \(\widetilde Z_\alpha\),
	\[
	\widetilde\gamma_\alpha^*\mathcal L_\alpha^\bullet
	\simeq
	L_\alpha^\bullet\otimes\mathbb Q_{\widetilde Z_\alpha}
	\]
	as graded variations of Hodge structure, where
	\(L_\alpha^\bullet\) is pure. Hence, componentwise, the intermediate 
	extension of $\widetilde\gamma_\alpha^*\mathcal L_\alpha^\bullet$ is the 
	tensor product of a pure vector space with the IC-sheaf of $Z_{\alpha}$:
	\[
	IC\left(
	Z_\alpha,
	\widetilde\gamma_\alpha^*\mathcal L_\alpha^\bullet
	\right)
	\cong
	IC(Z_\alpha,\mathbb Q)\otimes L_\alpha^\bullet.
	\]
	Taking stalks, we can conclude that for every \(z\in\gamma_\alpha^{-1}(y)\),
	\[
	\mathcal H^\bullet
	\left(
	IC\left(
	Z_\alpha,
	\widetilde\gamma_\alpha^*\mathcal L_\alpha^\bullet
	\right)
	\right)_z
	\simeq
	\mathcal H^\bullet
	\bigl(IC(Z_\alpha,\mathbb Q)\bigr)_z
	\otimes L_\alpha^\bullet.
	\]
	By assumption \textup{(iii)}, the first factor on the right is
	pure, and since \(L_\alpha^\bullet\) is pure, the stalk on the left is
	pure as well.
	
	Pushing forward along the finite map
	\(\gamma_\alpha\), and taking the stalk at \(y\),  we obtain
	\[
	\mathcal H^\bullet
	\left(
	R\gamma_{\alpha *}
	IC\left(
	Z_\alpha,
	\widetilde\gamma_\alpha^*\mathcal L_\alpha^\bullet
	\right)
	\right)_y
	\simeq
	\bigoplus_{z\in\gamma_\alpha^{-1}(y)}
	\mathcal H^\bullet
	\left(
	IC\left(
	Z_\alpha,
	\widetilde\gamma_\alpha^*\mathcal L_\alpha^\bullet
	\right)
	\right)_z,
	\]
	from which we can conclude that the stalk on the left hand side is pure.
	
	Since
	\(\widetilde\gamma_\alpha\colon\widetilde Z_\alpha\to Y_\alpha\)
	is finite étale, the adjunction morphism
	\(
	\mathcal L_\alpha^\bullet
	\longrightarrow
	\widetilde\gamma_{\alpha *}
	\widetilde\gamma_\alpha^*
	\mathcal L_\alpha^\bullet
	\)
	is split by the normalized trace map, and thus
	\(\mathcal L_\alpha^\bullet\) is a direct summand of
	\(
	\widetilde\gamma_{\alpha *}
	\widetilde\gamma_\alpha^*
	\mathcal L_\alpha^\bullet
	\).
	Because finite pushforward commutes with intermediate extension,
	it follows that
	\(
	IC(\overline Y_\alpha,\mathcal L_\alpha^\bullet)\)
	 is a direct summand of 
	\(R\gamma_{\alpha *}
	IC\left(
	Z_\alpha,
	\widetilde\gamma_\alpha^*\mathcal L_\alpha^\bullet
	\right)
	\), and hence
	\[
	\mathcal H^\bullet
	\left(
	IC(\overline Y_\alpha,\mathcal L_\alpha^\bullet)
	\right)_y
	\text{ is a direct summand of }
	\mathcal H^\bullet
	\left(
	R\gamma_{\alpha *}
	IC\left(
	Z_\alpha,
	\widetilde\gamma_\alpha^*\mathcal L_\alpha^\bullet
	\right)
	\right)_y.
	\]
	We can thus conclude that
	\(
	\mathcal H^\bullet
	\left(
	IC(\overline Y_\alpha,\mathcal L_\alpha^\bullet)
	\right)_y
	\)
	is pure, being a direct summand of a pure graded mixed Hodge
	structure.
	
	Every boundary-stratum summand in \eqref{eq:cohsheaves} is therefore
	pure. Since \(L_y\) is pure by Saito's theory, it follows that
	\(H^\bullet(f^{-1}(y),\mathbb Q)\) is pure as well.
\end{proof}

\subsection{Summary of the main results}\label{sec:summary}
Our central result is a detailed description of the cohomology 
of the fibers of the map $\pi$, as well as the  local systems guaranteed 
by the Decomposition Theorem for $\pi$.

For each partition $\rho=[r_1,\ldots,r_k]\vdash r,$
we define the \textit{abelian stratum}
$$\modtor := \{ V=V_1\oplus\cdots\oplus V_k\in \ms_0(r): \;\;  V_i\in
\ms_0(r_i)^{\mathrm{st}}, \; V_i\not\cong V_j\text{ for }i\neq j \}\subset
\ms_0(r),$$
where $\ms_0(r_i)^{\mathrm{st}}$ denotes the locus of stable bundles in
$\ms_0(r_i)$.
We can construct this stratum as follows: consider the projective variety
\[
\Sprod_\rho:=\ms_0(r_1)\times\cdots\times\ms_0(r_k),
\]
which carries the action of the automorphism group of the partition
\[
\aut(\rho):=\{\sigma\in\Sigma_{k}|\, r_i=r_{\sigma(i)},\,i=1,\dots,k\}.
\]
 Then the  direct-sum morphism
\begin{equation}\label{eq:gammarho}
	\gamma_\rho\colon \Sprod_\rho\longrightarrow\overline{\modtor},
	\qquad (V_1,\ldots,V_k)\longmapsto V_1\oplus\cdots\oplus V_k,
\end{equation}
is finite and $\aut(\rho)$-invariant. Over the open stratum, it restricts to a
finite \'etale cover
$\widetilde\gamma_\rho\colon\modtir\to\modtor$ with deck group $\aut(\rho)$, 
where
\begin{equation}\label{eq:defmodtir}
	\modtir:=\gamma_\rho^{-1}(\modtor)=
	\{ (V_1,\ldots,V_k)\in \prod_{i=1}^k\ms_0(r_i)^{\mathrm{st}}:
	\;\;  V_i\not\cong V_j\text{ for }i\neq j \};
\end{equation}
in particular $\modtor\cong\modtir/\aut(\rho)$.

\begin{theorem}\label{thm:mainDTdescription}
Fix $p\in C$, and let $\mathcal{P}_0(r)$ be the moduli space
defined in \S\ref{sec:parmoduli}; explicitly:
\[
\mathcal{P}_0(r):=
\left\{
(W,\ell)\ \middle|\
\begin{array}{l}
W\text{ is semistable of rank }r\text{ and degree }0,\\
\ell\subset W_p\text{ is a line, and}\\
W'\subset W,\ \ell\subset W'_p\Longrightarrow\deg W'<0
\end{array}
\right\}.
\]
 Denote by $\pi\colon \mathcal{P}_0(r)\to \ms_0(r)$ the morphism that forgets 
 the line $\ell$. Then:
\begin{enumerate}[(i)]
\item The complement
$$\ms_0(r)\setminus\bigsqcup_{\rho\vdash r}\modtor $$
contains no support for $\pi$, and thus, for \(V\in\modtor\),
\eqref{eq:cohsheaves} takes the form 
\begin{equation}\label{eq:dtsheaves} 
	H^\bullet (\pi^{-1}(V))\cong \mathcal{L}_\rho|_V \oplus \bigoplus_{\substack{ \mu\vdash r,\ 
	\mu>\rho
}} 
	\mathcal{H}^\bullet (
	IC\bigl(\overline{\vect^\mathrm{ab}_{\mu}},
	\mathcal{L}_\mu))_V,
\end{equation}
where, for each partition $\mu\vdash r$, $\mathcal L_\mu$ is 
a graded local system on $M_\mu^{\mathrm{ab}}$.
\item For $V=V_1\oplus\cdots\oplus V_k\in\modtor$, the fiber 
\(F_V=\pi^{-1}(V)\) is a union of iterated projective bundles. Its cohomology 
ring is generated by degree-two classes
$x_1,\ldots,x_k\in H^2(F_V)$ naturally associated with the summands 
\(V_1,\ldots,V_k\), and admits a presentation 
\begin{equation}\label{eq:intermediatefiberpres}
H^\bullet (F_V) \cong \mathbb{Q}[x_1,\ldots,x_k]/\Idm,
\end{equation}
where \(\Idm\) is the ideal generated by the relations in 
\eqref{eq:HStrelationsFib}. As \(V\) varies in \(\modtor\), these cohomology 
groups form a graded local 
system, which is associated to the cover $\modtir\to\modtor$ and 
the natural action of $\aut(\rho)$ on \(\mathbb{Q}[x_1,\ldots,x_k]/\Idm\).

\item  Under the presentation \eqref{eq:intermediatefiberpres}, the stalk
\(\mathcal{L}_\rho|_V\subset \mathbb{Q}[x_1,\ldots,x_k]/\Idm\) is the principal 
ideal generated by the element 
\[ w_\rho = \prod_{1\leq i<j\leq k} (x_i-x_j)^{r_ir_j(g-1)}. \]
In addition, the following set forms a basis of \(\mathcal{L}_\rho|_V\):
\[
\left\{
w_\rho \prod_{i=1}^k x_i^{m_i}
\;\middle|\;
0\leq m_i<r_i,\quad i=1,\ldots,k
\right\},
\]
and, in particular,
\(
\operatorname{rank}(\mathcal{L}_\rho)
=
r_1\cdot r_2\cdots r_k.
\)
\end{enumerate}
\end{theorem}
Parts (i), (ii), and (iii) are proved, respectively, in Theorem
\ref{thm:nosupp}, Theorem \ref{thm:fiberring}, and Theorem \ref{thm:corL}.

As mentioned in the introduction, the tautological Hecke correspondence 
identifies the moduli space \(\mathcal{P}_0(r)\) with a 
\(\mathbb{P}^{r-1}\)-bundle over \(\mathcal{M}_1(r)\) (see Corollary 
\ref{cor:hecke}). This allows us to derive an elegant recursive formula for the 
intersection Betti numbers of $\ms_0(r)$. 

\subsection{The plethystic exponential and the recursion formula}
Before stating the formula, we introduce some standard notation.
Let $\mathcal{H}^{\bullet}=\oplus_i \mathcal{H}^i$ be a graded vector space  
with finite-dimensional graded components;  we define its Hilbert series as the 
formal sum
$$\operatorname{Hilb}_t(\mathcal{H}^{\bullet}) := \sum_i t^i\dim \mathcal{H}^i.$$
For a variety $X$, we introduce the two Poincaré polynomials:
\[
P_t(X)=\mathrm{Hilb}_t(H^\bullet(X))\quad\text{and}\quad
IP_t(X)=\mathrm{Hilb}_t(IH^\bullet(X)).
\]
Similarly, for a bigraded vector space $\mathcal{H}^{\bullet} = 
\oplus_{i,j}\mathcal{H}^{i,j}$, we set
$$\mathrm{Hilb}_{t,q}(\mathcal{H}^{\bullet})=\sum_{i,j}t^iq^j\dim\mathcal{H}^{i,j}.$$
Now, given a Laurent series $Q(t,q)=\sum_{i,j}a_{i,j}t^iq^j$,
whose support is contained in an open half-plane, the
\textit{plethystic exponential} of $Q$ is defined as the function
\begin{equation}\label{eq:defplexp0}
	\mathbb{E}\mathrm{xp}[Q]=\prod_{i,j}(1-t^iq^j)^{-a_{i,j}}
\end{equation}
expanded in the direction of this half-plane.
 This notation is standard in motivic Donaldson--Thomas theory (see, for 
 example, \cite{MR15}).

The plethystic exponential has the following useful interpretation. For a 
bigraded vector space $\mathcal H^\bullet$, introduce the graded polynomial ring
\begin{equation}\label{eq:qtring}
	\mathbb{Q}_{[t]}[\mathcal{H}^{\bullet}]=
	\mathrm{Sym}^{\bullet}\left(\mathcal{H}^{\mathrm{even}}\right)\otimes
	\wedge^{\bullet}\left(\mathcal{H}^{\mathrm{odd}}\right)\!,
\end{equation}
where parity is determined by the $t$-degree.  With the usual Koszul sign
convention, one has then
\[
\mathbb{E}\mathrm{xp}\left[\mathrm{Hilb}_{-t,q}(\mathcal{H}^{\bullet})\right]=\mathrm{Hilb}_{-t,q}(\mathbb{Q}_{[t]}[\mathcal{H}^{\bullet}]).
\]

With this notation established, we can state our recursion:
\begin{theorem}\label{thm:main}
The Poincar\'e polynomials of the  moduli spaces $\ms_1(r)$ and  the intersection Poincar\'e polynomials of the  moduli spaces $\ms_0(r)$ are related by the following
identity:
	\begin{multline}\label{eq:recursion}
		1+\sum_{r\geq1}P_{-t}(\PP^{r-1})\cdot P_{-t}(\ms_1(r))\cdot
		(-t)^{(1-g)r^2}q^r=\\
		\mathbb{E}\mathrm{xp}\left[\sum_{r\geq 1}P_{-t}(\PP^{r-1})\cdot
		IP_{-t}(\ms_0(r))\cdot
		(-t)^{(1-g)r^2}q^r\right]\!.
\end{multline}
\end{theorem}
\begin{remark} 
	Once the Poincar\'e polynomials $P_{-t}(\ms_1(r))$ are known, formula 
	\eqref{eq:recursion} determines $IP_{-t}(\ms_0(r))$ recursively in the rank $r$.
\end{remark}

\subsection{Roadmap for the proofs}\label{sec:plan}
In \S\ref{sec:strat}, we introduce the stratification of $\ms_0(r)$ by strata ${\vect}_{\bro}$ indexed by multipartitions $\bro$ of $r$, which is adapted to $\pi$. The Decomposition Theorem for $\pi$ therefore gives
\begin{equation}\label{eq:decomp_nonab}
R\pi_*\Q_{\pms_0(r)}
\cong
\bigoplus_{\bro\Vdash r}
IC(\overline{\vect}_{\bro},\ls_{\bro}).
\end{equation}

In \S\ref{sec:nosupport}, we prove a \textit{no-support theorem}  showing that 
all summands indexed by non-abelian multipartitions vanish. Consequently, 
\eqref{eq:decomp_nonab} reduces to the abelian strata:
\begin{equation}\label{eq:dtforrho}
R\pi_*\Q_{\pms_0(r)} \cong \bigoplus_{\rho\vdash r} IC(\overline{\modtor},\ls_\rho).
\end{equation}

For any point \(V\in\modtor\), let \(F_V=\pi^{-1}(V)\) denote the fiber over $V$, and set \(L_V=(\ls_\rho)_V\subset H^\bullet(F_V)\). Taking the stalk at \(V\) in \eqref{eq:dtforrho} gives the decomposition
\begin{equation}\label{eq:cohsheavesrho}
H^\bullet(F_V)
\cong
L_V\oplus
\bigoplus_{\mu>\rho}
\mathcal H^\bullet
\left(IC(\overline{\vect^{\mathrm{ab}}_\mu},\ls_\mu)\right)_V,
\end{equation}
where the direct sum ranges over all partitions $\mu\vdash r$ strictly coarser 
than $\rho$. This proves Theorem \ref{thm:mainDTdescription}(i).

The proof of the cohomology-ring presentation in Theorem 
\ref{thm:mainDTdescription}\textup{(ii)} relies on the analysis of the geometry 
of the fibers over the abelian strata carried out in 
\S\ref{sec:parabolicfiber}, where they are described as unions of iterated 
projective bundles, their cohomology is shown to be pure, and the global 
degree-two classes $x_i$ are constructed. These results are used in 
\S\ref{sec:fiberring} to obtain the required presentation of $H^\bullet(F_V)$.

To establish part (iii), we use two complementary
methods: the subtraction argument of \S\ref{sec:combsub} determines the Hilbert
series of $L_V$, while the refined intersection-form calculation of
\S\ref{sec:locsys} provides the additional geometric information needed to identify
$L_V$ as an ideal of $H^\bullet(F_V)$. 
This allows us to calculate the monodromy of $\ls_{\rho}$, and, at this point, we are 
ready for
the final step: in \S\ref{sec:final}, we apply the global Decomposition Theorem
to obtain the recursion of Theorem \ref{thm:main}.

\section{The case of rank 2}\label{sec:rank2}
While not all of the difficulties present in higher rank arise when $r=2$, this 
case serves as a useful model for our method. The resulting formula 
\eqref{eq:mainrank2} was first obtained by Kirwan \cite{K86vect} (see also  \cite{K06}).

\subsection{The formula in rank 2}

Comparing the coefficients of $q^{2}$ in \eqref{eq:recursion}, we obtain a formula 
of 
the form
\[ 
(1+t^2) P_{-t}(\ms_{1}(2)) = (1+t^2) IP_{-t}(\ms_{0}(2)) + 
t^{2(g-1)}\cdot[\text{contribution of }r=1].
\]
As we see from \eqref{eq:defplexp0}, this latter contribution depends crucially on 
the 
sign of $t$ on the right-hand side of \eqref{eq:recursion}. This sign is $+1$ if 
$g$ is odd, and $-1$ if $g$ is even, and thus, setting $H=H^\bullet(\ms_0(1))$, 
the contribution of $r=1$ is 
\[ 
\begin{cases}
	\mathrm{Hilb}_{t}(\mathrm{Sym}^2H^{\mathrm{even}}\oplus 
	\bigwedge^2H^{\mathrm{odd}}\oplus (H^{\mathrm{even}}\otimes
	H^{\mathrm{odd}})) &\text{if }g\text{ is odd},\\
	\mathrm{Hilb}_{t}(\bigwedge^2H^{\mathrm{even}}\oplus 
	\mathrm{Sym}^2H^{\mathrm{odd}}\oplus (H^{\mathrm{even}}\otimes
	H^{\mathrm{odd}})) &\text{if }g\text{ is even}.
\end{cases}
\]

As $\ms_0(1)=\operatorname{Jac}^0(C)$ is the Jacobian 
variety of $C$,  we have $P_{-t}(\ms_{0}(1))=(1-t)^{2g}$, and, using basic  character theory, we arrive at
\begin{multline}\label{eq:mainrank2}
	(1+t^2) P_{-t}(\ms_{1}(2)) = (1+t^2) IP_{-t}(\ms_{0}(2)) \\
	+\frac12t^{2(g-1)}\left[ 
	(1-t)^{4g} +(-1)^{g-1}(1-t^2)^{2g}\right].
\end{multline}

\subsection{The stratification of $\ms_0(2)$ and the normal slices}
The moduli space $\ms_0(2)$ 
is stratified by three smooth subvarieties (cf. Theorem \ref{thm:moduliM}):
\begin{enumerate}
	\item $\vect_{[2]}\cong \ms_0(2)^\mathrm{st}$, the set of stable bundles $V$ of
	$\dim	= 4g-3$,
	\item $\vect_{[1,1]}^{\mathrm{ab}}\cong (\operatorname{Jac}^0(C)^2\!\setminus\!\Delta)/\Sigma_2$,  the set $\{L_1\oplus L_2|\, 
	L_1\not\cong 
	L_2\in\operatorname{Jac}^0(C)\}$ of $\dim= 2g$,
	\item $\vect_{[1,1]}^{\mathrm{nab}}\cong \operatorname{Jac}^0(C)$, the locus $\{L\oplus L|\, 
	L\in\operatorname{Jac}^0(C)\}$ of  $\dim =g$.
\end{enumerate}

The set of isomorphism classes of  semistable bundles which are 
$S$-equivalent to a point in each 
stratum looks as follows:
\begin{enumerate}
	\item just one element, the bundle $V$ itself;
 \item the set of non-trivial extensions parameterized by the points of $\PP\Ext^{1}(L_1,L_2)$ and $\PP\Ext^{1}(L_2,L_1)$, which are two copies of $\PP^{g-2}$, and the single element, the polystable bundle $L_1\oplus L_2$;
	\item the set of non-trivial extensions parameterized by 
	$\PP\Ext^{1}(L,L)\cong\PP^{g-1}$, and the polystable bundle $L\oplus L$. 
\end{enumerate}

\subsection{The fibers of the map $\pi$}\label{sec:fibers2}
\begin{figure}[H]
	\centering
	\begin{tikzpicture}[scale=1, transform shape] 
		\draw[thick] (0, 0) ellipse (3 and 0.7); 
		\node[scale=0.7] at (-1.45, 0.2) {$\vect_{[2]}$}; 
		\draw [thick, orange, rounded corners=7] (-2, 0) --(-1,-0.3)-- (0,0)--  
		(1,-0.3) -- (2, 0); 
		\node[scale=0.7, orange] at (2.4, 0) {$\vect^{\mathrm{ab}}_{[1,1]}$}; 
		\node[red] at (0, -0.04) { \textbullet }; 
		\node[scale=0.7, red] at (0, -0.3) {$\vect^{\mathrm{nab}}_{[1,1]}$}; 
		\node[draw, circle, minimum size=1cm, fill=gray!15] at (-1.85, 2.2) {}; 
		\node[scale=0.6] at (-1.85, 2.2) {\textbf{$\mathbb{P}^1$}}; 
		\node[draw, ellipse, minimum width=0.7cm, minimum height=1.7cm, 
		fill=orange!30] at (1.3, 2.2) {}; 
		\node[draw, ellipse, minimum width=0.7cm, minimum height=1.7cm, 
		fill=orange!30] at (2, 2.2) {};  
		\node[scale=0.6] at (1.65, 2.2) {\textbf{$\mathbb{P}^{g-1} \lor 
				\mathbb{P}^{g-1}$}}; 
		\node[draw, ellipse, minimum width=0.7cm, minimum height=1.7cm, 
		fill=red!30] at (0, 2.2) {}; 
		\node[scale=0.6] at (0, 2.2) {\textbf{$\mathbb{P}^{g-1}$}};  
		\draw[->, dashed] (-1.85, 1.6) -- (-1.85, 0.2); 
		\draw[->, dashed] (0, 1.25) -- (0, 0.06); 
		\draw[->, dashed] (1.65, 1.45) -- (1.65, -0.09); 
	\end{tikzpicture}
	\vskip 14pt
	\setlength{\belowcaptionskip}{-8pt}\caption{The fibers of 
		$\pi:\mathcal{P}_0(2)\to \mathcal{M}_0(2)$.}\label{fig:rank2}
\end{figure}
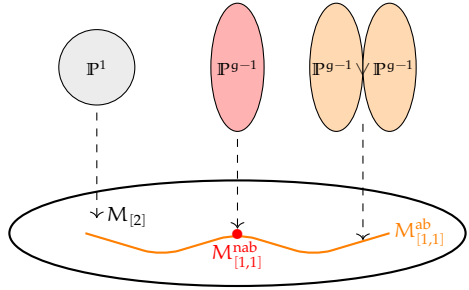

Fix $p\in C$. The smooth projective variety $\pms_0(2)$ parametrizes pairs $(W,\ell)$, where $W$ is a semistable vector bundle of rank $2$ and degree $0$, and $\ell\subset W_p$ is a line such that $\ell\neq L_p$ for every degree-$0$ line subbundle $L\subset W$ (see \S\ref{sec:parmoduli}). We analyze the Decomposition Theorem for the forgetful map
\[
\pi:\pms_0(2)\longrightarrow\ms_0(2),
\qquad
(W,\ell)\longmapsto[W].
\]
The fibers of $\pi$ are described below (see Figure \ref{fig:rank2}).
\begin{enumerate}
	\item For $V$ stable, $\pi^{-1}(V)\cong\PP(V_p)\cong\PP^1$.
	\item For $L_1\not\cong L_2$, the fiber $\pi^{-1}(L_1\oplus L_2)$ consists 
	of 
	\begin{itemize}
		\item a line bundle over $\PP\Ext^{1}(L_1,L_2)$, corresponding to the 
		extra data of $\ell\neq(L_2)_p$,
		\item a line bundle over $\PP\Ext^{1}(L_2,L_1)$,
		\item and a single point corresponding to  $L_1\oplus L_2$.
	\end{itemize}
Hence  $\pi^{-1}(L_1\oplus L_2) \cong \PP^{g-1}\lor_{\mathrm{pt}} \PP^{g-1}$
with cohomology ring
\begin{equation}\label{eq:cohringrank2}
	H^{\bullet}(\pi^{-1}(L_1\oplus L_2))
	\cong
	\Q[x_1,x_2]/\langle x_1x_2,x_1^g,x_2^g\rangle.
\end{equation}
	\item Finally, $\pi^{-1}(L\oplus L) \cong \PP\Ext^{1}(L,L)\cong\PP^{g-1}$.
\end{enumerate}

In \S\ref{sec:parabolicfiber}, we give a more systematic 
description of the fibers in arbitrary rank.

\subsection{Application of the Decomposition Theorem and subtraction}
The local and global forms \eqref{eq:cohsheaves} and
\eqref{eq:globalDT} of the Decomposition Theorem give rise to a subtraction
procedure, which we now apply to
$\pi\colon\pms_0(2)\to\ms_0(2)$, and recover
\eqref{eq:mainrank2}.

First,  consider the left-hand side of \eqref{eq:globalDT}.
By the Hecke correspondence (see Corollary \ref{cor:hecke}), $\pms_0(2)$ is identified with a  $\PP^1$-bundle over $\ms_1(2)$, yielding
\[
P_{-t}(\pms_0(2))
=P_{-t}(\ms_1(2))\cdot P_{-t}(\PP^1),
\]
which matches the left-hand side of \eqref{eq:mainrank2}.

We now turn to the right-hand side of \eqref{eq:globalDT}.
Over the  \textit{stratum $\vect_{[2]}$ of smooth points}, the local system 
$\mathcal{L}_{[2]}$ is trivial with fiber 
$H^\bullet(\mathbb{P}^1)=\mathbb{Q}\oplus\mathbb{Q}[-2]$. Its contribution to 
\eqref{eq:globalDT} is therefore $$IH^\bullet(\ms_0(2))\otimes 
H^\bullet(\mathbb{P}^1),$$ with Hilbert series $(1+t^2)\,IP_{-t}(\ms_0(2)).$
This recovers the first term on the right-hand side of \eqref{eq:mainrank2}.

To calculate the contribution of the \textit{abelian stratum} $\vect^{\mathrm{ab}}_{[1,1]}$, we first  identify the local system $\ls_{[1,1]}$. For a point $L_1\oplus L_2$ in this stratum, decomposition \eqref{eq:cohsheaves} yields
 \begin{equation}\label{eq:eqbrank2}
	H^\bullet(\pi^{-1}(L_1\oplus L_2))
	\cong
	(\ls_{[1,1]})_{L_1\oplus L_2}
	\oplus
	(\Q\oplus\Q[-2])\otimes IH^\bullet(\slice_{L_1\oplus L_2}),
\end{equation}
where $\slice_{L_1\oplus L_2}$ is the normal slice to $\vect^{\mathrm{ab}}_{[1,1]}$ at $L_1\oplus L_2$. This slice is the cone over
\[ 
\hat{\slice}_{L_1\oplus L_2}\cong \PP\Ext^{1}(L_2,L_1)\times 
\PP\Ext^{1}(L_1,L_2),
\]
obtained by contracting the zero-section of the line bundle  $\OO(-1)\boxtimes  
\OO(-1)$ to a point. 

On the left-hand side of \eqref{eq:eqbrank2}, we have the fiber cohomology ring \eqref{eq:cohringrank2}, equipped with the action of the monodromy group $\Sigma_2$ that interchanges $x_1$ and $x_2$. To determine $(\ls_{[1,1]})_{L_1\oplus L_2}$, it remains to compute $IH^\bullet(\slice_{L_1\oplus L_2})$. By Proposition \ref{prop:ICcones}, this intersection cohomology coincides with the primitive cohomology of $\hat{\slice}_{L_1\oplus L_2}$ with respect to multiplication by $c_1(\mathcal{O}(1)\boxtimes \mathcal{O}(1))$.

Explicitly, let 
 $x_1=c_1\bigl(\mathcal{O}_{\mathbb{P}\mathrm{Ext}^{1}(L_2,L_1)}(1)\bigr)$ and 
 $x_2=c_1\bigl(\mathcal{O}_{\mathbb{P}\mathrm{Ext}^{1}(L_1,L_2)}(1)\bigr)$. The 
 nontrivial element of $\Sigma_2$ interchanges $L_1$ and $L_2$, thus  swapping 
 $x_1$ and $x_2$. We have
\[
H_{\mathrm{prim}}^{\bullet}\bigl(\mathbb{P}\mathrm{Ext}^{1}(L_2,L_1)\times\mathbb{P}\mathrm{Ext}^{1}(L_1,L_2)\bigr)
\cong \mathbb{Q}[x_1,x_2]/\langle x_1^{g-1},x_2^{g-1},x_1+x_2 \rangle, 
\]
which gives the following identification of graded $\Sigma_2$-modules:
\begin{equation}\label{eq:simpleprim}
	IH^{\bullet}(\slice_{L_1\oplus L_2})
	\cong
	\Q[x_1-x_2]/\langle(x_1-x_2)^{g-1}\rangle.
\end{equation}

Now, using \eqref{eq:cohringrank2}, \eqref{eq:eqbrank2} and \eqref{eq:simpleprim},  we can identify the stalk of $\ls_{[1,1]}$ by ``subtracting'' the graded $\Sigma_2$-module
\[
(\Q\oplus\Q[-2])\otimes
\Q[x_1-x_2]/\langle(x_1-x_2)^{g-1}\rangle\quad\text{from}\quad
\Q[x_1,x_2]/\langle x_1x_2,x_1^g,x_2^g\rangle.
\]
The difference is the one-dimensional graded module
$\Q\cdot(x_1-x_2)^{g-1}$. Consequently, the local system is
\[\ls_{[1,1]} = (\operatorname{Jac}^0(C)^2\setminus\Delta) \times_{\Sigma_2} \left[\sgn[-2]\right]^{\otimes(g-1)}, \]
where $\sgn[-2]$ denotes the one-dimensional sign representation of $\Sigma_2$ placed in degree $2$.

We note that the cohomology of $\operatorname{Jac}^0(C)^2$ decomposes into invariant and anti-invariant parts under the action of $\Sigma_2$:
\[
H^\bullet(\operatorname{Jac}^0(C)^2)=H^\bullet(\operatorname{Jac}^0(C)^2)^+\oplus H^\bullet(\operatorname{Jac}^0(C)^2)^-,
\]
and therefore the contribution of the stratum $\vect_{[1,1]}^\mathrm{ab}$ to \eqref{eq:globalDT} is given by
\[
IH^\bullet(\overline{\vect}_{[1,1]},\ls_{[1,1]})\cong
\begin{cases}
	H^\bullet(\operatorname{Jac}^0(C)^2)^+[2-2g],&g\text{ odd},\\
	H^\bullet(\operatorname{Jac}^0(C)^2)^-[2-2g],&g\text{ even}.
\end{cases}
\]
An elementary calculation shows that the Hilbert series of this intersection cohomology  group is given by the second summand of the right-hand side of \eqref{eq:mainrank2}.

Finally, we show that the non-abelian stratum $\vect_{[1,1]}^{\mathrm{nab}}$ does not contribute to \eqref{eq:globalDT}. By part 0 of Theorem \ref{thm:DT}, any non-trivial local system on this stratum could only contribute to $H^i(\pi^{-1}(L\oplus L))$ 
 for degrees $i$ in the range
\[ [\dim \pms_0(2)-\dim 
\vect_{[1,1]}^{\mathrm{nab}}-\boldsymbol{r}(\pi), \dim \pms_0(2)-\dim 
\vect_{[1,1]}^{\mathrm{nab}}+\boldsymbol{r}(\pi)] = [3g-3, \, 3g-1].
\]
However, the top degree of the fiber cohomology is $2\dim \pi^{-1}(L\oplus L) = 2\dim \mathbb{P}^{g-1} = 2g-2$. For $g \ge 2$, this bound is strictly less than $3g-3$, so no contribution can arise from this stratum.

\begin{remark}
	We prove a general no-support theorem for arbitrary rank $r$ in 
	\S\ref{sec:nosupport} using completely different methods.
\end{remark}

	\section{Moduli spaces and the parabolic  projection 
	map}\label{sec:moduliandparmap}
In this section, we introduce the stratification of $\ms_0(r)$,
the local geometry along each stratum, and the parabolic projection map used throughout the paper.

\subsection{Stratification of $\ms_0(r)$}\label{sec:strat}
For a vector bundle $V$ on $C$, write
$$\mu(V):=\deg(V)/\operatorname{rk}(V)$$
for its \textit{slope}.
The bundle $V$ is \emph{stable}
(respectively, \emph{semistable}) if every proper nonzero subbundle
$V'\subset V$ satisfies $\mu(V')<\mu(V)$ (respectively,
$\mu(V')\leq\mu(V)$).
A vector bundle is called \emph{polystable} if it is isomorphic to a direct sum 
of stable bundles of the same slope.

 Every semistable bundle $V$ admits a 
Jordan--H\"older filtration
$0=W_0\subset W_1\subset\cdots\subset W_k=V$ whose quotients $W_i/W_{i-1}$ are 
stable and have the same slope as $V$. The associated graded bundle

	$$ V \mapsto \mathrm{Gr}(V):= \bigoplus_i W_i/W_{i-1}$$
	is well defined up to isomorphism, and we say that semistable vector 
	bundles $V$ and $V'$ are \textit{$S$-equivalent} (and write  $V\sim_sV'$)	
	if
	$ \mathrm{Gr}(V)\cong \mathrm{Gr}(V').$
In particular, every $S$-equivalence class contains a unique polystable bundle up to isomorphism.

Geometric invariant theory provides us with a construction of  a coarse moduli 
space  parametrizing $S$-equivalence classes of rank-$r$, degree-$d$ semistable 
vector bundles on $C$ (see \cite{MF82}).

\begin{theorem}\label{thm:moduliM}
	Let $C$ be a smooth complex projective curve of genus $g\geq 2$. There
	 is a  projective algebraic variety $\ms_d(r)$ of dimension $
		(g-1)r^2+1$, which serves as a coarse moduli space parametrizing rank-$r$, degree-$d$ semistable bundles up to $S$-equivalence.
\begin{itemize}
    \item As a set, $\ms_d(r)$ is in one-to-one correspondence with the set of isomorphism classes of polystable bundles;
    \item the subset $\ms_d(r)^{\mathrm{st}}$ of $\ms_d(r) $
		parametrizing isomorphism classes of stable bundles is a smooth, Zariski
		open subvariety; in particular, $\ms_d(r)$ is smooth if $\mathrm{gcd}(r,d)=1$. In this case, there
		exists a universal bundle $\mathcal{U}_d(r)$ over $\ms_d(r)\times C$. 
\end{itemize}
\end{theorem}
\begin{remark}\label{rem:betti-smooth}
The Poincar\'e polynomials of $\ms_d(r)$ for $\mathrm{gcd}(r,d)=1$  have been computed in
several ways (see \cite{HN74,DR75, AB83, EK00}).
\end{remark}

	From now on, we focus on the case $d=0$, and, using Theorem
	\ref{thm:moduliM}, we identify $\ms_0(r)$ with the set of isomorphism classes of
	polystable bundles.

	To every polystable bundle $V\in\ms_0(r)$, we associate the partition $\rho(V)$ of $r$ determined by the ranks of its stable direct summands. This induces a decomposition of $\ms_0(r)$ into a disjoint union of locally closed subsets:
\begin{equation}
		\label{eq:partialstrat} \ms_0(r) = \bigsqcup\limits_{\rho\vdash r}
		\vect_\rho,
	\end{equation}
where, for a partition $\rho=[r_1,\ldots,r_k]\vdash r,$
we define
\[ \vect_\rho= \{V\cong V_1\oplus\ldots\oplus V_k|\,  V_j\in
	\ms_0(r_j)^{\mathrm{st}},\,
	j=1,\ldots,k
	\}/\cong.
\]
In particular, $\vect_{[r]}=\ms_0(r)^{\mathrm{st}}$ is the stable locus.
 Given two partitions
		$\mu,\rho\vdash r$, we write
			$\mu>\rho$ if $\rho$ is a refinement of $\mu$; in this case, the closure of $\vect_\mu$ contains $\vect_\rho$.

The decomposition \eqref{eq:partialstrat} is not a stratification, however, since the automorphism group of a polystable bundle
$V\cong V_1\oplus\cdots\oplus V_k$
depends on whether some of the stable summands $V_i$ are isomorphic (see 
Proposition \ref{prop:autgroup}).
This motivates the following refinement of the notion of \textit{partition}.

A \textit{multi-partition} of $r$ is a sequence of positive integer pairs
\begin{equation}\label{eq:multipart}
	\bro=[(\br_1,m_1),(\br_2,m_2),\dots,(\br_s,m_s)],
\end{equation}
satisfying
\begin{itemize}
	\item 	 $\br_1\ge \br_2\ge\dots\ge \br_s$,
	\item    $m_1\br_1+\dots+m_s\br_s=r$, and
	\item    $\br_i=\br_{i+1}  \Rightarrow  m_i\ge m_{i+1}$.
\end{itemize}
Clearly, each multi-partition $\bro$ induces a partition by simply repeating
 $m_i$ times the part $\br_i$. We 
introduce the following notation:
\[
\rho\vdash r\text{ partition,}\quad \bro\Vdash r
\text{ multi-partition,}
\quad \bro \rightsquigarrow \rho
\text{ multi-partition induces partition}.
 \]

Now, for a multi-partition $\bro$ given as in
\eqref{eq:multipart}, we define
	\[
\vect_{\bro}=
\left\lbrace V_1^{\oplus m_1}\oplus\ldots \oplus V_s^{\oplus m_s}\in
\ms_0(r)
\mid V_i\in\ms_0(\overline{r}_i)^\mathrm{st}, V_i\not\cong V_j \text{ for
}1\le i\neq j\le s \right \rbrace.
\]
These sets do define a stratification of $\ms_0(r)$, and we have
\begin{equation}\label{eq:defstrat}
	\ms_0(r) =\bigsqcup\limits_{\bro\Vdash r}
\vect_{\bro},\quad\text{and}\quad \vect_{\rho} = \bigsqcup\limits_{\bro\rightsquigarrow
\rho}
\vect_{\bro}.
\end{equation}

	\subsection{Local structure of $\ms_0(r)$}\label{sec:localstructure} Our next goal is to describe
	the local structure of the singularities of the moduli space $\ms_0(r)$.

\begin{proposition}\label{prop:autgroup} Let $\bro$ be a multi-partition as
	in
	\eqref{eq:multipart}, and  let $V\in \vect_{\bro}$ be a polystable vector bundle.
	Then
	\begin{equation}\label{eq:autv}
	\mathrm{Aut}(V)\cong\prod_{i=1}^{s} \GL(m_i,\C).
	\end{equation}
	In particular, $\mathrm{Aut}(V)$ is abelian if and only if  $m_i=1$ for every $i$.
\end{proposition}
\begin{proof}See \cite[Propositions 4.3 and 4.4]{NS65}.
\end{proof}

\begin{definition}\label{def:abel}
	We call a polystable vector bundle $V$ \textit{abelian} if no two of its direct summands are isomorphic, i.e. if the corresponding
	multi-partition has multiplicities all equal to 1.
\end{definition}
	This
concept
induces a decomposition into abelian and non-abelian parts:
\begin{equation}\label{eq:abnab}
	\vect_\rho=\modtor\sqcup\vect_\rho^{\mathrm{nab}}.
\end{equation}
and here $\modtor$ is the abelian locus associated with the partition 
$\rho=[r_1,\ldots,r_k]\vdash r$, as defined in \S\ref{sec:summary}. It is dense 
in $\vect_\rho$, and hence
\begin{equation}\label{eq:remab}
	\overline{\modtor}=\overline{\vect}_\rho.
\end{equation}

Now we recall some well-known results about the infinitesimal structure of the
moduli space $\ms_0(r)$. For the proofs, we refer the reader to \cite{L96,N09}.
\begin{proposition} \label{prop:localet}
	Let $V \in \ms_0(r)$.  The  germ of the affine GIT quotient
	$$\locmod_V= \mathrm{Ext}^1(V,V)\sslash \mathrm{Aut}(V)$$
	at the origin is étale-locally isomorphic to the germ of \(\ms_0(r)\) at 
	\(V\).
\end{proposition}

\begin{remark}
Observe that for $t\in\C^{*}$, the element $(t,t,\dots,t)\in\aut(V)$, thought of as $t\cdot\mathrm{id}\in\prod_{i=1}^{s}\GL(m_i,\C)$, acts trivially on $ \mathrm{Ext}^1(V,V)$. Hence the action of $\aut(V)$ factors through the quotient $\aut(V)/\C^*$. Moreover, the induced action of $\aut(V)/\C^*$ is effective.
\end{remark}

\begin{remark}\label{rem:rescale}
On the other hand, the rescaling of $\mathrm{Ext}^1(V,V)$ by
$\lambda\in\C^*$ commutes
with the action of $\mathrm{Aut}(V)$, and thus descends to an action of $\C^*$
on $\locmod_V$.
\end{remark}

Combining Remark \ref{rem:rescale} with Proposition \ref{prop:purecones} gives 
the following statement.
\begin{proposition}\label{prop:corsemiproj}
 For $V\in\ms_0(r)$, endowed with the $\C^*$-action described in
	Remark \ref{rem:rescale}, the local model $\locmod_V$ is a semiprojective variety.
	In particular, its intersection cohomology $I\hbull(\locmod_V)$  has pure Hodge 
	structure.
\end{proposition}

We now describe the local model $\locmod_V$ for $\ms_0(r)$ at an abelian point $V \cong \bigoplus_{i=1}^kV_i$ more explicitly. In this setup,
$$\mathrm{Ext}^1(V,V)\cong \bigoplus_{1\leq i,j\leq k}\mathrm{Ext}^1(V_i,V_j) \text{\; and \;} \aut(V)\cong(\C^*)^{k}.$$
\begin{itemize}

  \item An element $(t_1,\ldots,t_k)\in(\mathbb C^*)^k$ acts on
$V_1\oplus\cdots\oplus V_k$ by scalar multiplication by $t_i$ on the summand $V_i$, and
    \item this induces the rescaling action by $ t_i^{-1} t_j$ on $\mathrm{Ext}^1(V_i,V_j)$.
 \end{itemize}

 Hence the local model $\locmod_V$ can be described as follows (see \cite[\S2]{ZhP}). 
\begin{proposition}\label{prop:cone}
		Let $\rho=[r_1,\ldots,r_k]\vdash r$ be a partition of $r$ and let $V\cong V_1\oplus\ldots\oplus V_k\in \modtor$.
		Then $\locmod_V$ is the product of an affine space and an affine toric variety:
		$$ \locmod_V \cong \C^{l}\times \slice_V,$$
		where $l=\sum_{i=1}^k (r_i^2(g-1)+1)$ and
		\begin{equation}\label{eq:afftoric}
			\slice_V=\bigoplus_{1\leq i\neq j\leq k}\mathrm{Ext}^1(V_i,V_j)\sslash
		(\C^*)^{k}.
		\end{equation}
The variety $\slice_V$ is an affine cone over a projective toric variety, and
is a normal slice to $\modtor$ in  $\ms_0(r)$ at $V$.
	\end{proposition}

\subsection{The moduli space $\pms_0(r)$}\label{sec:parmoduli}

Let $\pms_0(r)$  denote the moduli space parametrizing isomorphism classes of pairs
$(W,\ell),$
where $W$ is a rank-$r$ degree-$0$ vector bundle on $C$, and 
$\ell\subset W_p$ is a line, and the  pair 
$(W,\ell)$ satisfies the following stability condition:
\begin{itemize}
\item  for every proper nonzero subbundle $W'\subset W$, $\deg W'\leq 0,$
\item and the stronger inequality  $\deg W'<0$ holds whenever $\ell \subset W'_p$.
\end{itemize}

In the language of parabolic bundles (see, for example, \cite{MehtaSeshadri,BodenH,OTASz}), $\pms_0(r)$ is the moduli space of stable parabolic bundles whose parabolic structure at $p$ is given by the partial flag
$0\subset \ell\subset W_p,$
with weights $\frac{1}{2r}$ and  $\frac{1-r}{2r}$, attached respectively to $W_p/\ell$ and $\ell$. The corresponding vector of parabolic weights is  thus
$$(1/2r,\ldots,1/2r,(1-r)/2r)\in\mathbb{R}^r.$$ 
Here and below, we use the sign convention of \cite[\S2]{OTASz}, extended in 
the standard way to partial flags by allowing weights with multiplicities.

\begin{theorem}[\cite{MehtaSeshadri}]\label{thm:moduliP} The moduli space 
$\pms_0(r)$ is a smooth projective variety of dimension $r^2(g-1)+r.$ There 
exists a universal bundle $\mathcal U$ over $\pms_0(r)\times C,$ together with 
a universal line subbundle $$\mathcal{F}_1\subset {\mathcal{U}}_p, \quad 
{\mathcal{U}}_p:=\mathcal{U}\big{|}_{\pms_0(r)\times\{p\}},$$ satisfying the 
standard tautological properties. The pair $(\mathcal U,\mathcal F_1)$ is 
unique up to tensoring by the pullback of a line bundle from $\pms_0(r)$.
	In particular, there exists a unique normalized universal bundle 
	$\mathcal{U}$, characterized by the condition $\mathcal{O}\cong 
	\mathcal{F}_1\subset \mathcal{U}_p$.
\end{theorem}

We will also use the moduli space $\pms_1(r)$, whose points are represented 
	by isomorphism classes of pairs $(W,H)$, where $W$ is a stable vector 
	bundle on $C$ of rank $r$ and degree $1$, and $H\subset W_p$
	is a codimension one subspace. In terms of parabolic bundles, the parabolic 
	structure is given by the partial flag $0\subset H\subset W_p$
	and the vector of parabolic weights
	$$((1+r)/2r,1/2r,\ldots,1/2r)\in \mathbb{R}^r.$$

These two moduli spaces of stable 
parabolic bundles satisfy the following remarkable properties.
\begin{proposition}\label{prop:forgetmap}
	Let $\pms_0(r)$ and $\pms_1(r)$ be the moduli spaces introduced above. 
	\begin{enumerate}[(i)]
		\item Let $\mathcal U_1(r)$ be a universal bundle over $\ms_1(r)\times 
		C$, then 
		\begin{equation}\label{eq:pi1}
			\pms_1(r)\cong \PP\left(\mathcal 
			U_1(r)\big|_{\ms_1(r)\times\{p\}}^\vee\right)\longrightarrow\ms_1(r).
		\end{equation}
		\item If $(W,\ell)$ represents a point of $\pms_0(r)$, then $W$ is 
		semistable. Consequently, forgetting the line $\ell$ defines a proper 
		morphism of projective algebraic varieties
		\begin{equation}\label{eq:defpi}
			\pi\colon \pms_0(r)\to \ms_0(r).
		\end{equation}
		For every stable bundle $V\in \ms_0(r)^{\mathrm{st}}$, the fiber 
		$\pi^{-1}(V)$ 
		is isomorphic to the projective space $\PP(V_p)$.
		
		\item The moduli spaces $\pms_0(r)$ and $\pms_1(r)$ are isomorphic.
	\end{enumerate}
\end{proposition}
\begin{proof}
	The first two statements are immediate from the definition, while the third 
	follows 
	from the tautological Hecke correspondence (see \cite[Proposition 
	7.1]{OTASz}).
\end{proof}
\begin{corollary}\label{cor:hecke}
The moduli space $\pms_0(r)$ is isomorphic to a $\PP^{r-1}$-bundle over $\ms_1(r)$.
\end{corollary}

\section{The cohomology of the abelian fibers of the map 
$\pi$}\label{sec:parabolicfiber}

Fix a partition $\rho=[r_1,\ldots, r_k]\vdash r$, and a bundle 
$V=V_1\oplus\ldots\oplus V_k\in\modtor$, where $\rk(V_i)=r_i$, $i=1,\dots,k$. 
In this Section, we study the 
cohomology $H^{\bullet}(F_V)$ of the fiber $\pi^{-1}(V)$ of the map $\pi$.

\subsection{Towers of projective bundles}\label{sec:towerproj}
 For convenience, we fix an ordering of the direct summands of $V$ reflected by 
 the indices; thus we have a sequence \((V_1,\ldots,V_k)\) of pairwise 
 non-isomorphic stable vector 
bundles of degree $0$.
We define the subset \(\pat_{(V_1,\ldots,V_k)}\subset \pms_0(r)\) as the 
locus of pairs \((W,\ell)\) for which $W$ admits a filtration
\[
W=W_1\supset W_2\supset\cdots\supset W_{k+1}=0
\]
whose quotients satisfy
$ W_j/W_{j+1}\cong V_j$, for $j=1,\ldots,k.$

\begin{proposition}\label{prop:projtower}
The variety \(\pat_{(V_1,\ldots,V_k)}\) is an iterated projective bundle. More precisely, there is a sequence of projective bundles
\begin{equation}\label{eq:projtower}
\pat_{(V_1,\ldots,V_k)}
\xrightarrow{\psi_{k-1}}
\pat_{(V_1,\ldots,V_{k-1})}
\longrightarrow\cdots\longrightarrow
\pat_{(V_1,V_2)}
\xrightarrow{\psi_1}
\pat_{(V_1)},
\end{equation}
where $\pat_{(V_1)}\cong \PP(V_{1|p})$.
For every \(i=1,\ldots,k-1\), the morphism \(\psi_i\) is the projectivization of a vector bundle of rank
$(r_1+\cdots+r_i)r_{i+1}(g-1)+r_{i+1}.$
In particular,
$$\dim \pat_{(V_1,\ldots,V_k)} = \sum_{1\leq i<j\leq k}r_i r_j(g-1)+r-k.
$$
\end{proposition}
\noindent Informally, we will  call the varieties \(\pat_{(V_1,\ldots,V_k)}\) 
\emph{towers}.

Before proving this statement, we introduce the parabolic Hom and Ext spaces 
that will be used in the arguments.

Let $(W,\ell)$ be a vector bundle equipped with a line $\ell \subset W_p$, and 
let $\iota_p \colon \{p\} \hookrightarrow C$ denote the inclusion map. For any 
vector bundle $U$ on $C$, define the sheaf 
\[
\mathrm{\mathcal{P}ar\mathcal{H}om}((W,\ell),U):=\ker\bigl(\mathrm{\mathcal{H}om}(W,U)\longrightarrow
\iota_{p*}\operatorname{Hom}(\ell,U_p)\bigr),
\]
where the map is induced by restriction to the fiber at $p$ and then to the 
line $\ell$. The sections of $\mathrm{\mathcal{P}ar\mathcal{H}om}((W,\ell),U)$ 
are 
homomorphisms $\varphi\colon W\to U$ satisfying $\varphi_p(\ell)=0$.
We set
$$\operatorname{ParHom}((W,\ell),U) := 
H^0\bigl(C,\mathrm{\mathcal{P}ar\mathcal{H}om}((W,\ell),U)\bigr) $$
and $$\operatorname{ParExt}^1((W,\ell),U):= 
H^1\bigl(C,\mathrm{\mathcal{P}ar\mathcal{H}om}((W,\ell),U)\bigr).$$
When the line $\ell$ is clear from the context, we abbreviate these to
$\operatorname{ParHom}(W,U)$ and $\operatorname{ParExt}^1(W,U)$.
By definition, there is a short exact sequence of sheaves
\begin{equation}\label{eq:parhom-ses}
	0
	\longrightarrow
	\mathrm{\mathcal{P}ar\mathcal{H}om}((W,\ell),U)
	\longrightarrow
	\mathrm{\mathcal{H}om}(W,U)
	\longrightarrow
	{\iota_{p*}\operatorname{Hom}(\ell,U_p)}
	\longrightarrow
	0,
\end{equation}
and the associated long exact sequence in cohomology reads as
\begin{multline}\label{eq:parhom-les}
	0 \longrightarrow
	\operatorname{ParHom}((W,\ell),U)
	\longrightarrow
	\operatorname{Hom}(W,U)
	\longrightarrow
	\operatorname{Hom}(\ell,U_p)
	\longrightarrow \\
	\operatorname{ParExt}^1((W,\ell),U)
	\longrightarrow
	\operatorname{Ext}^1(W,U)
	\longrightarrow
	0.
\end{multline}

We will need  the following observation.
\begin{lemma}\label{lem:vanishing}
	Let $(W,\ell)$ represent a point in $\pms_0(r)$, and let $U$ be a stable 
	vector bundle of degree $0$. Then $\operatorname{ParHom}((W,\ell),U)=0$.
\end{lemma}
\begin{proof}
	We note that any nonzero homomorphism
	$\varphi\colon W\to U$
	must be surjective: its image is simultaneously a quotient of the 
	semistable bundle $W$ and a subbundle of the stable bundle $U$, both of 
	slope $0$.
	Consequently, $\ker\varphi$ is a degree-zero subbundle of $W$.
	If $\varphi_p(\ell)=0$, then $\ell\subset (\ker\varphi)_p$, contradicting 
	the strict inequality in the stability condition for $(W,\ell)$. Thus there 
	is no nonzero parabolic homomorphism from $(W,\ell)$ to $U$.
\end{proof}
\begin{proof}[Proof of Proposition \ref{prop:projtower}]
Set $\underline{i}:=\{1,\ldots,i\}$, and introduce the simplified notation
\[
\pat_{\underline i}:=\pat_{(V_1,\ldots,V_i)} \;\;\; \text{and} \;\;\; 
s_i:=r_1+\cdots+r_i\;\text{ for }1\leq i\leq k.
\]

We construct $\pat_{\underline i}$ inductively.
For \(i=1\), the underlying bundle is \(V_1\), which is stable, and thus the only additional datum is the choice of a line in \(V_{1|p}\) (see Proposition \ref{prop:forgetmap}(ii)). We thus have
$\pat_{\underline 1} = \PP(V_{1|p}).$
It is a simple exercise to show that the universal bundle on  
\(\pat_{\underline 1}\times C\), normalized in the sense of Theorem 
\ref{thm:moduliP}, is given by
\begin{equation}\label{eq:univW1}
\mathcal{W}_{\underline 1} := \mathcal O_{\pat_{\underline 1}}(1)\boxtimes V_1.
\end{equation}

Suppose now that the tower \(\pat_{\underline i}\) and the normalized universal bundle
\[
\mathcal{W}_{\underline i}
\longrightarrow
\pat_{\underline i}\times C
\]
have been constructed. In particular, for every point \(w\in \pat_{\underline 
i}\), the restriction
$\mathcal{W}_{\underline i}\big|_{{w}\times C}$ is the parabolic bundle represented by $w\in\pms_0(s_i)$.
Let
\[ q_i\colon \pat_{\underline i}\times C\longrightarrow \pat_{\underline i},
\qquad
u_i\colon \pat_{\underline i}\times C\longrightarrow C
\]
denote the projections, and let 
\[
\iota_i\colon \pat_{\underline i}\hookrightarrow \pat_{\underline i}\times C,
\qquad
w\longmapsto (w,p),
\]
be the natural inclusion.
Since \(\mathcal{W}_{\underline i}\) is normalized, one has the short exact 
sequence
\begin{equation}\label{eq:SESfamparhom}
0
\to
 \mathrm{\mathcal{P}ar\mathcal{H}om}
\left(\mathcal{W}_{\underline i},u_i^*V_{i+1}\right)
\to
\mathrm{\mathcal{H}om}\left(\mathcal{W}_{\underline i},u_i^*V_{i+1}\right)
\to
\iota_{i*}\left(\mathcal O_{\pat_{\underline i}}\otimes V_{i+1|p}\right)
\to
0.
\end{equation}

For every \(w\in \pat_{\underline i}\), Lemma \ref{lem:vanishing} gives
\begin{equation}\label{eq:zeroR}
H^0\left(
C, \mathrm{\mathcal{P}ar\mathcal{H}om}
\bigl(\mathcal{W}_{\underline i},u_i^*V_{i+1}\bigr)\big|_{{w}\times C}
\right)=0,
\end{equation}
and thus $$\mathcal{E}_i := R^1q_{i*}\mathrm{\mathcal{P}ar\mathcal{H}om}
\bigl(\mathcal{W}_{\underline i},u_i^*V_{i+1}\bigr) $$
is a vector bundle, whose fiber at \(w\in \pat_{\underline i}\) is
$
(\mathcal{E}_i)_w \cong \operatorname{ParExt}^1\bigl(\mathcal{W}_{\underline i}\big|_{w}, V_{i+1}\bigr).
$
We now define $$\pat_{\underline{i+1}}:= \PP(\mathcal{E}_i) 
\xrightarrow{\psi_i}  \pat_{\underline{i}}. $$

It follows from \cite{L83} that over \(\pat_{\underline{i+1}}\times C\) there is a universal short exact sequence
\begin{equation}\label{eq:univext}
0
\longrightarrow
q_{i+1}^*\mathcal O_{\pat_{\underline{i+1}}}(1)
\otimes u_{i+1}^*V_{i+1}
\longrightarrow
\mathcal W_{\underline{i+1}}
\longrightarrow
(\psi_i\times\operatorname{id})^*\mathcal W_{\underline{i}}
\longrightarrow
0.
\end{equation}
A direct check of the stability condition shows that, for every \(w\in 
\pat_{\underline{i+1}}\), the restriction  
$\mathcal{W}_{\underline{i+1}}\big|_{{w}\times C}$ represents a point of 
$\pms_0(s_{i+1})$.   Its marked line is induced from that of 
$(\psi_i\times\operatorname{id})^*\mathcal W_{\underline{i}}\big|_{{w}\times 
C}$, so $\mathcal W_{\underline{i+1}}$ is normalized.

Iterating this construction produces the sequence of projective
bundles \eqref{eq:projtower}. It remains only to compute the rank of the bundle at
each stage.
Let $w\in\pat_{\underline{i}}$; by the Riemann--Roch theorem we have
$$\chi (
\mathrm{\mathcal{H}om} (
\mathcal W_{\underline{i}}\big|_{{w}\times C},
V_{i+1}))
=s_i r_{i+1}(1-g).$$
Using \eqref{eq:parhom-les}, we then obtain
$\operatorname{rk}(\mathcal E_i) = s_i r_{i+1}(g-1)+r_{i+1}.$
This concludes the proof of Proposition \ref{prop:projtower}.
\end{proof}
\begin{corollary}
	The Poincar\'e polynomial of the tower $\pat_{(V_1,\ldots,V_k)}$ is given by
	\begin{equation}\label{eq:pttower}
			P_{t}\left(\pat_{(V_1,\ldots,V_k)}\right) = \prod_{i=1}^k
		\frac{t^{2(r_i(r_1+\ldots+r_{i-1})(g-1)+r_i)}-1}{t^2-1}.
	\end{equation}
\end{corollary}

\subsection{The cohomology ring of the projective tower}
	The following proposition describes the cohomology ring of the projective tower \eqref{eq:projtower}.
	\begin{proposition}\label{prop:ringtower}
	There is an isomorphism of graded rings
	\begin{equation*}
H^{\bullet}\left(\pat_{(V_1,\ldots,V_k)}\right)\cong \mathbb{Q}[x_1,\ldots,x_k]/I_{(r_1,\ldots,r_k)},
\end{equation*} 
  where $\deg x_i=2$, for  $i=1,\ldots,k$, and $I_{(r_1,\ldots,r_k)}$ is the ideal generated by
		$$\phi_1(x_1,\ldots,x_k)
		=  x_1^{r_1} \text{\,\,\, and\,\,\,} \phi_i(x_1,\ldots,x_k) =
		x_i^{r_i}\cdot\prod_{j=1}^{i-1}(x_i-x_j)^{r_ir_j(g-1)} \text{\,\,
		for\,\,\,} 2\leq i\leq k.$$
	\end{proposition}
	\begin{proof} It is easy to see that
	$\mathbb{Q}[x_1,\ldots,x_k]/I_{(r_1,\ldots,r_k)}$ has Hilbert function
	\eqref{eq:pttower}, and thus it is sufficient to identify the relations
	$\phi_{i}$, $i=1,\dots,k$, in the cohomology ring of $\pat_{(V_1,\ldots,V_k)}$.

We retain the notation introduced in the proof of Proposition \ref{prop:projtower} and write	 $\pat_{\underline i}:=\pat_{(V_1,\ldots,V_i)}$.
For $i=1,\ldots,k$, we then define
			\begin{equation*}
				x_{i}=c_1\bigl(\mathcal{O}_{\pat_{\underline{i}}}(1)\bigr).
			\end{equation*}
With a slight abuse of notation, we will use the same symbol for the class $x_i\in H^2(\pat_{\underline{i}})$ and for its pullback to each subsequent stage of the tower, and in particular to $\pat_{\underline{k}}=\pat_{(V_1,\ldots,V_k)}$.

	Since $\pat_{\underline{1}}\cong \PP(V_{1|p})$, we have $x_1^{r_1}=0$. Recall that for $i\geq 1$,
\begin{equation}\label{eq:vbmiddle}
\pat_{\underline{i+1}} \cong \PP(\mathcal E_i) \to \pat_{\underline{i}}, \text{\; where \;} \mathcal{E}_i := R^1q_{i*}\mathrm{\mathcal{P}ar\mathcal{H}om} \bigl(\mathcal{W}_{\underline i},u_i^*V_{i+1}\bigr).
\end{equation}

Our goal now is to calculate the Chern character of $\mathcal{E}_i$.
Using the short exact sequence \eqref{eq:SESfamparhom} together with the vanishing \eqref{eq:zeroR}, we can write
$$
\ch(\mathcal{E}_i) = - \ch\left(q_{i!} \mathrm{\mathcal{P}ar\mathcal{H}om}
\left(\mathcal{W}_{\underline i},u_i^*V_{i+1}\right)\right)= -\ch\left(q_{i!} \mathrm{\mathcal{H}om}
\left(\mathcal{W}_{\underline i},u_i^*V_{i+1}\right)\right)+r_{i+1}.
$$
Applying the Grothendieck--Riemann--Roch theorem we obtain
\begin{multline*}
-\ch\left(q_{i!} \mathrm{\mathcal{H}om} \left(\mathcal{W}_{\underline i},u_i^*V_{i+1}\right)\right) = -q_{i*}\left(\ch\left(\mathrm{\mathcal{H}om} \left(\mathcal{W}_{\underline i},u_i^*V_{i+1}\right)\right)(1-(g-1)\omega)\right),
\end{multline*}
where  $\omega\in H^2(C)$ is the class of the curve.
Note that the short exact sequence \eqref{eq:univext} gives the equality
\begin{equation*}
\ch\left( \mathcal{W}_{\underline{i}}\right)= \sum_{j=1}^i
 \exp({x_{j}})\ch(u_i^*V_{j}),
 \end{equation*}
and thus
 \begin{multline}
\ch(\mathcal{E}_i) = r_{i+1} - \sum_{j=1}^iq_{i*}\big((1-(g-1)\omega)\cdot u_i^*\ch\big(V_j^\vee\otimes V_{i+1} \big)\cdot \exp(-x_j)\big) = \\
r_{i+1} - \sum_{j=1}^iq_{i*}\big((1-(g-1)\omega)\cdot r_jr_{i+1}\cdot \exp(-x_j)\big) = r_{i+1} + \sum_{j=1}^i r_jr_{i+1}(g-1)\exp(-x_j).
 \end{multline}
Hence we arrive at the following statement.
\begin{lemma}
The formal Chern roots of the vector bundle $\mathcal{E}_i$ are:
   \begin{itemize}
\item $-x_j \text{ with multiplicity } r_jr_{i+1}(g-1)\text{ for }j=1,\ldots,i$;
\item $0\text{ with multiplicity }r_{i+1}$.
\end{itemize}
\end{lemma}
Applying the projective bundle formula to \eqref{eq:vbmiddle}, we arrive at the relation $$x_{i+1}^{r_{i+1}}\prod_{j=1}^i(x_{i+1}-x_j)^{r_jr_{i+1}(g-1)}=0,$$
which is precisely the equation for $\phi_{i+1}$.
\end{proof}
A similar calculation shows the following formula, which we will use in \S\ref{sec:locsys}.
	\begin{proposition}\label{prop:chTP}
		The Chern character of the tangent bundle of $\pat_{(V_1,\ldots,V_k)}$ is
			equal to $$\ch\left(T\pat_{(V_1,\ldots,V_k)}\right) = \sum_{1\leq i<j\leq
		k}r_ir_j(g-1)\exp({x_j-x_i})+\sum_{1\leq i\leq k}r_i\exp({x_i})-k.$$
	\end{proposition}
Proposition \ref{prop:ringtower} also yields the following
integral formula.
\begin{corollary}\label{cor:towerintegral}
	Consider the variables
	$x_1,\dots,x_k$ as elements of the second cohomology of the tower $\pat_{(V_1,\ldots,V_k)}$, as described above.
	Then, for any polynomial $\phi\in\mathbb{Q}[x_1,\ldots,x_k]$,
	\begin{multline}\label{eq:towerint}
		\operatornamewithlimits{\int}_{\pat_{(V_1,\ldots,V_k)}}\phi(x_1,\dots,x_k)=\\
		\operatornamewithlimits{\mathrm{IRes}}_{x_1,\dots ,
		x_k=0}\frac{\phi(x_1,\dots,x_k)\;dx_1\dots dx_k}{x_1^{r_1}\cdot
			x_2^{r_2}(x_2-x_1)^{r_2r_1(g-1)}\cdot\ldots\cdot
			x_k^{r_k}\prod_{i<k}(x_k-x_i)^{r_kr_i(g-1)}},
	\end{multline}
where the iterated residue is defined by the normalized nested-contour integral
\begin{equation*}
\operatornamewithlimits{\mathrm{IRes}}_{x_1,\dots,x_k=0}
\psi(x_1,\dots,x_k):=\frac{1}{(2\pi i)^k}\int_{|x_1|=\varepsilon_1}\!\cdots \int_{|x_k|=\varepsilon_k}\psi(x_1,\dots,x_k),
\end{equation*}
where $0<\varepsilon_1\ll\varepsilon_2\ll\dots\ll\varepsilon_k.$
\end{corollary}

\subsection{Intersection of projective towers} \label{sec:inttowers}
For a permutation $\sigma\in\Sigma_k$, set 
$$\pat_{\sigma}:= \pat_{(V_{\sigma(1)},\ldots,V_{\sigma(k)})}.$$
 The following statement is straightforward.
\begin{proposition}\label{prop:intuni} Let $V\cong V_1\oplus\ldots\oplus V_k \in
	\modtor$. Then
	\begin{equation*}
		\pi^{-1}(V) = \bigcup_{\substack{ \sigma\in
				\Sigma_k}}\pat_{(V_{\sigma(1)},\ldots,V_{\sigma(k)})}	
				\text{\; and \;}
		\bigcap_{\substack{ \sigma\in \Sigma_k}}
		\pat_{(V_{\sigma(1)},\ldots,V_{\sigma(k)})}  \cong
		\prod_{i=1}^k\mathbb{P}(V_{i|p}).
	\end{equation*}
\end{proposition}
Now fix $\sigma,\tau\in\Sigma_k$; our next goal is to describe the 
intersection $\pat_\sigma\cap\pat_\tau$, which will be a key input in 
\S\ref{sec:locsys}.

We begin by introducing some combinatorial data. Define a partial order $\lle$ 
on $\uka =\{1,\ldots,k\}$ as follows:
\[
j\lle i \;\;\; \text{if and only if }\;\;\;   \tau^{-1}(j)<\tau^{-1}(i)\text{ and }\sigma^{-1}(j)<\sigma^{-1}(i).
  \]
We will call a subset \(A\subset\uka\)  a \textit{lower ideal} if $a\in A$ and 
$j\lle a$ imply $j\in A$.

Next, we decompose the poset $(\uka,\lle)$ into successive layers of minimal 
elements: $\uka=\Pi_1\sqcup\cdots\sqcup\Pi_n.$
We define $\Pi_1$ to be the set of minimal elements of $\uka$, and, 
inductively, $\Pi_m$ is defined to be the set of minimal elements of 
$\uka\setminus(\Pi_1\sqcup\cdots\sqcup\Pi_{m-1})$;
we set $\Pi_{\le m}:=\Pi_1\sqcup\cdots\sqcup\Pi_m.$

We will construct a smooth projective variety $\ipat(\sigma, \tau)$ (see 
Proposition \ref{prop:inttowers} below) by induction on $m$, adding a 
layer at a time.
The variety $\ipat(\sigma, \tau)$ will come equipped with a family of vector 
bundles 
$\mathcal{W}_A$ over $\ipat(\sigma, \tau) \times C$, indexed by the lower 
ideals $A \subseteq \uka$. These vector bundles form a compatible system of 
quotients: for every inclusion $B\subset A$, there is a natural surjection 
$f_{A,B}\colon\mathcal W_A\twoheadrightarrow\mathcal W_B,$ and for $C\subset 
B\subset A$, $f_{A,C}=f_{B,C}\circ f_{A,B}.$

We begin with $\Pi_{\leq1}=\Pi_1,$ where every subset is a lower ideal.  Set
$$\ipat^{\le1}:=\prod_{i\in\Pi_1}\PP(V_{i|p}),$$ and, for $i\in\Pi_1$, let $ \mathcal W_i :=\mathcal O_{\PP(V_{i|p})}(1)\boxtimes V_i$
be the normalized universal bundle, as in \eqref{eq:univW1}. After pulling these
bundles back to $\ipat^{\le1}\times C$, define
$$\mathcal W_A:=\bigoplus_{i\in A}\mathcal W_i, \;\;\;  A\subset\Pi_1.$$
For $B\subset A$, the quotient $\mathcal W_A\twoheadrightarrow\mathcal W_B$ is the natural projection.

Suppose now that the construction has been carried out up to  level $m$.
For each element $i\in\Pi_{m+1}$,  define
$$D_i^\circ:=\{j\in\uka:j\lle i\}, \;\;\;  D_i:=D_i^\circ\sqcup\{i\}.$$
Thus $D_i$ is the principal lower ideal generated by $i$, and $D_i^\circ\subset\Pi_{\le m}$. Let
$$\ipat^{\le m}\xleftarrow{q_m} \ipat^{\le m}\times C \xrightarrow{u_m}C$$
denote the natural projections. By Lemma \ref{lem:vanishing} and equation
\eqref{eq:zeroR}, the sheaf
	$$\mathcal E_i^{\sigma\tau} := R^1q_{m*}\mathrm{\mathcal{P}ar\mathcal{H}om}
\left(\mathcal W_{D_i^\circ},u_m^*V_i\right)$$
is a vector bundle on $\ipat^{\le m}$; we denote 
$Y_i:=\PP(\mathcal E_i^{\sigma\tau}).$
On $Y_i\times C$, the universal extension sequence
\begin{equation}\label{eq:extDi}
0\to \mathcal O_{Y_i}(1)\boxtimes{V_i} \to \mathcal W_{D_i} \to \mathcal W_{D_i^\circ} \to 0
\end{equation}
defines the vector bundle $\mathcal W_{D_i}$. 
Writing $\Pi_{m+1} = \{i_1, \ldots, i_s\}$, we define the space
$$\ipat^{\le m+1} := Y_{i_1}\times_{\ipat^{\le m}}\cdots \times_{\ipat^{\le m}}Y_{i_s}.$$

Now let $A\subset\Pi_{\le m+1}$ be a lower ideal, and write
$$A=A^\circ\sqcup A', \;\;\;  A^\circ:=A\cap\Pi_{\le m}, \;\;\; A':=A\cap\Pi_{m+1}.$$
For any $i\in A'$, we have $D_i^\circ\subset A^\circ\subset \Pi_{\leq m}$, and hence the quotient map $\mathcal W_{A^\circ} \twoheadrightarrow\mathcal W_{D_i^\circ}$ is already defined. We set
\begin{equation*}
\mathcal W_{A^\circ\cup\{i\}} := \mathcal W_{A^\circ} \times_{\mathcal W_{D_i^\circ}} \mathcal W_{D_i}.
\end{equation*}
Finally, writing $A'=\{i_1,\ldots,i_t\}$, we define
\begin{equation*}
\mathcal W_A := \mathcal W_{A^\circ\cup\{i_1\}} \times_{\mathcal W_{A^\circ}} \cdots \times_{\mathcal W_{A^\circ}} \mathcal W_{A^\circ\cup\{i_t\}}.
\end{equation*}

It is a simple exercise to show that the quotient maps associated with inclusions of lower ideals are induced by the construction above and are compatible with composition.

We now identify the resulting variety $\ipat(\sigma,\tau)$ 
with the scheme-theoretic intersection of the two towers,  
$\pat_\sigma$ and 
$\pat_\tau$.

\begin{proposition}\label{prop:inttowers}
Let $\rho\vdash r$ be a partition, and let $V\cong V_1\oplus\ldots\oplus V_k$ be a point in $\modtor$. For $\sigma,\tau\in\Sigma_k$, the scheme-theoretic intersection of the two
towers $\pat_{\tau}$ and $\pat_{\sigma}$   is isomorphic to
the space constructed above:
$\pat_\sigma\cap\pat_\tau\cong\ipat(\sigma,\tau).$
In particular, $\pat_\sigma\cap\pat_\tau$ is a smooth iterated projective bundle, and the towers $\pat_\sigma$ and $\pat_\tau$ intersect
cleanly.
\end{proposition}
\begin{proof}
Let $\VV_{\uka}$ be the restriction of the normalized universal bundle  on $\pms_0(r)\times C$ to $(\pat_\sigma\cap\pat_\tau)\times C$. By the definition of the two towers,
it carries two filtrations  \begin{equation}\label{eq:twofilt}
\VV_{\uka}=\VV^\sigma_1\supset\VV^\sigma_2\supset\cdots \supset\VV^\sigma_{k+1}=0, \;\;\; 
\VV_{\uka}=\VV^\tau_1\supset\VV^\tau_2\supset\cdots \supset\VV^\tau_{k+1}=0,
\end{equation}
whose successive quotients $\VV^\sigma_s/\VV^\sigma_{s+1}$ and $\VV^\tau_s/\VV^\tau_{s+1}$ have fibers $V_{\sigma(s)}$ and $V_{\tau(s)}$, respectively.

\begin{lemma}\label{lem:quotients}
The two filtrations \eqref{eq:twofilt} determine, for every lower
ideal $A\subset\uka$, a quotient $f_{\uka, A}\colon\VV_{\uka}\twoheadrightarrow\VV_{A}.$
These quotients are compatible with inclusions of lower ideals. Moreover, for $s=0,\ldots,k$, the lower ideals
$$A^\sigma_s:=\{\sigma(1),\ldots,\sigma(s)\} \text{\; and \;}  A^\tau_s:=\{\tau(1),\ldots,\tau(s)\}$$
give precisely the quotients determined by the two filtrations
 \eqref{eq:twofilt}.
\end{lemma}
\begin{proof}
For $i\in\uka$, recall that $D_i$ denotes the principal lower ideal generated by $i$. 
We note that for $s=\sigma^{-1}(i)$ and $t=\tau^{-1}(i)$,  $$D_i=A_s^\sigma\cap A_t^\tau \text{\; and define \;}
 \VV_{D_i}:= \VV_{\uka}/(\VV^\sigma_{s+1}+\VV^\tau_{t+1}).$$ For an arbitrary lower ideal $A\subset \uka$, using 
$A=\bigcup_{i\in\max(A)}D_i,$ set 
$$\mathcal{K}_A := \bigcap_{i\in\max(A)}\ker(\VV_{\uka}\to\VV_{D_i}), \text{\; and \;} 
\VV_A:=\VV_{\uka}/\mathcal K_A.$$ 

The vanishing $\operatorname{Hom}(V_i,V_j)=0$ for $i\neq j$ implies that $\VV_A$ is the unique quotient associated with $A$.  If $B\subset A$, then $\mathcal K_A\subset \mathcal K_B$, giving a
natural quotient map $\VV_A\twoheadrightarrow\VV_{B}.$ These maps are compatible with composition, and for $A=A_s^\sigma$ or $A=A_s^\tau$, the construction recovers the corresponding
quotient in  \eqref{eq:twofilt}.
\end{proof}
The compatible quotients $\VV_A$  define a morphism
\begin{equation}\label{eq:onedir}
\pat_\sigma\cap\pat_\tau\longrightarrow\ipat(\sigma,\tau).
\end{equation}
Indeed, suppose inductively that a morphism to $\ipat^{\le m}$ has been constructed. For each $i\in\Pi_{m+1}$, the quotient $\VV_{D_i}\twoheadrightarrow\VV_{D_i^\circ}$ determines an extension of $\VV_{D_i^\circ}$ by $u_m^*V_i$, and hence, by the universal property of
$Y_i=\PP(\mathcal E_i^{\sigma\tau}),$ a morphism to $Y_i$. Since these morphisms lie over the same map to $\ipat^{\le m}$, they together define a morphism to $\ipat^{\le m+1}$. Proceeding inductively, we obtain  \eqref{eq:onedir}.

Conversely, the quotient bundles $\mathcal W_{\uka}\twoheadrightarrow \mathcal W_{A^\sigma_s}$ and $\mathcal W_{\uka}\twoheadrightarrow \mathcal W_{A^\tau_s}$, for $s=0,\ldots,k$, endow the universal family $\mathcal W_{\uka}$ on $\ipat(\sigma,\tau)\times C$ with a
$\sigma$- and a $\tau$-filtration. By the definition of the two towers, they therefore determine a morphism
 \begin{equation}\label{eq:twodir}
\ipat(\sigma,\tau)\to\pat_\sigma\cap\pat_\tau.
\end{equation}

By construction, the morphisms \eqref{eq:onedir} and \eqref{eq:twodir} are inverse to one another. Hence
$\pat_\sigma\cap\pat_\tau\cong\ipat(\sigma,\tau)$
as schemes, and in particular the scheme-theoretic intersection $\pat_\sigma\cap\pat_\tau$ is smooth. Thus $\pat_\sigma$ and $\pat_\tau$ intersect cleanly.
\end{proof}

\begin{remark}\label{rem:IntTowers}
 The argument of Proposition \ref{prop:inttowers} extends to the intersection of any finite collection of towers   $\bigcap_{\sigma\in S}  \pat_{\sigma}$,  $S\subseteq \Sigma_k$. One simply replaces the partial order $\lle$ by the common partial order  $j\prec_S i \Longleftrightarrow  \sigma^{-1}(j)<\sigma^{-1}(i)$ for every $\sigma\in S$. 
\end{remark}

		\subsection{The Poincar\'e polynomial of the abelian fiber}
We now pass from the description of the individual towers to the abelian 
fiber  itself.

The purity of the cohomology of the abelian fiber $F_V\!=\!\pi^{-1}(V)$ is a 
key fact in our argument.
\begin{proposition} \label{prop:pureabelianfibers}
	Let $\rho=[r_1,\ldots,r_k]\vdash r$ be a partition and let
	$V\in\modtor$. Then the mixed Hodge structure on $\hbull(F_V)$ is pure.
\end{proposition}
\begin{proof}
	We argue by induction on the length \(k\) of \(\rho\), using Proposition
	\ref{prop:purefibers}. If \(k=1\), then \(F_V\cong\PP^{r-1}\), and the
	result is immediate.
	
	Assume \(k>1\), and let \(\mu=[m_1,\ldots,m_n]\) be a strict coarsening
	of \(\rho\). Thus \(n<k\), and, for
	\(U=U_1\oplus\cdots\oplus U_n\in\vect_\mu^{\mathrm{ab}}\), the induction
	hypothesis implies that \(H^\bullet(F_U)\) is pure.
	
	We apply Proposition \ref{prop:purefibers} to the stratum
		$Y_\mu=\vect_\mu^{\mathrm{ab}}$, taking $Z_\mu=\Sprod_\mu$ and 
		$\gamma_\mu$
		as in \eqref{eq:gammarho}; then
		$\widetilde Z_\mu=\gamma_\mu^{-1}(\vect_\mu^{\mathrm{ab}})
		=\mathcal S_\mu^{\mathrm{ab}}$, and condition \textup{(i)} holds because
		$\widetilde\gamma_\mu$ is a finite \'etale cover.
	
	We next show that
	\(\widetilde\gamma_\mu^*\ls_\mu\) is constant.
	For each ordering of the summands of \(U\), the corresponding tower
	maps properly and birationally onto the associated irreducible component
	of \(F_U\) (see Proposition \ref{prop:intuni}). Since \(H^\bullet(F_U)\) is
	pure, Corollary
	\ref{cor:pureinject} therefore gives an injection
	\begin{equation}\label{eq:towerinject}
		H^\bullet(F_U)\hookrightarrow
		\bigoplus_{\sigma\in\Sigma_n}
		H^\bullet\left(
		\pat_{(U_{\sigma(1)},\ldots,U_{\sigma(n)})}
		\right).
	\end{equation}
	
	After pulling back to the ordered cover
	\(\mathcal S_\mu^{\mathrm{ab}}\), these injections
	assemble into an injection of local systems. By the projective-bundle
	description of the towers, the local system on the right is constant,
	hence the pullback to \(\mathcal S_\mu^{\mathrm{ab}}\) of 
	the cohomology local system
	with stalk \(H^\bullet(F_U)\) is constant, and so is its summand
	\(\widetilde\gamma_\mu^*\ls_\mu\). This verifies condition
	\textup{(ii)} of Proposition \ref{prop:purefibers}.
	
	Finally, a point \( z\in\gamma_\mu^{-1}(V)\) determines an ordered
	decomposition
	\[
	\{1,\ldots,k\}=I_1\sqcup\cdots\sqcup I_n,
	\qquad
	\sum_{i\in I_a}r_i=m_a.
	\]
	Setting \(V_{I_a}=\bigoplus_{i\in I_a}V_i\), the normal slice at
	\(z\) decomposes as the product of toric varieties described in Proposition
	\ref{prop:cone}, and its intersection cohomology is pure 
	by
	Proposition \ref{prop:corsemiproj} together with the Künneth formula. This
	verifies condition \textup{(iii)}.
	
	Thus all the conditions of Proposition \ref{prop:purefibers} hold for
	every strict coarsening \(\mu\) of \(\rho\), and consequently
	\(H^\bullet(F_V)\) is pure.
\end{proof}

We will need a slight extension of this result. For \(i=1,\ldots,k\), set
	\begin{equation}\label{eq:introducecover}
		\patf_i:=
		\bigcup_{\substack{\sigma\in\Sigma_k\\ \sigma(k)=i}}
		\pat_{(V_{\sigma(1)},\ldots,V_{\sigma(k)})}
		\subset F_V
	\end{equation}
to be the set of parabolic bundles containing \(V_i\) as a
	subbundle, and, for every nonempty subset \(J\subset\uka=\{1,\ldots,k\}\), 
	define
	\[
	\patf_J:=\bigcap_{j\in J}\patf_j,
	\qquad
	V_J:=\bigoplus_{j\in J}V_j,
	\qquad
	r_J:=\sum_{j\in J}r_j.
	\]
	
	\begin{lemma}\label{lem:pureintersections}
		For every nonempty subset \(J\subset\uka\), the mixed Hodge structure on
		\(H^\bullet(\patf_J)\) is pure.
	\end{lemma} 
	\begin{proof}
Let $\mathcal U_J :=\mathcal U\big|_{F_{V/V_J}\times C}$ be the restriction of the normalized universal bundle on \(\pms_0(r-r_J)\times C\), and consider the projections 
$$F_{V/V_J}\ \xleftarrow{\ q\ } \ F_{V/V_J}\times C \ \xrightarrow{\ u\ }\ C.$$ For \(j\in J\), set $$\mathcal E_{J,j}:=R^1q_*\mathrm{\mathcal{P}ar\mathcal{H}om}(\mathcal U_J,u^*V_j).$$  By Lemma \ref{lem:vanishing}, this is a vector bundle of rank $r_j(r-r_J)(g-1)+r_j.$ Writing \(J=\{j_1,\ldots,j_s\}\), the  tower construction of Proposition \ref{prop:projtower} yields an identification
\begin{equation}\label{eq:FJproj}
\patf_J\cong\PP(\mathcal E_{J,j_1})\times_{F_{V/V_J}}\cdots\times_{F_{V/V_J}}\PP(\mathcal E_{J,j_s}).
\end{equation}
Here, when \(J=\uka\), we use the conventions \(F_{V/V_J}=\mathrm{pt}\) and \(\mathcal E_{J,j}=V_{j|p}\).
	
By Proposition \ref{prop:pureabelianfibers} the fiber \(F_{V/V_J}\) has pure cohomology. Applying the projective bundle formula successively to the description \eqref{eq:FJproj}, we conclude that $\patf_J$ also has pure cohomology. \end{proof}
We now use these purity properties to study the \textit{Mayer--Vietoris complex} associated with the covering
$F_V=\bigcup_{i=1}^k \patf_i.$

 Let $e_1,\ldots,e_k$ denote the standard basis of $\Q^k$, and for $J=\{j_1<\cdots<j_s\}$ write $e_J=e_{j_1}\wedge\cdots\wedge e_{j_s}$.
The corresponding augmented Mayer--Vietoris complex is
\begin{equation}\label{eq:MV}
	0\to\hbull(F_V)\to \bigoplus_{|J|=1}e_J\otimes\hbull(\patf_J) \to \bigoplus_{|J|=2}e_J\otimes\hbull(\patf_J) \to\cdots\to e_{\uka}\otimes\hbull(\patf_{\uka}) \to0,
\end{equation}
with differential $d_{\mathrm{MV}}=\sum_{j=1}^k\wedge e_j\otimes\iota_j^*$, where $\iota_j^*$ denotes  the canonical restriction map.
\begin{proposition}\label{prop:MVexact}
The complex \eqref{eq:MV} is exact.
\end{proposition}
\begin{proof}
Its $E_1^{p,q}$-terms are direct sums of
the groups $H^q(\patf_J)$ with $|J|=p+1$, and hence are pure of weight
$q$. By strictness, the same is true for the $E_2^{p,q}$-terms. Since all
higher differentials change the weight, they vanish, and the spectral sequence
degenerates at $E_2$.

On the other hand, the abutment $H^{p+q}(F_V)$  is pure of weight $p+q$ by
Proposition \ref{prop:pureabelianfibers}. Thus $E_2^{p,q}=E_\infty^{p,q}=0$ for $p>0.$
Hence each row of the augmented $E_1$-complex is exact, which is precisely
the exactness of \eqref{eq:MV}.
\end{proof}
As a corollary, we obtain the following statement. 

\begin{proposition}\label{cor:pureven} 
For $V\in\modtor$, the cohomology $H^\bullet(F_V)$ vanishes in  odd degrees, 
and the weight polynomial  of $F_V$ (see \S\ref{sec:wp}) coincides with its 
Poincar\'e polynomial: $E_t(F_V)=P_t(F_V)$. Similarly, for every nonempty 
subset \(J\subset\uka\), the cohomology $\hbull(\patf_J)$ vanishes in  odd 
degrees, and  $E_t(\patf_J)=P_t(\patf_J)$.
\end{proposition}
\begin{proof}
We argue by induction on $k$. For $k=1$, we have $F_V\cong\patf_{\{1\}}\cong \PP^{r-1}$, so the statement is immediate. 

Assume now that $k>1$. For every nonempty $J\subset\uka$, the quotient ${V/V_J}$ involves fewer than $k$ summands. Hence, by the inductive hypothesis, the cohomology of $F_{V/V_J}$ is concentrated in even degrees. Applying the projective bundle formula successively to the description \eqref{eq:FJproj}, we conclude that the same is true for $\hbull(\patf_J)$.

Since every term following  $H^\bullet(F_V)$ in the exact Mayer--Vietoris complex \eqref{eq:MV} is concentrated in even degrees, it follows that $H^\bullet(F_V)$ is concentrated in even degrees as well. Proposition \ref{prop:pureabelianfibers} and Lemma \ref{lem:pureintersections} then imply that, for $F_V$ and every $\patf_J$, the weight polynomial coincides with the Poincar\'e polynomial.
\end{proof}

With these preparations, we are now ready to prove a crucial recursion on the Poincar\'e polynomials of abelian fibers.
	\begin{theorem}\label{thm:PoinFib}
		Let $\rho=[r_1,\ldots,r_k]\vdash r$ be a partition of $r$, $V\in\modtor$ and let $f(\rho;t)$ be the polynomial such that
		$f(\rho;t^2)=P_t(F_V)$.
For any subset
		$J=\{j_1,j_2,\ldots,j_s\}\subset
		\{1,2,\ldots,k\}$, denote by $\overline{J}=\{1,\ldots,k\}\setminus J$ its
		complement, and by $\rho|J = [r_{j_1},\ldots,r_{j_s}]$ the corresponding
		partition of $r_J=\sum_{j\in J} r_j$.
		Setting $f(\rho_{\emptyset};t)=1$, and $r_\emptyset=0$, we have  the recursion
\begin{equation}\label{eq:frec-JT}
		f(\rho;t) = \sum_{\emptyset\neq J\subset\{1,\ldots,k\}} (-1)^{|J|-1}
		f(\rho|{\overline{J}};t)\prod_{j\in J}\frac{t^{r_jr_{\overline{J}}(g-1)+r_j}-1}{t-1}.
\end{equation}
\end{theorem}
\begin{proof}
For every nonempty $J\subset\uka$, the projective-bundle description \eqref{eq:FJproj}, together with the purity and evenness of $\hbull(\patf_J)$ noted in Proposition \ref{cor:pureven},  give
\[
E_t(\patf_J)=P_t(\patf_J) = P_t(F_{V/V_J}) \prod_{j\in J} \frac{t^{2(r_jr_{\overline J}(g-1)+r_j)}-1}{t^2-1}.
\]
Since the sets $\patf_i$ cover $F_V$, inclusion--exclusion for the $E$-polynomials yields 
\begin{equation}\label{eq:Etincexc}
E_t(F_V) = \sum_{\emptyset\ne J\subset\{1,\ldots,k\}} (-1)^{|J|-1}E_t(\patf_J).
\end{equation}

By Proposition \ref{cor:pureven}, we have $E_t(F_V) = P_t(F_V)=f(\rho; t^2)$. Substituting this into \eqref{eq:Etincexc}, using  $P_t(F_{V/V_J}) = f(\rho|{\overline J}; t^2)$, and replacing $t^2$ by $t$, we arrive at \eqref{eq:frec-JT}.
\end{proof}

\subsection{Global degree-two classes on abelian fibers}\label{sec:cohringsfib}
Repeating the argument from the proof of Proposition \ref{prop:inttowers}, for every $\sigma\in\Sigma_k$ and every $s=0,\ldots,k$, we obtain vector bundles
$$\VV_{A^\sigma_s}\longrightarrow\pat_\sigma\times C, \text{\; where \;}  A^\sigma_s=\{\sigma(1),\ldots,\sigma(s)\}.$$
For $\sigma=\mathrm{id}$, these are precisely the bundles $\mathcal W_{\underline{s}}$ appearing in the construction of the tower $\pat_{\mathrm{id}}$ in \S\ref{sec:towerproj}. We shall therefore identify the two notations.
 Set
\[\mathcal E^\sigma_1:=V_{\sigma(1)|p}, \;\; \mathcal E^\sigma_i := R^1q_* \mathrm{\mathcal{P}ar\mathcal{H}om} \left(\VV_{A^\sigma_{i-1}},u^*V_{\sigma(i)}\right), \;\; 2\leq i\leq k,\]
where $q$ and $u$ denote the natural projections from $\pat_\sigma\times C$, and define
\[ x^\sigma_{\sigma(i)} := c_1\left(\mathcal O_{\PP(\mathcal E^\sigma_i)}(1)\right), \;\; 1\leq i\leq k.\]
By Proposition \ref{prop:ringtower}, these classes generate the cohomology ring of $\pat_\sigma$.

\begin{proposition}\label{prop:xsigma}
There exist classes
$ x_1,\ldots,x_k\in H^2(F_V)$
such that, for every $\sigma\in\Sigma_k$ and every $i=1,\ldots,k$,
$x_i\big|_{\pat_\sigma}=x^\sigma_i.$
\end{proposition}

\begin{proof}
Consider the Mayer--Vietoris spectral sequence associated with the cover $F_V=\bigcup_{\sigma\in\Sigma_k}\pat_\sigma.$ By Proposition \ref{prop:inttowers} and Remark \ref{rem:IntTowers}, all multiple intersections of the towers are iterated projective bundles and hence have pure cohomology. As in the proof of Proposition \ref{prop:MVexact}, the  spectral sequence therefore degenerates at $E_2$.
		
For each $i$, the classes $x_i^\sigma$, $\sigma\in\Sigma_k$, determine an 
element of $E_1^{0,2}=\bigoplus_{\sigma\in\Sigma_k}H^2(\pat_\sigma).$ 
To prove our statement, it is enough to show  that, for every 
$\sigma,\tau\in\Sigma_k$,
\begin{equation}\label{eq:d1closed}
			x^\sigma_i\big|_{\pat_\sigma\cap\pat_\tau}= x^\tau_i\big|_{\pat_\sigma\cap\pat_\tau}.
\end{equation}
Indeed, this shows that the element $\{x_i^\sigma\}_{\sigma\in \Sigma_k}$ is a class in $E_2^{0,2}=E_\infty^{0,2}$, and hence the restriction of a class from $H^2(F_V)$.

Fix $\sigma,\tau\in\Sigma_k$ and set $s=\sigma^{-1}(i)$,  $t=\tau^{-1}(i),$  $D_i^\circ=A^\sigma_{s-1}\cap A^\tau_{t-1}.$
Suppose first that $D_i^\circ\neq\varnothing$. The quotient maps $$ \VV_{A^\sigma_{s-1}}\to \VV_{D_i^\circ}, \;\;\;   \VV_{A^\tau_{t-1}}\to\VV_{D_i^\circ}$$ fit into short exact sequences whose kernels satisfy the $R^0q_*\mathrm{\mathcal{P}ar\mathcal{H}om}(\,\cdot\,,u^*V_i)$-vanishing by Lemma \ref{lem:vanishing}. The associated long exact sequences therefore give injections
\[\PP(\mathcal E^\sigma_{\sigma^{-1}(i)}) \hookleftarrow 
		\PP(R^1q_*\mathrm{\mathcal{P}ar\mathcal{H}om}(\VV_{D_i^\circ},u^*V_i))\!\hookrightarrow
		\PP(\mathcal E^\tau_{\tau^{-1}(i)}).\]

If $D_i^\circ=\varnothing$, then $i$ is minimal for $\lle$, and  \eqref{eq:parhom-les} gives  the corresponding maps
$$\PP(\mathcal{E}^\sigma_{\sigma^{-1}(i)}) \hookleftarrow  \PP(V_{i|p})\hookrightarrow\PP(\mathcal{E}^\tau_{\tau^{-1}(i)}).$$ 

By the inductive construction of  $\ipat(\sigma, \tau)\cong \pat_\sigma\cap\pat_\tau$, the projective-bundle stage corresponding to $i$ is, in the first case, $\PP\left(R^1q_*\mathrm{\mathcal{P}ar\mathcal{H}om}(\VV_{D_i^\circ},u^*V_i)\right),$ and, in the second case, $\PP(V_{i|p})$. Thus the tautological line bundles defining $x_i^\sigma$ and $x_i^\tau$ both restrict to the same tautological line bundle on $ \pat_\sigma\cap\pat_\tau$; hence \eqref{eq:d1closed} holds.
\end{proof}

For completeness, we also give a full description of cohomology rings of
abelian fibers in terms of the classes introduced above.
This result is not used elsewhere in the paper; its proof is given in the supplementary subsection
\ref{sec:fiberring}.

\begin{theorem}\label{thm:fiberring} Let $\rho=[r_1,\ldots,r_k]\vdash r$ be a partition of $r$, and  let $V$ be a point in $\modtor$.
		The cohomology ring $H^\bullet(\pi^{-1}(V))$ is generated by $k$ 
		elements $x_1,\ldots,x_k\in H^2(\pi^{-1}(V))$ with the only relations 
		\begin{equation}\label{eq:HStrelationsFib}
			\phi_{A\sqcup  B}(x_1,\ldots,x_k):=\prod_{i\in A}x_i^{r_i}\cdot\prod_{\substack{i\in
					A, j\in B}}(x_i-x_j)^{r_{i}r_j(g-1)},
		\end{equation} where $A\sqcup  B$  runs  	over all  
		decompositions of the set $\{1,2,\ldots,k\}$, such that 
		$A\neq\emptyset$.
	\end{theorem}

 We end this section with a statement that we will need later, in \S\ref{sec:locsys}. This formula uses Proposition \ref{prop:xsigma}, and is a result of a simple calculation, analogous to the one in the proof of Proposition  \ref{prop:ringtower}.
\begin{proposition}\label{prop:chTInt}
Let $\sigma, \tau\in\Sigma_k$ and let $\pat_\sigma, \pat_\tau$ be two towers of projective bundles defined in \S\ref{sec:inttowers}. 
The Chern character of the tangent bundle of $\pat_{\tau}\cap \pat_\sigma$ is
$$\ch(T(\pat_{\tau}\cap \pat_\sigma))=\sum_{\substack{\tau^{-1}(i)<\tau^{-1}(j) \\ \sigma^{-1}(i)<\sigma^{-1}(j)}} r_{i}r_{j}(g-1)\exp(x_{j}-x_{i})  +\sum_{1\leq i\leq k}r_i\exp(x_i) -k.$$
\end{proposition}

	\section{The no-support theorem for non-abelian strata}
\label{sec:nosupport}

Recall the setup of \S\ref{sec:strat}: let $ \bro=[(\br_1,m_1),(\br_2,m_2),\dots,(\br_s,m_s)]$ be a multipartition of $r$, let $\wV\in\vect_{\bro}$, and let $\mathcal C$ be a sufficiently small analytic representative of the normal slice to
 $\vect_{\bro}$ in $\ms_0(r)$ at $\wV$. Set \(\ccirc=\mathcal C\setminus\{\wV\}\), and let \({\jmap}\colon\pi^{-1}(\ccirc)\hookrightarrow\pi^{-1}(\mathcal C)\) be the inclusion map. Via the identification $\hbull(F_{\wV})\simeq \hbull(\pi^{-1}(\mathcal C))$, we define
\[
L_{\bro} = \ker\left(\hbull(F_{\wV})\overset{\jmap^*}\longrightarrow\hbull(\pi^{-1}(\ccirc))\right).
\]
Locally near $\wV$, a neighborhood of the stratum $\vect_{\bro}$ decomposes as the product of $\mathcal C$ with a contractible neighborhood in $\vect_{\bro}$, ensuring that the local restriction map remains unchanged after passing to the normal slice. By \eqref{eq:cohrestr}, $L_{\bro}$ is therefore precisely the stalk at $\wV$ of the graded local system supported on $\vect_{\bro}$ that appears in the Decomposition Theorem \ref{thm:DT} for the map $\pi$.

The goal of this section is to prove the following statement.
\begin{theorem}\label{thm:nosupp}
If $\bro$ is a non-abelian multipartition, i.e. \(m_i>1\) for some \(i\), then
$L_{\bro}
=0$, and thus $\vect_{\bro}$ carries no support of
	$R\pi_*\Q_{\pms_0(r)}$.
\end{theorem}

Recall from \S\ref{sec:strat} that the multipartition $\bro$ determines a partition $\rho=[r_1,\dots,r_k]$, $k=m_1+\ldots+ m_s$.
Let $V\in\modtor\cap\ccirc,$ and consider the composition
\begin{equation}\label{eq:composed}
	\hbull(F_{\wV})
\xrightarrow{\jmap^*}
\hbull(\pi^{-1}(\ccirc))
\longrightarrow
\hbull(F_V),
\end{equation}
where both arrows are induced by the natural inclusions. To prove that $\jmap^*$ is injective, it is enough to show that the composed map \eqref{eq:composed} is injective.

Recall from \S\ref{sec:summary} that $\modtor$ is obtained as the quotient of $\modtir$ by the action of the group $\aut(\rho)\subset\Sigma_k$. By a slight abuse of notation, $\aut(\rho)$ also acts on the cohomology of the fiber \(F_V\).
We define the subgroup
\[
\autbro\simeq\Sigma_{m_1}\times\cdots\times\Sigma_{m_s}\subset\aut(\rho)
\subset\Sigma_k,
\]
which permutes the summands within each of the $s$ blocks. In particular, if 
\(\bro\) is abelian, then \(\autbro=\{\mathrm{id}\}\). 

Since permuting the  summands within a block does not affect their specialization to \(\wV\), the image of the composed restriction map \eqref{eq:composed} is invariant under $\autbro$.  We thus obtain a natural map
\begin{equation}
	\spmap\colon\hbull(F_{\wV})\longrightarrow \hbull(F_V)^{\autbro}.
\end{equation}

\begin{proposition}\label{prop:ginjective}
	The map $\spmap$ is injective.
\end{proposition}
Theorem \ref{thm:nosupp} clearly follows from this statement.
A motivic argument establishing that $\spmap$ is surjective --- and therefore an 
isomorphism --- is given in Appendix \ref{sec:suppisom}.

The remainder of the section is devoted to the proof of
Proposition \ref{prop:ginjective}.

\subsection{Injectivity and the generalized Mayer--Vietoris complex}\label{sec:injcohmap}

As in the abelian case, a key fact is the purity
of the fiber cohomology.

\begin{proposition}\label{lem:purenonab} The mixed Hodge structure on $H^\bullet(F_{\wV})$ is pure.
\end{proposition}
\begin{proof}
		We argue by induction on the closure order of the
		non-abelian strata. In a sufficiently small normal slice to
		$\vect_{\bro}$, the possible supports are abelian strata and
		non-abelian strata $\vect_{\overline{\mu}}$ satisfying $ 
		\vect_{\bro}\subsetneq\overline{\vect}_{\overline{\mu}}$, and the 
		latter 
		carry no
		support by induction. For each  abelian support, the finite
		morphism \eqref{eq:gammarho} trivializes the relevant local systems on the
		ordered cover. Their stalks are pure by Saito's decomposition theorem, and
		the  normal slices on the cover have pure intersection cohomology by
		Proposition \ref{prop:corsemiproj}. Thus the hypotheses of 
		Proposition \ref{prop:purefibers} are satisfied, and we conclude that 
		$H^\bullet(F_{\wV})$, for $\wV\in\vect_{\bro}$,  is pure. 
\end{proof}
To simplify the notation, we describe in detail the
case
\[
\bro=[({\br},k)],\quad {\br} k=r,\;\; \text{and}\;\; \rho=[{\br},\ldots,\br].
\]
Thus, the corresponding stratum $\vect_{\bro}$ consists of bundles  $\wV=U^{\oplus k}$, where \(U\) is a stable bundle of rank $\br$ and degree zero.  We will prove the injectivity of
\begin{equation}\label{eq:gmap}
\spmap\colon \hbull(F_{U^{\oplus k}}) \to \hbull(F_V)^{\Sigma_k},
\end{equation}
 where $V\in\modtor$, by induction on $k$. The case \(k=1\) is immediate.

Now let $V=V_1\oplus\cdots\oplus V_k\in\modtor$.  By Proposition
\ref{prop:MVexact}, the Mayer--Vietoris complex \eqref{eq:MV} associated to 
the covering $F_V=\bigcup_j\patf_j$ is exact, and the differential
$d_{\mathrm{MV}}$ is $\Sigma_k$-equivariant. Taking $\Sigma_k$-invariants therefore gives the
following.
\begin{corollary}\label{cor:MVembed}
The first map in \eqref{eq:MV} induces an embedding
\begin{equation}
\hbull(F_V)^{\Sigma_k}
\hookrightarrow
\left[
\bigoplus_{j=1}^k\hbull(\patf_j)
\right]^{\Sigma_k}.
\end{equation}
Note that we omitted the tensor factors \(e_j\).
\end{corollary}

Next, we introduce a key object in our argument: the moduli space of parabolic extensions
\[
\widehat{\mathcal E}
=
\{ 0\to V''\to \widehat V\to V'\to0\mid
V''\in\ms_0({\br}),\,
V'\in\pms_0((k-1){\br}),\,
\widehat V\in\pms_0(k{\br})
\}.
\]
This space has two natural morphisms 
\begin{equation}\label{eq:structure} 
\pms_0((k-1){\br})
\overset{f'}{\longleftarrow}
\widehat{\mathcal E}
\overset{h'}{\longrightarrow}
\pms_0(k{\br})
\overset{\pi}{\longrightarrow}
\ms_0(k \br),
\end{equation}
where $f'$ sends an extension to its quotient $V'$, while $h'$ sends it to its middle term.
Set
\(\tip:=\pi\circ h'. \)
Note that for an abelian point \(V\in\modtor\) its preimage decomposes as
\(\tip^{-1}(V)=\bigsqcup_{j=1}^k\patf_j\), where $\patf_j$ is the component consisting of parabolic bundles containing $V_j$  as a subbundle.

Now we consider the non-abelian point \(\wV=U^{\oplus k}\), and set
\[
\ParExtSp
=
\tip^{-1}(\wV)
=
\{0\to U\to W\to W'\to 0
\mid
W'\in F_{U^{\oplus(k-1)}},\
W\in\pms_0(k\br)
\}.
\]
The restriction of  \eqref{eq:structure} to this fiber gives
\begin{equation}\label{eq:structurebar}
	F_{U^{\oplus(k-1)}}
	\overset{f}\longleftarrow
	\ParExtSp
	\overset{h}\longrightarrow
	F_{U^{\oplus k}}\overset{\pi}\longrightarrow
	\{U^{\oplus k}\}.
\end{equation}
It follows from Lemma \ref{lem:vanishing} that  $f$ is a projective bundle 
morphism. By induction, $\ParExtSp$ is irreducible. Since every point of
$F_{U^{\oplus k}}$ contains a subbundle isomorphic to $U$ arising from its
Jordan--H\"older filtration, $h$ is surjective and
$F_{U^{\oplus k}}$ is irreducible. The locus of bundles admitting more than
one such subbundle is proper; hence a generic point admits a unique one, and
$h$ is birational.

Since $H^\bullet(F_{U^{\oplus k}})$ is pure by Proposition \ref{lem:purenonab}, Corollary \ref{cor:pureinject} implies that the pullback map
\[
h^*:\hbull(F_{U^{\oplus k}})\longrightarrow\hbull(\ParExtSp)
\]
is injective; the projective-bundle formula and induction then
show that $H^\bullet(F_{U^{\oplus k}})$ is concentrated in even degrees.

 Recall  that $\modtor$ is the quotient of the ordered cover 
$\modtir$ by $\aut(\rho)\subset\Sigma_k$. 
The cohomology groups of the fibers of $\widehat\pi$ form a local system on $\modtor$. After pulling this local system back to $\modtir$, its fibers
are naturally identified with $\bigoplus_{j=1}^k H^\bullet(\patf_j),$ and the deck group $\aut(\rho)\cong\Sigma_k$ acts by permuting the summands.

Let $B$ be a sufficiently small stratified neighborhood of $\wV$, and let $V\in B\cap\modtor$. Since $\widehat\pi^{-1}(B)$ retracts onto the special fiber $\widehat\pi^{-1}(\wV)=\ParExtSp$,  restriction to the fiber over $V$ gives a map
\[H^\bullet(\ParExtSp) \simeq H^\bullet\bigl(\widehat\pi^{-1}(B)\bigr)
\longrightarrow H^\bullet\bigl(\widehat\pi^{-1}(V)\bigr),\]
whose image lies in the local monodromy invariants.  Since the pulled-back local system on the ordered cover is constant, the monodromy downstairs is given by the deck-group action. Thus the specialization map factors as 
\begin{equation}\label{eq:lastmap}
\Spmap\colon H^\bullet(\ParExtSp)
\longrightarrow
\left[
\bigoplus_{j=1}^k H^\bullet(\patf_j)
\right]^{\Sigma_k}\!.
\end{equation}
This map fits into the following commutative diagram:
\begin{equation}\label{eq:nonabdiagram}
	\begin{tikzcd}
		\hbull(F_{U^{\oplus k}}) \arrow[r,"\spmap"] \arrow[d,hook,"h^*"]
		&
		\hbull(F_V)^{\Sigma_k}
		\arrow[d,hook,"d_{\mathrm{MV}}"]
		\\
		\hbull(\ParExtSp) \arrow[r,"\Spmap"]
		&
		\left[\bigoplus_{j=1}^k \hbull(\patf_j)\right]^{\Sigma_k}\!.
	\end{tikzcd}
\end{equation}

\begin{lemma}\label{lem:Gmap-injective}
The map $\Spmap$ defined in \eqref{eq:lastmap} is injective.
\end{lemma}
\begin{proof}
We observe that composing $\Spmap$ with the projection onto the first summand gives a morphism
$$\hbull(\ParExtSp) \to  \hbull(\patf_1)^{\Sigma_{k-1}}, \text{\; where\; }
	\Sigma_{k-1}\cong\mathrm{Stab}_{\Sigma_k}(1).$$
Indeed,
\(\bigoplus_{j=1}^k \hbull(\patf_j)\)
	is the representation of \(\Sigma_k\) which is induced from the \(\Sigma_{k-1}\)-representation \(\hbull(\patf_1)\). Thus, by
	Frobenius reciprocity,
	\[
	\left[
	\bigoplus_{j=1}^k \hbull(\patf_j)
	\right]^{\Sigma_k}\simeq \hbull(\patf_1)^{\Sigma_{k-1}}.
	\]

The spaces \(\ParExtSp\) and \(\patf_1\) are projective bundles over
\(F_{U^{\oplus(k-1)}}\) and \(F_{V/V_1}\), respectively, and these
projective-bundle structures are compatible under specialization. By the
inductive hypothesis, the map
\[
\hbull(F_{U^{\oplus(k-1)}})
\longrightarrow
\hbull(F_{V/V_1})^{\Sigma_{k-1}}
\]
is injective. The projective bundle formula then shows that
\(
\hbull(\ParExtSp)
\longrightarrow
\hbull(\patf_1)^{\Sigma_{k-1}}
\)
is injective. Under the identification above, this is precisely the map $\Spmap$.
\end{proof}
It now follows from Lemma \ref{lem:Gmap-injective} and the commutative diagram \eqref{eq:nonabdiagram} that the map $\spmap$ is injective. This proves Theorem \ref{thm:nosupp} for \(\bro=[(\br,k)]\).

For a general point $\wV = \bigoplus_{a=1}^s U_a^{\oplus m_a}\!\in\! M_{\bro}$, where $\bro$ is a multipartition of $r$, we argue by induction on $m_1+\cdots+m_s$. Fix a block $U_a^{\oplus m_a}$, and let $\bro_{a}$ be obtained from $\bro$  by replacing $m_a$ with $m_a-1$.
The corresponding special and nearby spaces are projective bundles over the fibers associated with $\bro_a$, and these
projective-bundle structures are compatible under specialization, as in the proof of Lemma \ref{lem:Gmap-injective}. The stabilizer of the chosen copy is precisely $\Sigma[\bro_a].$ Hence Frobenius reciprocity, the inductive hypothesis, and the projective-bundle formula apply exactly as in the proof of Lemma \ref{lem:Gmap-injective}, giving the injectivity of $\spmap$ for an arbitrary multipartition, and thereby completing the proof of Theorem \ref{thm:nosupp}.

\section{Combinatorics of subtraction}\label{sec:combsub}

With the vanishing of all non-abelian supports established, we can now turn to 
the study of the stalks $L_\rho$ over the abelian strata.
Recall our simplified notation from \S\ref{sec:plan}.
	\begin{itemize}
\item We denoted the fiber 
$\pi^{-1}(V)$ at $V\in\modtor$   by  $F_V$; when we want to 
emphasize the  partition $\rho$ in the notation, we will write $F_\rho$ instead.
	\item Similarly, we will write $L_\rho=(\ls_\rho)_V$ for the stalk of the 
	local system $\ls_\rho$ supported on $\modtor$; we thus have $
	L_{\rho}\subset\hbull(F_{\rho})$.
	\item  We also introduce the notation $\slice_{\rho}$ for the normal
	slice  to $\modtor$ in  $\ms_0(r)$ at $V$  (cf. Proposition 
	\ref{prop:cone}). 
\end{itemize}
We can then write
\eqref{eq:cohsheavesrho}  as
\begin{equation}\label{eq:frhoeq}
	H^\bullet(F_\rho)=L_\rho\oplus \bigoplus_{\mu>\rho}
	IH^{\bullet}(\slice_{\mu\rho},\ls_\mu),
\end{equation}
where $\slice_{\mu\rho}$ denotes the normal slice to $\vect_\rho$ in $\vect_\mu$
at the point $V\in\modtor$.

In this section, we use the combinatorial models of \cite{ZhP} associated to
these varieties to calculate the
Hilbert function of the graded
vector space $ L_{\rho}\subset\hbull(F_{\rho}) $.
Our goal is to prove the following result.

\begin{theorem}\label{thm:HilbofL}
	Let $\rho= [r_1,\ldots ,r_k]\vdash r$.
	    The Hilbert function of $L_\rho$ is given by
	\begin{equation}\label{eq:HilbofL} \mathrm{Hilb}_t(L_\rho) =
	t^{2d(\rho)}\prod_{i=1}^k \frac{t^{2r_i}-1}{t^2-1},
	 \end{equation}
	 where $d(\rho) = \sum_{i<j}r_ir_j(g-1)$.
\end{theorem}

\subsection{The subtraction algorithm}\label{sec:subalg}

As a first step, we derive a recursion (the \textit{subtraction algorithm}) for
$\hilb_t(L_{\rho})$ based on \eqref{eq:frhoeq}. To describe the normal slice $\slice_{\mu\rho}$, we first introduce the following combinatorial device. 

\begin{definition}\label{def:decomposition}
A \textit{decomposition}
$\underline{\lambda}=\{\la_1,\ldots,
\la_n\}$ of $\underline{k}=\{1,\ldots, k\}$ is a representation 
of $\underline{k}$ as a disjoint union of subsets; we will use the notation
$\ula\vdash\uka$ in this case. Given $\underline{\lambda}=\{\la_1,\ldots, \la_n\}\vdash\uka$ and a partition $\rho=[r_1,\ldots,r_k]\vdash r$, we can compose $\ula$ and $\rho$  to obtain a
new partition $\ula\circ\rho=\mu=[m_1,\ldots,m_n]\vdash r$ such that
$m_j=\sum_{i\in \la_j}r_i$
for  $j=1,\ldots, n$; we will denote this partition of $m_{j}$ by
$\rho|\la_{j}$.
\end{definition}

\begin{remark}
	Note that if $\ula=\ukad$ is the decomposition of $\uka$ into single elements,
	then we have $\ukad\circ\rho=\rho$, while
$\ula\neq\ukad$ implies $\mu=\ula\circ\rho>\rho$, and any coarsening of $\rho$ may be obtained this way.
\end{remark}

For a fixed \(\mu>\rho\), the components of the normal slice
\(\slice_{\mu\rho}\) are indexed by decompositions \(\ula\) satisfying
\(\ula\circ\rho=\mu\). Its normalization may be written as
\[
\widetilde{\slice}_{\mu\rho}
\cong
\bigsqcup_{\ula:\,\ula\circ\rho=\mu}
\prod_{j=1}^n \slice_{\rho|\la_j}.
\]
Since the fibers \(F_\mu\) form a trivial bundle near \(V\), the local system
\(\mathcal L_\mu\) is trivial there, and hence
\[
IH^\bullet(\slice_{\mu\rho},\mathcal L_\mu)
=
\bigoplus_{\ula:\,\ula\circ\rho=\mu}
L_\mu\otimes
\bigotimes_{j=1}^n
IH^\bullet(\slice_{\rho|\la_j}).
\]
Then, taking the Hilbert function of the two sides of \eqref{eq:frhoeq}, we obtain
$$ P_t(F_\rho)= \mathrm{Hilb}_t\left(L_\rho\right)+
\sum_{\mu>\rho}\sum_{\ula:\,\mu=\ula\circ\rho}
\mathrm{Hilb}_t\left(L_\mu\right)\cdot\prod_{j=1}^n
IP_t(\slice_{\rho|\la_j}).$$

Finally, eliminating $\mu$, we arrive at the following formula, which forms
the basis of our
\textit{subtraction algorithm}.
\begin{proposition}\label{prop:subtract}
Given an abelian point $V\in\modtor$, the Hilbert function of the stalk of the
local system $\mathcal{L}_\rho$ at $V$ may be expressed via the following recursive formula:
\begin{equation}\label{eq:subtralg}
	 \mathrm{Hilb}_t\left(L_\rho\right)=
P_t(F_\rho)-	\sum_{\ula\vdash\uka,\ \ula\ne\ukad}
	\mathrm{Hilb}_t\left(L_{\ula\circ\rho}\right)\cdot\prod_{j}
	IP_t(\slice_{\rho|\la_j})	.
\end{equation}
\end{proposition}
This recursion requires two numerical inputs: the Poincar\'e polynomial of the
abelian fiber $F_\rho$ and the intersection Poincar\'e polynomial of the
normal slice $\slice_\rho$. Below, we will use the type-A toric graph calculus
of \cite{ZhP} to represent these quantities combinatorially.

\subsection{Graphs and intersection cohomology}\label{sec:defFrho} We start by recalling
the setup and some of the results of \cite{ZhP}.
Let $G$ be an oriented graph on a finite vertex set $V(G)$, possibly
with multiple edges, but without loops; we denote the set of edges of $G$ by $E(G)$, and their number  by  $e(G)$. Given a subset $J\subseteq V(G)$, we denote by $G|J$ the subgraph of $G$ induced by the vertices in $J$.
We say that $G$ is \textit{rooted} in $v\in V(G)$ if every vertex
has an oriented path to $v$, and $G$ is \textit{acyclic} if it has no oriented cycles.  For an oriented graph $\Gamma$, set
\[
\acyc[\Gamma,v,n]=
\{G\subset\Gamma\mid V(G)=V(\Gamma),\ G\text{ is acyclic and rooted in }v,\ e(G)=n\}.
\]
\noindent Note that the edges of $\Gamma$ emanating from the root $v$ cannot
belong to any such subgraph $G$.
For a partition $\rho=[r_1,\ldots,r_k]\vdash r$, we denote by 
\begin{itemize}
	\item  $\bugf$ the directed graph with vertex set
	$\{1,\ldots,k\}$, and $r_ir_j(g-1)$ oriented edges from $i$ to $j$ for
	$0<i\neq j\leq k$;
	\item $\augf$ the augmented graph on the vertex set
	$\{0,1,\ldots,k\}$, such that $\augf|\uka=\bugf$, and there are $r_i$ edges
	between $0$ and $i$ in either direction for $i\in\uka$.
\end{itemize}

In the notation of \cite{ZhP}, the affine toric variety $\slice_\rho$ (cf. Proposition \ref{prop:cone})
is the type-A toric variety associated with the graph $\bugf$. Its intersection Poincar\'e polynomial admits the following combinatorial description.\footnote{ In \cite[Theorem 4.5]{ZhP}, the
formula is stated with the first vertex as root. However, the
root-changing construction \cite[Construction 5.6]{ZhP} gives a
bijection between acyclic subgraphs rooted at the first vertex and those
rooted at any prescribed vertex, preserving the number of edges. Thus
the same formula holds for an arbitrary choice of the root vertex.}

\begin{theorem}[{\cite[Theorem 4.5]{ZhP}}]\label{thm:zhep}
Let $\rho=[r_1,\ldots,r_k]\vdash r$, and set $g(\rho;t^2)=IP_t(\slice_\rho)$.  For
any root vertex $v\in\{1,\ldots,k\}$,
\begin{equation}\label{eq:gpoly}
 g(\rho;u+1)=\sum_{m=0}^{d(\rho)-k+1} \left|\acyc[\bugf,v,k-1+m]\right|u^m.
\end{equation}
\end{theorem}
\begin{remark}\label{rem:indeproot} 
	In particular, the right-hand side of \eqref{eq:gpoly} is independent of the choice of root $v$.
\end{remark}
\subsection{Poincar\'e polynomials of the fibers in terms of graphs}
We have the following graph-theoretic presentation for the Poincar\'e 
polynomial of the fiber $F_\rho$, which rhymes with Theorem \ref{thm:zhep}.
		\begin{theorem}\label{thm:fibers}
			Let $\rho=[r_1,\ldots,r_k]\vdash r$ be a partition of $r$,  and recall the notation $f(\rho;
			t^2)=P_t(F_\rho)$ from
			Theorem \ref{thm:PoinFib}. Then
		\begin{equation}\label{eq:frec-graph}
				f(\rho; u+1) =
			\sum_{m=0}^{d(\rho)+r-k}
			\left|\acyc[\augf,0,k+m]\right|\cdot
			u^{m}.
		\end{equation}
		\end{theorem}
\begin{proof}
We prove the theorem using induction on $k$. A simple calculation shows that the statement is true for $k = 1$.

It follows from \cite[Theorem 4.5]{ZhP} that the right-hand side of \eqref{eq:frec-graph} is the intersection Poincar\'e polynomial of the type-A affine toric variety determined by the vertex set
 \(\{0,1,\ldots,k\}\) and the multiplicities
	\[
	r_{0i}=r_{i0}=r_i,\qquad
	r_{ij}=r_ir_j(g-1)\quad (i,j>0,\ i\ne j).
	\]
Theorem 4.19 of \cite{ZhP} provides a recursion for the intersection Poincar\'e polynomials of these type-A affine toric varieties. After substituting the above multiplicities, this recursion agrees exactly with the recursion \eqref{eq:frec-JT} satisfied by $f(\rho;u+1)$.
 Moreover, the two recursions have the same initial value for $k=1$. The inductive hypothesis therefore yields \eqref{eq:frec-graph}.
\end{proof}
\begin{remark}
In particular, the Betti numbers of $F_\rho$ coincide with the intersection
cohomology Betti numbers of an
affine toric variety.
\end{remark}

Now we are ready to prove the main combinatorial result of this section.

		\begin{proposition}\label{prop:fibers-zhp}
			 Let $\rho=[r_1,\ldots,r_k]\vdash r$
			be a partition of $r$. Recall the notation introduced above:
			for $\ula=\{\la_{1},\dots,\la_{n}\}\vdash\uka$, we denote by
			 $\ula\circ\rho$ the corresponding composed partition of $r$. We set $r_{\la_j}=\sum_{i\in \la_j}r_i$ and denote by $\rho|\la_j$ the partition of $r_{\la_j}$ induced by $\rho$.
             Finally, we set $d(\ula\circ\rho) = \sum_{i<j}r_{\la_i}r_{\la_j}(g-1) $.
			Then
			\begin{equation}\label{eq:substrgraphst}
				f(\rho; t)= \sum_{\underline{\lambda} \vdash \underline{k}}
				t^{d(\ula\circ\rho)}
				\prod_{i}g(\rho|{\la_i};t)\frac{t^{r_{\la_i}}-1}{t-1},
			\end{equation}
			where the sum is taken over all decompositions
			$\underline{\lambda}=\{\la_{1},\dots,\la_{n}\}$ of the set
			$\underline{k}=\{1,2,\ldots,k\}$.
		\end{proposition}

\begin{proof}
Substituting $t=u+1$   into \eqref{eq:substrgraphst}, we obtain:
\begin{equation}\label{eq:subsu}
	f(\rho;u+1)=
\sum_{\ula\vdash\uka}
(u+1)^{d(\ula\circ\rho)}
\cdot\prod_i g(\rho|\la_i;u+1)\frac{(u+1)^{r_{\la_i}}-1}{u}.
\end{equation}

The proof rests on a weight-preserving bijection between the graph-theoretic models on the two sides of \eqref{eq:subsu}. We first interpret each side as a weighted count and then construct the bijection explicitly.

According to \eqref{eq:frec-graph}, the \textit{left-hand side} counts the acyclic subgraphs $G\subset \augf$ on the vertex set $\{0, 1,\ldots,k\}$,  rooted in $0$, with weight $u^{e(G)-k}$.

To interpret the \textit{right-hand side}, we use the elementary identities
	\[
	(u+1)^m=\sum_{S\subset\underline{m}}u^{|S|} \text{ \; and \; } 	\frac{(u+1)^m-1}{u}=\sum_{\emptyset\ne S\subset\underline{m}}u^{|S|-1}.
	\]
Together with Theorem \ref{thm:zhep} and Remark \ref{rem:indeproot}, these identities show that the right-hand side of \eqref{eq:subsu} counts the following choices: 
\begin{enumerate}
	\item Choose a decomposition \(\ula\vdash \uka\) with $n$ blocks and, for each block, a nonempty set of ``red'' edges from its vertices to $0$. 
\item Label the blocks and the corresponding sets of red edges by $\la_1,\ldots,\la_n$ and $R_1,\ldots,R_n$ so that, with $v_i= \min\{v\in\lambda_i\mid v\rag 0\in R_i\}$, one has $ v_1<\cdots<v_n.$
	\item For each \(i=1,\ldots,n\), choose an acyclic subgraph $G_i\in \acyc[\bugf|\lambda_i,v_i,|\lambda_i|-1+m_i]$ for some integer \(m_i\geq0\).
	\item  Finally, choose an arbitrary subset $E$ of the oriented edges of $\bugf$ which point from a vertex of $\lambda_i$ to a vertex of $\lambda_j$ with $i<j$. 
The number of  possible inter-block edges  is \(d(\ula\circ\rho)\); since the multiplicity of the edges between two blocks is symmetric, this number is independent of the ordering of the blocks. 
\end{enumerate}
Denoting by underlined symbols the objects depending on \(i=1,\ldots,n\), the resulting tuple $(\ula,\underline G,\underline R,E)$ has weight $ u^{|E|-k+\sum_{i=1}^n |R_i| +\sum_{i=1}^n e(G_i)}.$

We now construct the bijection underlying equality \eqref{eq:subsu}. To each acyclic subgraph $G\subset \augf$ rooted at $0$, we associate a tuple $(\ula, \underline{G}, \underline{R}, E)$ satisfying the conditions above.
\begin{notation}
  We write $v \rag w$ if there is an edge $e$ directed from $v$ to $w$ (allowing multiple such edges), and  $v \tap w$ if there exists an oriented path in $G$ from $v$ to $w$. By convention, $v \tap v$ for every vertex $v$.
\end{notation} \vspace{-15pt}
\begin{itemize}
	\item Set $R=\{v\rag0\in E(\Gr)\}$ to be the set of red edges.
	\item Let   $v_1=\min \{v\in \uka|\,v\rag0\in R\}$,
	$\lambda_1=\{w\in\uka|\,w\tap v_1\}$,  and set $G_1=\Gr|\lambda_1$, $R_1=\{v\rag0\in R|\,v\in \lambda_1\}$.
         \item  Having constructed \(\lambda_1,\ldots,\lambda_{i-1}\), set $v_i=\min\{v\in\uka\setminus\bigcup_{j<i}\lambda_j|\,v\rag 0\in R\},$ and $\lambda_i=\{w\in\uka\setminus\bigcup_{j<i}\lambda_j\; |\: w\tap v_i\}$, together with $G_i=\Gr|{\la_i},$ $R_i=\{v\rag 0\in R|\,v\in\la_i\}.$ We note that the graph $G_i$ is acyclic and rooted at $v_i$.
\item We continue until the process terminates. Since \(\Gr\) is rooted at  \(0\), the construction exhausts \(\uka\), and hence $\uka=\lambda_1\sqcup\cdots\sqcup\lambda_n.$ 
	\item Finally, set $E=\{v\rag w\in E(\Gr)|\, v\in\lambda_i,w\in\lambda_j,i\neq j\}$. By construction, whenever $v\rag w$ belongs to $E$, with $v\in\lambda_i$ and $w\in\lambda_j$, we  have $i<j$.
\end{itemize} 

Conversely,  given a tuple $(\ula,\underline G,\underline R,E)$ as above, we recover  $\Gr$ by taking the union of the internal graphs $G_i$, the red edges $R_i$, and the inter-block edges $E$. The ordering condition on $E$ and the definition of the roots $v_i$ from the sets $R_i$ ensure that the resulting graph is acyclic and rooted at $0$, and that the successive construction above recovers the same blocks $\lambda_i$. Thus the two constructions are inverse. 

 Finally, \[ e(\Gr)-k = |E| +\sum_{i=1}^n e(G_i) +\sum_{i=1}^n |R_i| -k, \] so the bijection preserves weights. This completes the proof of Proposition \ref{prop:fibers-zhp}. \end{proof}

We are now ready to present the proof of Theorem \ref{thm:HilbofL}.
\begin{proof}
Substituting  $g(\rho|\la_j;t^2)=IP_t(\slice_{\rho|\la_j})$ and $f(\rho; t^2)=P_t(F_{\rho})$ (see Theorems \ref{thm:zhep} and \ref{thm:fibers}) into
\eqref{eq:subtralg}, we obtain the identity
\begin{equation}\label{eq:recursionHilb}
	\mathrm{Hilb}_t\left(L_\rho\right) = f(\rho;t^2) -
\sum_{\ula\vdash\uka,\ \ula\ne\ukad} \mathrm{Hilb}_t\left(L_{\ula\circ\rho}\right)\cdot\prod_{j}
	g(\rho|\la_{j};t^2).
\end{equation}
This formula gives a recursion for the Hilbert series $\mathrm{Hilb}_t(L_\rho)$ with respect to the refinement ordering of partitions $\rho \vdash r$. The base of the recursion is the trivial partition $\rho = [r]$, where
$\mathrm{Hilb}_t\left(L_{[r]}\right) =   P_t(F_{[r]})$ is the Poincar\'e
polynomial of the generic fiber of the map $\pi$.  Since this generic fiber is isomorphic to $\mathbb{P}^{r-1}$, we have $\mathrm{Hilb}_t\left(L_{[r]}\right) =\frac{t^{2r}-1}{t^2-1}$.

On the other hand, equation \eqref{eq:substrgraphst} can be rewritten as
\[
t^{2d(\rho)}
\prod_{i=1}^k
\frac{t^{2r_i}-1}{t^2-1} = f(\rho;t^2) -\!
\sum_{\ula\vdash\uka,\ \ula\ne\ukad}\!
t^{2d(\ula\circ\rho)}
\prod_j g(\rho|\la_j;t^2) \frac{t^{2r_{\la_j}}-1}{t^2-1}.
\]
Comparing this identity directly with the recursive formula \eqref{eq:recursionHilb}, we conclude that the solution to \eqref{eq:recursionHilb} is uniquely given by \eqref{eq:HilbofL}. This completes the proof of Theorem \ref{thm:HilbofL}.
\end{proof}

		\section{Refined intersection forms and the calculation of the local
		systems}\label{sec:locsys}

		In this section, we describe the local system $\ls_{\rho}$
		appearing	in \eqref{eq:dtforrho} by identifying its stalk
		$L_\rho$ (cf. \eqref{eq:cohsheavesrho}) as a subspace of
		$\hbull(F_\rho)$; see Theorem \ref{thm:corL} for the result.

		As a first step, we recast this problem in
		the framework of the calculation of certain enumerative quantities
		called \textit{refined intersection forms}. Let us briefly recall their definition in
		a somewhat more general context.

		\subsection{Refined intersection forms}\label{sec:refint}
		Let $\iota\colon Z\hookrightarrow X$	be a projective, possibly singular
		subvariety in a
		smooth, not necessarily compact variety $X$ of dimension  $m \le2\dim Z$.
		Let $S,T\subset Z$ be cycles representing classes
		$$[S]\in H_{n}(Z) \text{\; and \;} [T]\in H_{2m-n}(Z),$$ respectively.

Poincar\'e duality induces a pairing 
		$$\langle \cdot, \cdot\rangle_X\colon H_{n}({X}) \otimes 
		H^{\mathrm{BM}}_{2m-n}(X)\to\mathbb Q.$$  Let $\Psi\colon 
		H_{\bullet}({X}) \to H_{\bullet}^{\mathrm{BM}}(X)$ be the natural map from 
		ordinary to Borel--Moore homology. 
 The \emph{refined intersection form} on $H_{n}(Z)\otimes H_{2m-n}(Z)$ is then defined by
\begin{equation}\label{eq:refinedint}
[S]\cdot_{X}[T] := \langle\iota_{*}[S],
			\Psi(\iota_{*}[T])\rangle_{X}\in\Q.
		\end{equation}
Equivalently, this pairing is encoded by  an adjoint map
$$\Phi\colon H_n(Z)\to H^{2m-n}(Z), \text{\; such that \;}
 [S]\cdot_{X}[T] = \langle \Phi[S],[T]\rangle. $$
We now relate this map to the long exact sequence in
		Borel--Moore homology:
\begin{equation}\label{eq:les-BM}
			\cdots \to H^{\mathrm{BM}}_n(Z)  \to H^{\mathrm{BM}}_n(X) \to
			H^{\mathrm{BM}}_n(X\setminus Z)\to\cdots.
		\end{equation}

Since $Z$ is projective, we have $H_n^{\mathrm{BM}}(Z)=H_n(Z)$, and Poincar\'e 
duality for the smooth varieties $X$ and $X\setminus Z$ gives 
$$H_n^{\mathrm{BM}}(X)\cong H^{2m-n}(X) \text{\; and \;} 
H_n^{\mathrm{BM}}(X\setminus Z)\cong H^{2m-n}(X\setminus Z).$$
Thus the sequence \eqref{eq:les-BM} may be rewritten as
\begin{equation}\label{eq:les-X}
			\cdots \to H_n(Z) \to H^{2m-n}(X)
			\to H^{2m-n}(X\setminus Z) \to\cdots.
\end{equation}
The following statement is straightforward.

\begin{lemma}\label{lem:proprefL} Assume that, in addition, 
			$Z$ is a deformation retract of \(X\), and let 
		 $\jmap\colon X\setminus Z\hookrightarrow X$ be the inclusion. 
		Using the identification 
			$H^\bullet(X)\xrightarrow{\sim}H^\bullet(Z)$, set
			\[
L^\bullet=\ker\left(H^\bullet(Z)\xrightarrow{\jmap^*}
			H^\bullet(X\setminus Z)\right).
			\]
			Then the sequence \eqref{eq:les-X} can be written as
			\begin{equation}\label{eq:les-Z}
				\cdots \to H_n(Z) \xrightarrow{\ \Phi\ } H^{2m-n}(Z)  
				\overset{\jmap^{*}}\longrightarrow 		H^{2m-n}(X\setminus Z) 
				\to\cdots,
			\end{equation}
			where \(\Phi\colon H_n(Z)\to H^{2m-n}(Z)\) is the  adjoint map 
			defined above. In particular, for each $n$,  we have 
			\(
			L^{2m-n}=\operatorname{im}\Phi|_{H_n(Z)}.
			\)
		\end{lemma}
\begin{corollary}\label{cor:refform} The refined intersection form can be 
represented by an element of $H^\bullet(Z)\otimes H^\bullet(Z)$, and   every 
cohomology class obtained by contracting this element belongs to  $L^\bullet$. 
In particular, if the refined intersection form is represented by a 
decomposable tensor $\alpha\otimes\beta$ with $\alpha, \beta\neq 0$, then we 
have $ \alpha,\beta\in L^\bullet. $
\end{corollary}

	\subsection{The  refined intersection form for abelian points}
We now explain how the refined intersection forms introduced above are used to compute  $L_\rho\subset \hbull(F_\rho)$.
	
Let $V\in \modtor$ be a generic point, and recall the notation $\slice_{\rho}$ for the normal slice to $\modtor$ in  $\ms_0(r)$ at $V$. According to the Decomposition Theorem \ref{thm:DT},  $L_\rho$ can be obtained as the kernel of a cohomological restriction map \eqref{eq:cohrestr}.
Let $U\simeq \slice_\rho\times D$ 
be a sufficiently small product neighborhood of \(V\), where \(D\) is a
small disk in the stratum \(\modtor\). Then the pair
$$(\pi^{-1}(U),\pi^{-1}(U\setminus\modtor)) \text{\; retracts to \;} 
(\pi^{-1}(\slice_{\rho}),\pi^{-1}(\slice_{\rho}\setminus\{V\})),$$ while 
\(F_\rho\) is a deformation retract of \(\pi^{-1}(\slice_\rho)\). Therefore 
\(L_\rho\) may be identified with the
kernel of the map
\begin{equation}\label{eq:jmap}
\jmap^*\colon \hbull(F_\rho) \cong\hbull(\pi^{-1}(\slice_\rho))
\longrightarrow
\hbull(\pi^{-1}(\slice_\rho\setminus\{V\})).
\end{equation}

To relate this description of \(L_\rho\) to the refined intersection form, we use the same local product neighborhood. The projection \(U\simeq\slice_\rho\times D\to D\) induces a submersion $\pi^{-1}(U)\to D$. Since \(\pi^{-1}(U)\) is open in the smooth variety \(\pms_0(r)\), the fiber $X_\rho:=\pi^{-1}(\slice_\rho\times\{0\})$ is smooth. Therefore, by the discussion in \S\ref{sec:refint}, the pair $(X_\rho, Z_\rho:=F_\rho)$ carries a refined intersection form, which may be represented by a class in $\hbull(F_\rho)\otimes\hbull(F_\rho).$

 The following theorem gives an explicit description of this form in the degrees relevant to our computation.
\begin{theorem}\label{thm:intform}
Let $\rho=[r_1,\ldots,r_k]\vdash r$ be a partition, and recall the notation $d(\rho) = \sum_{i<j}r_ir_j(g-1)$. Let $x_1,\ldots,x_k\in H^2(F_\rho)$ be the cohomology classes defined in Proposition \ref{prop:xsigma} and set
$$w_\rho:=\prod_{i<j}(x_i-x_j)^{r_ir_j(g-1)} \text{ \; and \;} w_\rho^+:=(-1)^{d(\rho)} w_\rho\prod_{i=1}^k x_i^{r_i-1}.$$
Then the refined intersection form on $H_{2d(\rho)}(F_\rho)\otimes H_{2(d(\rho)+r-k)}(F_\rho)$ is  represented by $w_\rho\otimes w_\rho^+.$
		\end{theorem}

Combining \eqref{eq:jmap} with Corollary \ref{cor:refform}, we see that the classes $w_\rho$ and $w_\rho^+$ described in Theorem \ref{thm:intform} give explicit elements of \(L_\rho\) in their respective degrees noted above.

			\subsection{Excess intersection calculations}
			Before turning to the proof, we make a brief detour
			through the excess intersection formula. Its role here is geometric: the
			refined form is computed by intersecting cycles in $F_\rho$ inside the
			smooth ambient variety $X_\rho$. Although $F_\rho$ is singular, the
			relevant homology classes may be represented on its smooth
			projective-tower components, so the calculation reduces to Chern classes
			of excess bundles.
\begin{proposition}[Excess intersection formula, {\cite[\S 13.3]{EH16}}]\label{thm:excess}
Let \(S,T\subset X\) be smooth projective subvarieties of a smooth variety 
\(X\), meeting cleanly along the smooth subvariety $K:=S\cap T,$ $\iota_K\colon 
K\hookrightarrow X$.  Define the quotient bundle $E_{S,T} 
:=\bigl(N_{S/X}\big|_K\bigr)\big/\bigl(N_{K/T}\bigr)$ and set  
$e:=\operatorname{rk}E_{S,T}.$
	\begin{itemize}
	\item In Borel--Moore homology, the intersection class \([S]\cdot_X[T]\) is the pushforward to \(X\) of the Poincar\'e dual on \(K\) of the top Chern class of the excess bundle:
\begin{equation}\label{eq:excess}
[S]\cdot_X[T] = (\iota_K)_*\operatorname{PD}_K\left(c_e(E_{S,T})\right).
\end{equation}
\item If \(S\) and \(T\) have complementary dimensions in \(X\), then \(e=\dim K\), and hence
$$[S]\cdot_X[T]=\int_K c_e(E_{S,T}).$$
	
\item If \(S,T\subset Z\subset X\), where \(Z\) is smooth, and \(S\) and \(T\) meet transversely in \(Z\), then
	\[
	E_{S,T}\cong N_{Z/X}|_K.
	\]
\item In particular, if \(T=Z\) and \(S\) and \(Z\) have complementary dimensions in \(X\), then 
	\begin{equation}\label{eq:refsmooth1}
		[S]\cdot_X[Z]=\int_S c_{\dim S}(N_{Z/X}|_S).
	\end{equation}
	\end{itemize}
\end{proposition}
We now proceed to the proof of Theorem \ref{thm:intform} using this excess intersection formula. We begin with several preparatory results.

Given a permutation $\sigma\in\Sigma_k$, recall the notation $\pat_\sigma$ for 
the tower of projective bundles introduced in \S\ref{sec:inttowers}. It has 
dimension $d(\rho)+r-k$ and is an irreducible component of the fiber $F_\rho$ 
(see Proposition \ref{prop:intuni}). In what follows, we write $i\prec_\sigma j 
$ whenever  $\sigma^{-1}(i)<\sigma^{-1}(j),$ and set $$l(\sigma)= 
\sum_{\substack{1\le i<j\le k\\ j\prec_\sigma i}}r_i r_j(g-1).$$
		\begin{lemma}\label{lem:normal-tower} The top Chern class of the normal bundle of $\pat_\sigma$ in $X_\rho$ is
$$	c_{d(\rho)}(N_{\pat_\sigma/X_\rho})=(-1)^{l(\sigma)}w_\rho.$$
		\end{lemma}
		\begin{proof}
A simple calculation, similar to that of Proposition \ref{prop:ringtower}, shows that the Chern character of the tangent bundle of $X_\rho$, the inverse image of the normal slice, is equal to
			\[
			\ch(TX_\rho\big|_{\pat_\sigma})=
			\sum_{1\leq i\leq k} r_i\exp(x_i)+ \sum_{1\le i\ne j\le k}r_ir_j(g-1)\exp(x_i-x_j)-k.
			\]
			By Proposition \ref{prop:chTP}, we have
			\[ \ch(T\pat_\sigma)= \sum_{1\leq i\leq k} r_i\exp(x_i)+ \sum_{1\le i<j\le k}r_{\sigma(i)}r_{\sigma(j)}(g-1) \exp(x_{\sigma(j)}-x_{\sigma(i)})-k,	\]
and thus
			\begin{equation}\label{eq:chPPX}
				\ch(N_{\pat_\sigma/X_\rho})= (g-1)\sum_{1\le i<j\le k}r_{\sigma(i)}r_{\sigma(j)} \exp(x_{\sigma(i)}-x_{\sigma(j)}).
			\end{equation}
			Taking the top Chern class, we arrive at the stated result.
		\end{proof}
	Recall the notation $\ipat(\sigma,\tau)=\pat_\sigma\cap\pat_\tau$ for the 
	intersection of the two towers corresponding to the permutations	
	$\sigma, 
	\tau\in \Sigma_k$.
		We will also  write $\ipat_\square\cong\prod_{i=1}^k\PP^{r_i-1}$,  the 
		common intersection of all towers in $F_\rho$ (see Proposition 
		\ref{prop:intuni}).
				\begin{lemma}\label{lem:pd}
Let $w_{\sigma,\tau}$ be the top Chern class of the quotient  bundle
$$\bigl(N_{\pat_\sigma/X_\rho}\big|_{\ipat(\sigma,\tau)}\bigr)\big/\bigl(N_{\ipat(\sigma,\tau)/\pat_\tau}\bigr)$$
	on $\ipat(\sigma,\tau)$.  Then the Poincar\'e dual of $w_{\sigma,\tau}$  in	$\ipat(\sigma,\tau)$ is
			\((-1)^{l(\sigma)+l(\tau)+d(\rho)}[\ipat_\square].\)
		\end{lemma}
\begin{proof}
It follows from  Proposition \ref{prop:chTInt} that
			\[
			\ch(T\ipat(\sigma,\tau))-\ch(T\ipat_\square)=
			\sum_{\substack{i\prec_\tau j\\ i\prec_\sigma j}}
			r_ir_j(g-1)\exp(x_j-x_i),
			\]
hence $\ipat_\square$ is Poincar\'e dual  in $\ipat(\sigma,\tau)$  to
			\[
			\prod_{\substack{i\prec_\tau j\\ i\prec_\sigma j}}
			(x_j-x_i)^{r_ir_j(g-1)}.
			\]
			On the other hand,  Propositions \ref{prop:chTP} and  \ref{prop:chTInt}, together with  equation
			\eqref{eq:chPPX}, give
			\[
			w_{\sigma,\tau}=\prod_{\substack{j\prec_\tau i\\ j\prec_\sigma i}}
			(x_j-x_i)^{r_ir_j(g-1)}.
			\]
			Comparing the two products gives the sign
			$(-1)^{l(\sigma)+l(\tau)+d(\rho)}$, and the result follows.
		\end{proof}

		\begin{proof}[Proof of Theorem \ref{thm:intform}]
Since $F_\rho$ has dimension $d(\rho)+r-k$, and the towers  $\pat_\sigma$ are its irreducible components, their fundamental classes generate $H_{2(d(\rho)+r-k)}(F_\rho)$.  On the other hand, applying \eqref{eq:towerinject} to the partition $\rho$ and dualizing
its degree-$2d(\rho)$ part gives a surjection
\[ \bigoplus_{\tau\in\Sigma_k}H_{2d(\rho)}(\pat_\tau) \twoheadrightarrow H_{2d(\rho)}(F_\rho).
\] 
Hence, by bilinearity, it is enough to evaluate the refined intersection
form on pairs consisting of a $d(\rho)$-dimensional cycle $S\subset\pat_\tau$ and a tower $\pat_\sigma$.

If $\sigma=\tau$, putting together equation \eqref{eq:refsmooth1},  Lemma \ref{lem:normal-tower} and Corollary \ref{cor:towerintegral}, we obtain
			\[
			[S]\cdot_{X_\rho}[\pat_\sigma] \overset{\eqref{eq:refsmooth1}}{=} \int_{S} c_{d(\rho)}(N_{\pat_\sigma/ X_\rho}|_S) \overset{\ref{lem:normal-tower}}{=}  (-1)^{l(\sigma)}\int_S w_\rho \overset{\ref{cor:towerintegral}}{=} \int_{\pat_\sigma}w_\rho^+\cdot \int_Sw_\rho.
			\]
For general $\sigma$ and $\tau$, the excess intersection formula \eqref{eq:excess} shows that
the intersection of $\pat_\sigma$ and $\pat_\tau$ in $X_\rho$ is the pushforward of the Poincar\'e dual of the class $w_{\sigma,\tau}$ defined in Lemma \ref{lem:pd}. Hence Lemma \ref{lem:pd} gives 
$$[\pat_\tau]\cdot_{X_\rho}[\pat_\sigma] = (-1)^{l(\sigma)+l(\tau)+d(\rho)}[\ipat_\square].$$ 
Since \(S\subset\pat_\tau\), compatibility of intersection products gives
\[
[S]\cdot_{X_\rho}[\pat_\sigma] = [S]\cdot_{\pat_\tau} \left([\pat_\tau]\cdot_{X_\rho}[\pat_\sigma]\right)=(-1)^{l(\sigma)+l(\tau)+d(\rho)}
[S]\cdot_{\pat_\tau}[\ipat_\square].
\]
A calculation analogous to that in the proof of Lemma \ref{lem:normal-tower}, together with Proposition \ref{prop:chTP}, shows that the Poincar\'e dual of $\ipat_\square$ in $\pat_\tau$ is
$(-1)^{d(\rho)+l(\tau)}w_\rho.$
We thus obtain
$$
[S]\cdot_{X_\rho}[\pat_\sigma]= (-1)^{l(\sigma)} \int_S w_\rho= \int_S w_\rho\cdot \int_{\pat_\sigma}w_\rho^+, $$
where the last equality follows from Corollary \ref{cor:towerintegral}.
\end{proof}

		\subsection{The local system on an abelian stratum}
		We now identify the stalk of the local system $\ls_\rho$ with  a 
		subspace of $\hbull(F_\rho)$.

		\begin{theorem}\label{thm:corL}
 Let $\rho=[r_1,\ldots,r_k]\vdash r$, and let $x_1,\ldots,x_k\in H^2(F_\rho)$ 
 be the cohomology classes defined in Proposition \ref{prop:xsigma}.  Then 
 $L_\rho\subset H^\bullet(F_\rho)$ is the ideal  generated by $w_\rho = 
 \prod_{i<j}(x_i-x_j)^{r_ir_j(g-1)}\in\hbull(F_\rho)$, and has the following 
 basis:
\begin{equation}\label{eq:deflmu}
				L_\rho= \bigoplus_{\substack{0\leq m_i<r_i \\  i=1,\ldots,k}}
\Q\,w_\rho x_1^{m_1}\cdots x_k^{m_k}.
			\end{equation}
		\end{theorem}
		\begin{proof}
It follows from  Corollary \ref{cor:refform} and Theorem \ref{thm:intform} that 
the element $w_\rho$ belongs to $L_\rho$. By the Decomposition Theorem 
\ref{thm:DT}, $L_\rho$ is the kernel of the map \eqref{eq:cohrestr}, and is 
therefore an ideal in $H^\bullet(F_\rho)$; it thus contains the right-hand side 
of 
\eqref{eq:deflmu}.

To prove the linear independence of the classes in \eqref{eq:deflmu}, suppose 
that
			\[
			\sum_{\substack{0\le a_i<r_i\\ a_1+\cdots+a_k=s}}
			c_{a_1,\ldots,a_k}\,
			w_\rho x_1^{a_1}\cdots x_k^{a_k}=0
			\]
is a homogeneous relation. 	Fix a multi-index $(a_1,\ldots,a_k)$ appearing in 
the sum, multiply the relation  by
$\prod_{i=1}^k x_i^{r_i-1-a_i}, $ and integrate over any tower 
$\pat_\sigma$.  By Corollary \ref{cor:towerintegral}, all terms in the integral 
vanish except the one corresponding to $(a_1,\ldots,a_k)$, whose integral is 
nonzero.  Hence $c_{a_1,\ldots,a_k}=0$, and since the multi-index was	
arbitrary, 
all	coefficients vanish.

Finally, the Hilbert series of the right-hand side of \eqref{eq:deflmu} is 
obviously
\[ t^{2d(\rho)}\prod_{i=1}^k\frac{t^{2r_i}-1}{t^2-1}, \]
which, by Theorem \ref{thm:HilbofL}, is exactly the Hilbert series of $L_\rho$. 
Therefore the classes in \eqref{eq:deflmu} form a basis of $L_\rho$, and the 
result follows.
\end{proof}

\section{Application of the global Decomposition Theorem}
\label{sec:final}

We begin by summarizing the results obtained so far. For the
map $\pi\colon\pms_0(r)\to\ms_0(r)$, the Decomposition Theorem takes the
form
\begin{equation}\label{eq:sheafwow}
	R\pi_*\Q_{\pms_0(r)}
	\cong
	\bigoplus_{\rho\vdash r}
	IC\left(\overline{\modtor},\ls_\rho\right),
\end{equation}
where the local system $\ls_\rho$ is defined on the abelian
stratum $\modtor$, and its stalk is the vector space $L_\rho$ described in Theorem \ref{thm:corL}. In this final section, we convert \eqref{eq:sheafwow} into the recursive formula announced in Theorem \ref{thm:main}.

Our first step is to understand the monodromy of each local
system $\ls_\rho$.
Let \(\rho=[r_1,\ldots,r_k]\vdash r\), and recall from \eqref{eq:gammarho} the finite, \(\aut(\rho)\)-invariant morphism
\[
\gamma_\rho\colon {\Sprod_\rho}
\longrightarrow \overline{\modtor}.
\]

\begin{lemma}\label{lem:icohrho}
	We have the natural identification
	\begin{equation}\label{eq:icohrho}
		IH^\bullet(\overline{\modtor},\ls_\rho)
		\cong
		\left[
		IH^\bullet(\Sprod_\rho)\otimes L_\rho
		\right]^{\aut(\rho)}.
	\end{equation}
\end{lemma}

\begin{proof}
Over the abelian stratum, $\gamma_\rho$ restricts to the finite cover
\[
\widetilde{\gamma}_\rho\colon\modtir\longrightarrow\modtor
\]
with deck group \(\aut(\rho)\). The local system \(\ls_\rho\) is obtained from this cover and the $\aut(\rho)$-representation $L_\rho$. In our description of $L_\rho$ (see Theorem \ref{thm:corL}),  \(\aut(\rho)\) acts by permuting the variables $x_1,\ldots,x_k$ corresponding to equal parts of $\rho$.

Since $\gamma_\rho$ is finite, it is a small proper morphism. Hence,  $R\gamma_{\rho*}IC_{\Sprod_\rho}$ is the intermediate extension of its restriction to the abelian stratum. It follows that, for the diagonal action of $\aut(\rho)$,
$$IC(\overline{\modtor},\ls_\rho) \cong \left( R\gamma_{\rho*}IC_{\Sprod_\rho}\otimes L_\rho \right)^{\aut(\rho)}.$$
Taking hypercohomology and using the exactness of taking invariants over
$\Q$, we obtain \eqref{eq:icohrho}.
\end{proof}

With this description in hand, taking hypercohomology in
\eqref{eq:sheafwow} and using Lemma \ref{lem:icohrho}, we obtain
$$ \hbull(\pms_0(r)) \cong \bigoplus_{\rho\vdash r} \left[IH^\bullet(\Sprod_\rho)\otimes L_\rho \right]^{\aut(\rho)}.$$
Applying the K\"unneth formula for intersection cohomology, we arrive at the following decomposition.
\begin{corollary}\label{cor:decomp}
For every $r\geq 1$, there is a natural isomorphism
\begin{equation}\label{eq:maindecomp}
\hbull(\pms_0(r)) \cong \bigoplus_{\rho=[r_1,\ldots,r_k]\vdash r} \left[ IH^\bullet(\ms_0(r_1))\otimes\cdots\otimes IH^\bullet(\ms_0(r_k))\otimes L_\rho \right]^{\aut(\rho)}.
\end{equation}
\end{corollary}

We now describe the representation $L_\rho$ more explicitly by writing the 
partition $\rho$ in terms of its distinct parts and their multiplicities. We 
write 
\[
r=\sum_{a=1}^{{\nparts}}k_an_a,
\;\;\;
n_1>n_2>\cdots>n_{{\nparts}},
\]
where $n_a$ appears in $\rho$ with multiplicity $k_a$. Accordingly, the index set decomposes as
\[
\{1,\ldots,k\}=I_1\sqcup\cdots\sqcup I_{{\nparts}}, \;\;\; 
|I_a|=k_a,
\]
where \(i\in I_a\) if and only if \(r_i=n_a\). Then
$\aut(\rho)
=\Sigma_{k_1}\times\cdots\times\Sigma_{k_{{\nparts}}},$ and  the 
class $w_\rho$ from Theorem
\ref{thm:intform} can be written as
\begin{equation}\label{eq:wmu}
w_\rho=
\prod_{a<b}
\prod_{\substack{i\in I_a\\ j\in I_b}}
(x_i-x_j)^{n_an_b(g-1)}
\prod_{a=1}^{{\nparts}}\;
\prod_{\substack{i,j\in I_a\\ i<j}}
(x_i-x_j)^{n_a^2(g-1)}.
\end{equation}

Now for each $a$, let $\Q_a[-d]$ denote the trivial representation of $\Sigma_{k_a}$ placed in degree $d$, and let $\sgn_a[-d]$ denote the sign representation placed in degree $d$. Set $ d_a=(g-1)n_a^2. $

Note that the factors in \eqref{eq:wmu} joining two distinct blocks $I_a$ and $I_b$ are invariant under $\aut(\rho)$, while the factor internal to $I_a$  transforms under $\Sigma_{k_a}$ by $\sgn_a^{d_a}$.
Moreover, since each $x_i$ has cohomological degree 2, the total degree of $w_\rho$ is
$(g-1)(r^2-\sum_{a=1}^{{\nparts}}k_an_a^2).$
Hence we obtain the following identification.

\begin{lemma}\label{lem:identw}
As a graded one-dimensional representation of $\aut(\rho)$,
\begin{equation}\label{eq:identw}
		\Q w_\rho \cong \Q[(1-g)r^2]\otimes 
		\sgn_1^{d_1}[k_1d_1]\otimes\cdots\otimes 
		\sgn_{{\nparts}}^{d_{{\nparts}}}[k_{{\nparts}}
		 d_{{\nparts}}].
	\end{equation}
\end{lemma}

Combining this with the  identification
$\bigoplus_{j=0}^{n_a-1}\Q x^j \cong H^\bullet(\PP^{n_a-1}),$
we may rewrite the expression in Theorem \ref{thm:corL} as
the graded $\aut(\rho)$-module isomorphism 
\begin{equation}\label{eq:shiftedl}
	L_\rho[(g-1)r^2]
	\cong
	\bigotimes_{a=1}^{{\nparts}}
	\left(
	\sgn_a^{d_a}[d_ak_a]
	\otimes
	H^\bullet(\PP^{n_a-1})^{\otimes k_a}
	\right),
\end{equation}
where, for each \(a\), $\Sigma_{k_a}$ acts by graded permutations of the $k_a$ tensor factors.
Substituting this into \eqref{eq:maindecomp}, and shifting the left-hand side by
\([(g-1)r^2]\), we obtain
\begin{equation}\label{eq:longeq}
H^\bullet(\pms_0(r))[(g-1)r^2]
	\cong
	\bigoplus_{\rho\vdash r}
	\bigotimes_{a=1}^{{\nparts}}
	\left[
	\sgn_a^{d_a}[d_ak_a]
	\otimes
	\left(
	IH^\bullet(\ms_0(n_a))\otimes
	H^\bullet(\PP^{n_a-1})
	\right)^{\otimes k_a}
	\right]^{\Sigma_{k_a}}\!.
\end{equation}
In this formula, $\Sigma_{k_a}$ acts by the graded (Koszul)
permutation action on the $k_a$ copies of the entire factor
$IH^\bullet(\ms_0(n_a))\otimes H^\bullet(\PP^{n_a-1})$.

We now pass from the decomposition for a fixed partition to a
generating-series identity over all ranks by introducing the bigraded vector
space
\begin{equation}\label{eq:defH0}
\mathcal H_0:=\bigoplus_{r\ge1}(\mathcal H_0)_r, \;\;\;	(\mathcal H_0)_r:= IH^\bullet(\ms_0(r))\otimes H^\bullet(\PP^{r-1})[(g-1)r^2],
\end{equation}
where the summand $(\mathcal H_0)_r$ has \textit{rank degree} $r$ and retains its natural \textit{cohomological grading}.

For a graded vector space $\mathcal{H}$,  there is an isomorphism of graded
$\Sigma_{k}$-representations $$\sgn^{d}[d k]\otimes\mathcal{H}^{\otimes {k}}
\cong (\mathcal{H}[d])^{\otimes {k}}.$$
Applying this to $\mathcal H= IH^\bullet(\ms_0(n_a))\otimes 
H^\bullet(\PP^{n_a-1})$, we see that the factor
$ \sgn_a^{d_a}[d_a k_a] $
in \eqref{eq:longeq} converts the corresponding term into the $k_a$-th graded symmetric power of $(\mathcal H_0)_{n_a}$ (cf. \eqref{eq:qtring}). Since summing over all partitions of $r$ gives exactly the degree-$r$
 component of the polynomial algebra on $\mathcal H_0$ (see \eqref{eq:qtring}), we
 arrive at the following result.

\begin{theorem}\label{thm:onemain}
	Let $\mathcal{H}_0$ be the bigraded vector space defined in \eqref{eq:defH0} and set
	\[ \mathcal H_1:=\bigoplus_{r\ge0}(\mathcal H_1)_r, \;\;
	(\mathcal H_1)_0:=\Q \text{\; and \;}
	(\mathcal H_1)_r:=H^\bullet(\pms_0(r))[(g-1)r^2]
	\text{\; for \;}r\ge1.
	\]
	There is  an isomorphism of bigraded vector spaces $\mathcal H_1\cong \Q_{[t]}[\mathcal H_0],$ where $t$ records the cohomological degree.
	Consequently, using the plethystic  convention \eqref{eq:qtring}, we obtain
	\begin{equation*}
		\mathrm{Hilb}_{-t,q}(\mathcal H_1)
		=
		\mathrm{Hilb}_{-t,q}(\Q_{[t]}[\mathcal H_0])
		=
		\mathbb E\mathrm{xp}
		\left[
		\mathrm{Hilb}_{-t,q}(\mathcal H_0)
		\right];
	\end{equation*}  here $q$ records the rank degree.
\end{theorem}

It remains to rewrite this result in the form stated in Theorem \ref{thm:main}. By  Corollary \ref{cor:hecke}, the moduli space $\pms_0(r)$ is isomorphic to a projective bundle over $\ms_1(r)$ with fiber $\PP^{r-1}$. The projective bundle formula therefore gives
$$P_t(\pms_0(r)) = P_t(\ms_1(r)) P_t(\PP^{r-1}),$$
and thus
\[ \mathrm{Hilb}_{-t,q}(\mathcal H_1) = 1+ \sum_{r\geq 1} P_{-t}(\PP^{r-1})\cdot P_{-t}(\ms_1(r))\, (-t)^{(1-g)r^2}q^r. \]
Similarly, from  \eqref{eq:defH0} we have
\[ \mathrm{Hilb}_{-t,q}(\mathcal H_0) = \sum_{r\geq 1} IP_{-t}(\ms_0(r))\cdot P_{-t}(\PP^{r-1})\,
(-t)^{(1-g)r^2}q^r. \]
Theorem \ref{thm:onemain} therefore yields
\begin{multline}
	1+
	\sum_{r\geq1}
	P_{-t}(\PP^{r-1})P_{-t}(\ms_1(r))
	(-t)^{(1-g)r^2}q^r
	\\
	=
	\mathbb E\mathrm{xp}
	\left[
	\sum_{r\geq 1}
	P_{-t}(\PP^{r-1})IP_{-t}(\ms_0(r))
	(-t)^{(1-g)r^2}q^r
	\right],
\end{multline}
which completes the proof of Theorem \ref{thm:main}.

\begin{remark} 
	The same argument may be refined to Hodge polynomials, since all
	constructions
	above are compatible with mixed Hodge structures. For a smooth projective
	variety \(Y\), set
	\[
	h_{u,v}(Y)=\sum_{p,q}(-1)^{p+q}\dim H^{p,q}(Y)u^pv^q,
	\]
and define \(Ih_{u,v}(Y)\) analogously using intersection cohomology.
Working in the coefficient ring obtained by adjoining a formal
variable $z$ with $z^2=uv$, we obtain
	\begin{multline}\label{eq:hodgeref}
		1+
		\sum_{r=1}^{\infty}
		h_{u,v}(\PP^{r-1})h_{u,v}(\ms_1(r))
		(-z)^{(1-g)r^2}q^r
		\\
		=
		\mathbb E\mathrm{xp}
		\left[
		\sum_{r=1}^{\infty}
		h_{u,v}(\PP^{r-1})Ih_{u,v}(\ms_0(r))
		(-z)^{(1-g)r^2}q^r
		\right].
	\end{multline}
\end{remark}

\appendix
\section{Supplementary proofs}\label{sec:supplementary} 
\subsection{Proof of Theorem \ref{thm:fiberring}}\label{sec:fiberring}
We now show that for $V\in\modtor$, the cohomology ring of the fiber $F_V$ is 
isomorphic to the quotient ring $\mathbb{Q}[x_1,\ldots,x_k]/\Idm$, where $\Idm$ is  
the ideal generated by  \eqref{eq:HStrelationsFib}. 

We first verify that the relations \eqref{eq:HStrelationsFib} hold in 
$H^\bullet(F_V)$. 
\begin{lemma}\label{lem:phisatisfied}
Let $x_1,\ldots,x_k$ be the cohomology classes defined in Proposition 
\ref{prop:xsigma}. For every decomposition
$A\sqcup B$ of $\{1,\ldots,k\}$ with $ A\neq\emptyset$,
the class $\phi_{A\sqcup B}(x_1,\ldots,x_k)$ vanishes in $H^\bullet(F_V)$.
\end{lemma}
\begin{proof}
It follows from the proof of Proposition \ref{prop:xsigma} that the restriction map 
		\begin{equation}\label{eq:fibinject}
			\hbull(F_V)\to \bigoplus_{\sigma\in\Sigma_k}\hbull(\pat_\sigma)
		\end{equation}
		is injective. Hence   it is enough to prove that $\phi_{A\sqcup B}$ restricts to zero in the cohomology ring of the tower $\pat_\sigma$ for every permutation $\sigma\in\Sigma_k$. 

Fix $\sigma\in\Sigma_k$; it follows from  Propositions \ref{prop:ringtower} and \ref{prop:xsigma} that  the cohomology ring  $H^\bullet(\pat_\sigma)$ is generated by the restrictions of the classes $x_1,\ldots,x_k$, with relations
\begin{multline}\label{eq:phisigmatowers}
\phi_1^\sigma(x_1,\ldots,x_k):= x_{\sigma(1)}^{r_{\sigma(1)}}, \\ \phi_i^\sigma(x_1,\ldots,x_k):=x_{\sigma(i)}^{r_{\sigma(i)}}\cdot\prod_{j=1}^{i-1}(x_{\sigma(i)}-x_{\sigma(j)})^{r_{\sigma(i)}r_{\sigma(j)}(g-1)},\;\; i=2,\ldots,k.
\end{multline}

Let $a:=\min\{i:\sigma(i)\in A\};$ we then have $ \sigma(1),\ldots,\sigma(a-1)\in B,$ and thus  the polynomial 
$$\phi^\sigma_a(x_1,\ldots,x_k)= x_{\sigma(a)}^{r_{\sigma(a)}} \prod_{b<a} \bigl(x_{\sigma(a)}-x_{\sigma(b)}\bigr)^{r_{\sigma(a)}r_{\sigma(b)}(g-1)}$$ divides $\phi_{A\sqcup B}$. Hence the restriction of $\phi_{A\sqcup B}$  to $\pat_\sigma$ vanishes. Since $\sigma\in\Sigma_k$ was arbitrary, the result follows.
\end{proof}
For each $\sigma\in\Sigma_k$, denote by $I_\sigma\subset \mathbb{Q}[x_1,\ldots,x_k]$ the ideal generated by the polynomials \eqref{eq:phisigmatowers}. Set
\[
I_\Sigma:=\bigcap_{\sigma\in\Sigma_k}I_\sigma, \;\; R:=\mathbb{Q}[x_1,\ldots,x_k]/\Idm, \;\; S:=\mathbb{Q}[x_1,\ldots,x_k]/I_\Sigma.
\]
Lemma \ref{lem:phisatisfied} provides a natural surjective homomorphism
$\eta\colon R\twoheadrightarrow S.$
Moreover, it follows from  \eqref{eq:fibinject} that there is an injective homomorphism
$\theta\colon S\hookrightarrow H^\bullet(F_V),$ whose image is the subring generated by the classes $x_1,\ldots,x_k$.

Our goal is to prove that the composition
\begin{equation}\label{eq:compositionlem}
R\overset{\eta}{\twoheadrightarrow} S\overset{\theta}{\hookrightarrow}H^\bullet(F_V)
\end{equation}
is an isomorphism. We start by comparing the Hilbert series of both sides.
\begin{proposition}\label{lem:hilb-agree}
Assign  degree $2$ to each $x_i\in R$. Then $R$ is a graded ring whose Hilbert 
series coincides with the Poincar\'e polynomial of $F_V$.
\end{proposition}
\begin{proof}
First note that the polynomial $g(\rho;t)$ from \eqref{eq:gpoly} depends only on 
the 
graph $\mathbb{G}_\rho$, $g(\rho;t)=\hat{g}(\mathbb{G}_\rho,t)$, 
and one can define $\hat{g}(\Gamma,t)$ for any other directed graph on a vertex 
set $V(\Gamma)$ and with $r_{ij}=r_{ji}$ edges pointing from vertex $i$ to $j$ 
and $j$ to $i$.

The key point is an observation of  \cite[Theorem 6.3]{ZhP}: 
			the function $\hat{g}(\Gamma,t^2)$ 
			coincides with 
			the Hilbert function of the quotient algebra
	\begin{equation}\label{eq:quotalg}
			\mathbb{Q}[x_i-x_j,\, 0\leq i<j\leq k]/	
			\big{\langle} {\prod_{\substack{i\in 
						A, j\in B}}(x_i-x_j)^{r_{ij}}, \; A\sqcup 
					B\vdash V(\Gamma), \, A,B\neq\emptyset}\big{\rangle},
			\end{equation} where $A\sqcup 
			B$ runs over all nontrivial decompositions of the vertex set 
			$V(\Gamma)$. 

Now we note that the polynomial $f(\rho;t)$ from \eqref{eq:frec-graph} coincides with $\hat g(\augf,t)$, where $\augf$ is the graph introduced in \S\ref{sec:defFrho}.
It follows that $f(\rho;t^2)$ equals the Hilbert function of \eqref{eq:quotalg} for the data
$r_{ij}=r_{ji}=r_ir_j(g-1)$ for $1\leq i<j\le k$, and $r_{0i}=r_{i0}=r_i $ for
$i=1,\dots,k$.

Then noting that the  change of variables $x_i-x_0\mapsto x_i$, $i=1,\ldots,k$ 
converts the algebra 
\eqref{eq:quotalg} to the algebra \eqref{eq:HStrelationsFib}, we complete the proof of Proposition \ref{lem:hilb-agree}.
\end{proof}

The proof of Theorem \ref{thm:fiberring} now reduces to showing that $\eta$ is injective.

\begin{lemma}\label{lem:eta-injective}
The homomorphism $\eta$ is injective.
\end{lemma}
\begin{proof}
Let $N=\sum_{1\leq i<j\leq k}r_ir_j(g-1)+r-k =\dim F_V$. We start with proving 
that $\eta$ induces an isomorphism $R^{2N}\!\cong S^{2N}$ on top degree.  For 
each $\sigma\in\Sigma_k$, we define 
\[
\Psi_\sigma:= \prod_{i=1}^k x_{\sigma(i)}^{r_{\sigma(i)}-1+\sum_{j<i}{r_{\sigma(i)}r_{\sigma(j)}(g-1)}}\in H^{2N}(F_V).
\]
Using Corollary \ref{cor:towerintegral}, one can show that 
\[
\int_{\pat_\tau}\Psi_\sigma\big|_{\pat_\tau}=\begin{cases}
1,&\text{if }\tau=\sigma,\\
0,&\text{otherwise},
\end{cases}
\]
and thus the classes $\{\Psi_\sigma\}_{\sigma\in\Sigma_k}$ form a basis of
\[
H^{2N}(F_V) \cong\bigoplus_{\sigma\in\Sigma_k}H^{2N}(\pat_\sigma).
\]
In particular, we have $\dim H^{2N}(F_V)=k!.$
Since each monomial $\Psi_\sigma$ belongs to the image of $\theta$, the injective map $\theta$ induces an isomorphism  $S^{2N}\to H^{2N}(F_V)$, and consequently $\dim S^{2N}=k!$. 
On the other hand, Proposition \ref{lem:hilb-agree} gives 
$\dim R^{2N}= \dim H^{2N}(F_V)= k!$. Since $\eta\colon R^{2N}\to S^{2N}$ is surjective, it is therefore an isomorphism.

Next, we observe that under the change of variables used in the proof of Proposition \ref{lem:hilb-agree}, \(R\) can be identified with the cohomology ring of the core of the Lawrence toric variety associated to the graph $\mathbb{F}_\rho$ introduced in \S\ref{sec:defFrho} (see \cite[Theorem 6.3]{ZhP}). Then applying the levelness argument from the proof of  \cite[Theorem 7.1]{HSt}, we deduce that  $\operatorname{Soc}(R)=R^{2N}.$

Now suppose that $\ker(\eta)\neq 0$, and choose a nonzero homogeneous element $h\in\ker(\eta)$
of maximal degree. Since $\ker(\eta)$ is an ideal, we have $x_i h\in\ker(\eta)$  for every $i$. The  element $x_i h$ has degree strictly larger than that of $h$, so the maximality of $h$ implies that $x_ih=0$ for every $i$.
It follows that $h\in\operatorname{Soc}(R)=R^{2N}.$
However, $\eta$ is injective in degree $2N$, hence $h=0$, a contradiction. Therefore $\ker(\eta)=0$, and $\eta$ is injective.
\end{proof}
The composition \eqref{eq:compositionlem} is thus an isomorphism; this finishes the proof of Theorem \ref{thm:fiberring}.

\subsection{The isomorphism for non-abelian fibers} \label{sec:suppisom}
Recall that in \S\ref{sec:injcohmap} we constructed an
injective map from the cohomology of a non-abelian fiber to the monodromy
invariants in the cohomology of a nearby abelian fiber.  For completeness, we 
now prove
the stronger statement that this map is an isomorphism, for which it is 
sufficient to
compare the Hilbert series of its source and target.

As in \S\ref{sec:injcohmap}, to simplify the notation, we describe in detail the case $\wV=U^{\oplus k}\!,$ where \(U\) is a stable bundle of rank $\br$ and degree zero.
\begin{proposition}\label{prop:Enonab}
In the notation of \S\ref{sec:injcohmap},
\begin{equation}\label{eq:inveq}
\mathrm{Hilb}_t\bigl(\hbull(F_{U^{\oplus k}})\bigr) = \mathrm{Hilb}_t\bigl(\hbull(F_V)^{\Sigma_k}\bigr).
\end{equation}
Consequently, the map \eqref{eq:gmap} is an isomorphism.
\end{proposition}

\begin{proof}
To analyze the right-hand side of \eqref{eq:inveq}, we return to the Mayer--Vietoris exact sequence \eqref{eq:MV}. For
$J_a=\{1,\ldots,a\}$, we have
\[\bigoplus_{|J|=a} e_J\otimes 
H^\bullet(\patf_J)\simeq\mathrm{Ind}_{\Sigma_a\times\Sigma_{k-a}}^{\Sigma_k}\bigl(
 e_{J_a}\otimes H^\bullet(\patf_{J_a})\bigr).\]
Taking $\Sigma_k$-invariants and applying Frobenius reciprocity gives
$$\Bigl[ \bigoplus_{|J|=a} e_J\otimes H^\bullet(\patf_J)\Bigr]^{\Sigma_k}
\!\!\simeq\Bigl[e_{J_a}\otimes 
H^\bullet(\patf_{J_a})\Bigr]^{\Sigma_a\times\Sigma_{k-a}}.$$

The  description \eqref{eq:FJproj} of \(\patf_J\)  shows that \(\patf_{J_a}\) is a fiber product
over \(F_{V_{a+1}\oplus\cdots\oplus V_k}\) of \(a\) projective bundles,
each with fiber \(\mathbb P^{N_a-1}\), where $N_a=\br^2(k-a)(g-1)+\br$. Note that 
the factor \(e_{J_a}\) carries the
sign representation of \(\Sigma_a\), hence applying the inductive hypothesis to
the \(\Sigma_{k-a}\)-action and then taking \(\Sigma_a\)-invariants 
gives
\[
\Bigl[
\bigoplus_{|J|=a} e_J\otimes H^\bullet(\patf_J)
\Bigr]^{\Sigma_k}
\simeq
H^\bullet(F_{U^{\oplus(k-a)}})
\otimes
\bigwedge^a H^\bullet(\mathbb P^{N_a-1}).
\]
Taking the Euler characteristic of the invariant Mayer--Vietoris complex, and
using
$$\mathrm{Hilb}_t\Bigl(\bigwedge^a H^\bullet(\mathbb P^{N_a-1})
\Bigr) =  t^{a(a-1)}E_t\bigl(\mathbb G(a,N_a)\bigr),$$
where \(\mathbb G(a,N)\) is the Grassmannian of
\(a\)-dimensional subspaces of \(\mathbb C^N\), we find
\begin{equation}\label{eq:finetcompressed}
		\mathrm{Hilb}_t\left(\hbull(F_V)^{\Sigma_k}\right)	=	\sum_{a=1}^k	(-1)^{a+1}
		t^{a(a-1)} E_t\bigl(F_{U^{\oplus(k-a)}}\bigr)E_t\bigl(\mathbb 
		G(a,N_a)\bigr).
\end{equation}

We now turn to the non-abelian fiber, the left-hand side of
\eqref{eq:inveq}. Define the following generalization of the parabolic extension space \(\ParExtSp\):
	\begin{equation}\label{eq:parext-space}
		\ParExtSp^a := \{0\to U^{\oplus a}\to W\to W'\to 0 \mid
		W'\in F_{U^{\oplus(k-a)}},\ W\in\pms_0(k\br) \}/\!\sim\;,
	\end{equation}
where extensions are considered up to natural isomorphisms.

For fixed \(W'\in F_{U^{\oplus(k-a)}}\), an extension in
\eqref{eq:parext-space} is represented by an element of
\[
\operatorname{ParExt}^1(W',U^{\oplus a})
\cong \operatorname{ParExt}^1(W',U)\otimes\mathbb C^a
\cong\operatorname{Hom}\bigl((\mathbb C^a)^*\!, \,
\operatorname{ParExt}^1(W',U)\bigr).
\]
The parabolic stability condition is equivalent to the
corresponding linear map having rank $a$. Modulo
$\Aut(U^{\oplus a})\cong \GL_a$, such maps are therefore parametrized by
the Grassmannian of \(a\)-dimensional subspaces of \(\operatorname{ParExt}^1(W',U)\). 
As $W'$ varies, the spaces $\operatorname{ParExt}^1(W',U)$ form a vector bundle of rank $N_a$ over $F_{U^{\oplus(k-a)}}$. Thus $\ParExtSp^a$ is the
corresponding relative Grassmannian, and 
\begin{equation}\label{eq:pattsa-base}
E_t(\ParExtSp^a)=E_t\bigl(F_{U^{\oplus(k-a)}}\bigr)E_t\bigl(\mathbb 
G(a,N_a)\bigr).
\end{equation}

On the other hand, consider the decomposition
\[F_{U^{\oplus k}}=\Fiso_1\sqcup\cdots\sqcup\Fiso_k,\]
where $\Fiso_m$ is the locus of bundles $W$ whose maximal subbundle isomorphic 
to a sum of copies of $U$ is  $U^{\oplus m}$.
For $W\in\Fiso_m$, write this maximal subbundle as $U\otimes M_W,$ $\dim M_W=m.$ A point of $\ParExtSp^a$ lying over $W$ is then equivalent to the choice of
an $a$-dimensional subspace of $M_W$. Hence, over $\Fiso_m$, the map
$\ParExtSp^a\to F_{U^{\oplus k}}$ is a Zariski-locally trivial
Grassmannian bundle with fiber $\mathbb G(a,m)$, and hence
\begin{equation}\label{eq:pattsa-strata}
E_t(\ParExtSp^a)=\sum_{m\ge a}
E_t(\Fiso_m)E_t\bigl(\mathbb G(a,m)\bigr).
\end{equation}

Combining \eqref{eq:pattsa-base} and
\eqref{eq:pattsa-strata}, and substituting the resulting identity into
\eqref{eq:finetcompressed}, gives
\[	\mathrm{Hilb}_t\left(\hbull(F_V)^{\Sigma_k}\right)=\sum_{m=1}^k 	
E_t(\Fiso_m)\sum_{a=1}^m (-1)^{a+1} t^{a(a-1)}
E_t\bigl(\mathbb G(a,m)\bigr).\]
The $q$-binomial identity, with $q=t^2$, reads
$$\sum_{a=1}^m (-1)^{a+1}t^{a(a-1)}
E_t\bigl(\mathbb G(a,m)\bigr)=1,$$
therefore 
\[\mathrm{Hilb}_t\left(\hbull(F_V)^{\Sigma_k}\right)=\sum_{m=1}^kE_t(\Fiso_m)=
E_t\bigl(F_{U^{\oplus k}}\bigr).\]

Since \(F_{U^{\oplus k}}\) has  pure cohomology concentrated in even degrees, 
$E_t\bigl(F_{U^{\oplus k}}\bigr)=\mathrm{Hilb}_t\bigl(\hbull(F_{U^{\oplus 
k}})\bigr),$
which  completes the proof of Proposition \ref{prop:Enonab}.
\end{proof}
The same Hilbert-series argument applies successively to the
repeated blocks of an arbitrary multipartition, exactly as in the proof of
the general case in \S\ref{sec:injcohmap}. This completes the comparison for all
non-abelian fibers.


\begin{thebibliography}{10}

\bibitem{AB83}
M.~F. Atiyah and R.~Bott.
\newblock The {Y}ang-{M}ills equations over {R}iemann surfaces.
\newblock {\em Philos. Trans. Roy. Soc. London Ser. A}, 308(1505):523--615,
  1983.

\bibitem{BBD82}
A.~A. Be\u{\i}linson, J.~Bernstein, and P.~Deligne.
\newblock Faisceaux pervers.
\newblock In {\em Analysis and topology on singular spaces, {I} ({L}uminy,
  1981)}, volume 100 of {\em Ast\'{e}risque}, pages 5--171. Soc. Math. France,
  Paris, 1982.

\bibitem{BGL94}
E.~Bifet, F.~Ghione, and M.~Letizia.
\newblock On the {A}bel-{J}acobi map for divisors of higher rank on a curve.
\newblock {\em Math. Ann.}, 299(4):641--672, 1994.

\bibitem{BodenH}
H.~U. Boden and Y.~Hu.
\newblock Variations of moduli of parabolic bundles.
\newblock {\em Math. Ann.}, 301(3):539--559, 1995.

\bibitem{BDIKP}
C.~Bu, B.~Davison, A.~Ib\'{a}\~{n}ez N\'{u}\~{n}ez, T.~Kinjo, and
  T.~P\u{a}durariu.
\newblock Cohomology of symmetric stacks.
\newblock {\em arXiv:2502.04254}, 2025.

\bibitem{dCM05}
M.~A.~A. de~Cataldo and L.~Migliorini.
\newblock The {H}odge theory of algebraic maps.
\newblock {\em Ann. Sci. \'{E}cole Norm. Sup. (4)}, 38(5):693--750, 2005.

\bibitem{dCM09}
M.~A.~A. de~Cataldo and L.~Migliorini.
\newblock The decomposition theorem, perverse sheaves and the topology of
  algebraic maps.
\newblock {\em Bull. Amer. Math. Soc. (N.S.)}, 46(4):535--633, 2009.

\bibitem{D71}
P.~Deligne.
\newblock Th\'{e}orie de {H}odge. {II}.
\newblock {\em Inst. Hautes \'{E}tudes Sci. Publ. Math.}, 40:5--57, 1971.

\bibitem{DR75}
U.~V. Desale and S.~Ramanan.
\newblock Poincar\'{e} polynomials of the variety of stable bundles.
\newblock {\em Math. Ann.}, 216(3):233--244, 1975.

\bibitem{D95}
A.~H. Durfee.
\newblock Intersection homology {B}etti numbers.
\newblock {\em Proc. Amer. Math. Soc.}, 123(4):989--993, 1995.

\bibitem{EK00}
R.~Earl and F.~Kirwan.
\newblock The {H}odge numbers of the moduli spaces of vector bundles over a
  {R}iemann surface.
\newblock {\em Q. J. Math.}, 51(4):465--483, 2000.

\bibitem{EH16}
D.~Eisenbud and J.~Harris.
\newblock {\em 3264 and all that---a second course in algebraic geometry}.
\newblock Cambridge University Press, Cambridge, 2016.

\bibitem{FT}
C.~Felisetti and O.~Trapeznikova.
\newblock Intersection theory on singular moduli spaces of vector bundles: a
  parabolic approach.
\newblock {\em arXiv:2603.01946}, 2026.

\bibitem{GM83}
M.~Goresky and R.~MacPherson.
\newblock Intersection homology. {II}.
\newblock {\em Invent. Math.}, 72(1):77--129, 1983.

\bibitem{GS75}
P.~Griffiths and W.~Schmid.
\newblock Recent developments in {H}odge theory: a discussion of techniques and
  results.
\newblock In {\em Discrete subgroups of {L}ie groups and applications to moduli
  ({I}nternat. {C}olloq., {B}ombay, 1973)}, Tata Inst. Fundam. Res. Stud.
  Math., No. 7, pages 31--127. Published for the Tata Institute of Fundamental
  Research, Bombay by Oxford University Press, Bombay, 1975.

\bibitem{HN74}
G.~Harder and M.~S. Narasimhan.
\newblock On the cohomology groups of moduli spaces of vector bundles on
  curves.
\newblock {\em Math. Ann.}, 212:215--248, 1974/75.

\bibitem{HSt}
T.~Hausel and B.~Sturmfels.
\newblock Toric hyper{K}\"{a}hler varieties.
\newblock {\em Doc. Math.}, 7:495--534, 2002.

\bibitem{JKKW}
L.~Jeffrey, Y.H.~Kiem, F.~Kirwan, and J.~Woolf.
\newblock Intersection pairings on singular moduli spaces of bundles over a {R}iemann surface and their partial desingularisations.
\newblock {\em Transform. Groups}, 11(3):439--494, 2006.

\bibitem{K06}
Y.H. Kiem.
\newblock Intersection cohomology of representation spaces of surface groups.
\newblock {\em Internat. J. Math.}, 17(2):169--182, 2006.

\bibitem{K85}
F.~Kirwan.
\newblock Partial desingularisations of quotients of nonsingular varieties and
  their {B}etti numbers.
\newblock {\em Ann. of Math. (2)}, 122(1):41--85, 1985.

\bibitem{K86vect}
F.~Kirwan.
\newblock On the homology of compactifications of moduli spaces of vector
  bundles over a {R}iemann surface.
\newblock {\em Proc. London Math. Soc. (3)}, 53(2):237--266, 1986.

\bibitem{K86}
F.~Kirwan.
\newblock Rational intersection cohomology of quotient varieties.
\newblock {\em Invent. Math.}, 86(3):471--505, 1986.

\bibitem{L83}
H.~Lange.
\newblock Universal families of extensions.
\newblock {\em J. Algebra}, 83(1):101--112, 1983.

\bibitem{L96}
Y.~Laszlo.
\newblock Local structure of the moduli space of vector bundles over curves.
\newblock {\em Comment. Math. Helv.}, 71(3):373--401, 1996.

\bibitem{M23}
M.~Mauri.
\newblock Intersection cohomology of rank 2 character varieties of surface
  groups.
\newblock {\em J. Inst. Math. Jussieu}, 22(4):1615--1654, 2023.

\bibitem{MehtaSeshadri}
V.~B. Mehta and C.~S. Seshadri.
\newblock Moduli of vector bundles on curves with parabolic structures.
\newblock {\em Math. Ann.}, 248(3):205--239, 1980.

\bibitem{Meinhardt}
S.~Meinhardt.
\newblock {D}onaldson--{T}homas invariants vs. intersection cohomology for
  categories of homological dimension one.
\newblock {\em arXiv:1512.03343}, 2015.

\bibitem{MR15}
S.~Mozgovoy and M.~Reineke.
\newblock Intersection cohomology of moduli spaces of vector bundles over
  curves.
\newblock {\em Moduli}, 3:Paper No. e9, 2026.

\bibitem{MF82}
D.~Mumford, J.~Fogarty, and F.~Kirwan.
\newblock {\em Geometric invariant theory}, volume~34 of {\em Ergebnisse der
  Mathematik und ihrer Grenzgebiete (2) [Results in Mathematics and Related
  Areas (2)]}.
\newblock Springer-Verlag, Berlin, third edition, 1994.

\bibitem{NS65}
M.~S. Narasimhan and C.~S. Seshadri.
\newblock Stable and unitary vector bundles on a compact {R}iemann surface.
\newblock {\em Ann. of Math. (2)}, 82:540--567, 1965.

\bibitem{Newstead67}
P.~E. Newstead.
\newblock Topological properties of some spaces of stable bundles.
\newblock {\em Topology}, 6(2):241--262, 1967.

\bibitem{N09}
N.~Nitsure.
\newblock Deformation theory for vector bundles.
\newblock In {\em Moduli spaces and vector bundles}, volume 359 of {\em London
  Math. Soc. Lecture Note Ser.}, pages 128--164. Cambridge Univ. Press,
  Cambridge, 2009.

\bibitem{PS08}
C.~A.~M. Peters and J.~H.~M. Steenbrink.
\newblock {\em Mixed {H}odge structures}, volume~52 of {\em Ergebnisse der
  Mathematik und ihrer Grenzgebiete. 3. Folge. A Series of Modern Surveys in
  Mathematics [Results in Mathematics and Related Areas. 3rd Series. A Series
  of Modern Surveys in Mathematics]}.
\newblock Springer-Verlag, Berlin, 2008.

\bibitem{S89}
M.~Saito.
\newblock Introduction to mixed {H}odge modules.
\newblock {\em Ast\'{e}risque}, 179-180:145--162, 1989.
\newblock Actes du Colloque de Th\'{e}orie de Hodge (Luminy, 1987).

\bibitem{OTASz}
A.~Szenes and O.~Trapeznikova.
\newblock The parabolic {V}erlinde formula: iterated residues and
  wall-crossings.
\newblock {\em Geom. Topol.}, 28(5):2259--2311, 2024.

\bibitem{ZhP}
A.~Szenes and O.~Trapeznikova.
\newblock Intersection cohomology of type-{A} toric varieties.
\newblock {\em Algebraic Combinatorics}, 8(2):575--596, 2025.

\bibitem{W17}
G.~Williamson.
\newblock The {H}odge theory of the decomposition theorem.
\newblock {\em Ast\'{e}risque}, pages Exp. No. 1115, 335--367, 2017.
\newblock S\'{e}minaire Bourbaki. Vol. 2015/2016. Expos\'{e}s 1104--1119.
\end{thebibliography}
\end{document}